\documentclass[11 pt, reqno]{article}

\usepackage{amsmath,amsfonts,amssymb,amsthm}
\usepackage[colorlinks=true,hyperindex=true]{hyperref}
\usepackage{cancel,bbm}

\usepackage{mathabx} 
\usepackage{comment}
\usepackage{enumitem}
\usepackage{cite}

\usepackage{mathrsfs}

\usepackage{chngcntr}
\counterwithin*{equation}{section}

\usepackage[left=3cm, right=3cm, top=3 cm, bottom= 3 cm]{geometry}
\usepackage{color}
\usepackage{fancyhdr}
\usepackage{latexsym}

\usepackage{bm}

\newtheorem{ccounter}{ccounter}[section]
\newtheorem{thm}[ccounter]{Theorem}
\newtheorem{lem}[ccounter]{Lemma}
\newtheorem{cor}[ccounter]{Corollary}
\newtheorem{defn}[ccounter]{Definition}
\newtheorem{prop}[ccounter]{Proposition}
\newtheorem{ass}[ccounter]{Assumption}
\newtheorem{ex}[ccounter]{Example}

\newcommand{\pare}[1]{\left({#1}\right)}
\newcommand{\sqbrac}[1]{\left[ {#1}\right]}
\newcommand{\brac}[1]{\left\{{#1} \right\}}
\newcommand{\abs}[1]{\left|{#1}\right|}
\DeclareMathOperator{\Av}{Av}
\DeclareMathOperator{\Flat}{Flat}

\def\bet{\begin{thm}}
\def\eet{\end{thm}}
\def\bel{\begin{lem}}
\def\eel{\end{lem}}
\def\bas{\begin{ass}}
\def\eas{\end{ass}}
\def\bec{\begin{cor}}
\def\eec{\end{cor}}
\def\bed{\begin{defn}}
\def\eed{\end{defn}}
\def\bep{\begin{prop}}
\def\eep{\end{prop}}
\def\beq{\begin{equation}}
\def\eeq{\end{equation}}
\def\proof{\noindent {\bf Proof.}\ \ }
\def\bea{\begin{equation*}}
\def\eea{\end{equation*}}

\def\bex{\begin{ex}}
\def\eex{\end{ex}}
\def\remark{\noindent{\bf Remark. }}

\def\rr{\mathbb{R}}
\def\zz{\mathbb{Z}}
\def\cc{\mathbb{C}}
\def\1{\boldsymbol{1}}
\let\Im\relax
\DeclareMathOperator{\Im}{Im}
\def\im{\Im}
\let\Re\relax
\DeclareMathOperator{\Re}{Re}
\def\e{\mathrm{e}}
\def\i{\mathrm{i}}
\def\del{\partial}
\def\d{\mathrm{d}}
\def\eps{\varepsilon}
\renewcommand\leq\varleq
\renewcommand\geq\vargeq
\def\ee{\mathrm{E}}

\def\F{\mathcal{F}}
\def\O{\mathcal{O}}

\def\ee{\mathbb{E}}

\def\pp{\mathbb{P}}

\def\msc{m_{\mathrm{sc}}}
\def\rhosc{\rho_{\mathrm{sc}}}
\def\D{\mathcal{D}}
\def\I{\mathcal{I}}

\def\A{\mathcal{A}}

\def\dto{\downarrow}

\def\D{\mathcal{D}}

\def\be{\mathbf{e}}

\def\E{\mathcal{E}}

\def\mfc{\mathfrak{c}}

\def\gamkt{\gamma_k^t}
\def\bu{\bar{\mathfrak{u}}}
\def\sS{\mathscr{S}}
\def\US{\mathcal{U}^{\mathscr{S}}}
\def\pP{\mathscr{P}}
\def\bB{\mathscr{B}}
\def\sB{\mathscr{B}}
\def\UB{U^{\mathscr{B}}}
\def\lL{\mathscr{L}}
\def\be{\boldsymbol{e}}

\def\xnu{x^{(\nu)}}
\def\fu{\mathfrak{u}}
\def\hati{\hat{i}}
\def\hh{\mathbb{H}}
\def\Da{\mathcal{I}_{\mathrm{ad}}}
\def\Db{\mathcal{I}_{\mathrm{bulk}}}
\def\De{\mathcal{I}_{\mathrm{edge}}}
\def\cb{\mathfrak{c}_b}
\def\ce{\mathfrak{c}_e}
\def\blam{\boldsymbol{\lambda}}
\def\hblam{\boldsymbol{\hat{\lambda}}}

\def\hA{\hat{\mathcal{A}}}

\begin{document}

\begin{table}
\centering

\begin{tabular}{c}
\multicolumn{1}{c}{\Large{\bf Edge homogenization of Dyson Brownian motion}}\\ 
\multicolumn{1}{c}{\Large{\bf 
\phantom{blah}}}\\ 
\multicolumn{1}{c}{\Large{\bf and applications}} \\
\\
\\
\end{tabular}
\begin{tabular}{ c c c  }
Benjamin Landon
& \phantom{blah} & 
Tianhao Xian
 \\
 & & \\  
  \small{University of Toronto}&  &  \small{University of Toronto} \\
  \small{Department of Mathematics}&  & \small{Department of Mathematics}\\
 \small{\texttt{blandon@math.toronto.edu}} & & \small{\texttt{tianhao.xian@utoronto.ca}} \\
  & & \\
\end{tabular}
\\
\begin{tabular}{c}
\multicolumn{1}{c}{\today}\\
\\
\end{tabular}

\begin{tabular}{p{15 cm}}
\small{{\bf Abstract:} We prove a homogenization result for the difference of two coupled Dyson Brownian motions started from generalized Wigner matrix initial data. We prove an optimal order, high probability estimate that is valid throughout the spectrum, including up to the spectral edges. Prior homogenization results concerned only the bulk of the spectrum.

\phantom{blah} We apply our estimate to address the question of quantifying edge universality. Here, we have two results. We show that the Kolmogorov-Smirnov distance of the distribution of the gap between the largest two eigenvalues of a generalized Wigner matrix (with smooth entry distribution) and its GOE/GUE counterpart is $\O(N^{-1+\eps})$. On the other hand, we show that, for the distribution of the largest eigenvalue, there are Wigner matrices so that the analogous Kolmogorov-Smirnov distance is bounded below by $N^{-1/3-\eps}$. }
\end{tabular}
\end{table}

\tableofcontents

\section{Introduction}

Dyson Brownian motion (DBM) has proven an effective tool for studying the spectral universality of large random matrices \cite{local-relax}. Perhaps the most powerful approach to handling DBM is that of homogenization: here, the difference between two coupled DBMs is computed to leading order by developing a homogenization theory for the long-range parabolic equation that the difference satisfies. In particular, the strongest universality results are proved in part by using homogenization. This approach to DBM was initiated in the work \cite{fixed-wig}, 
which was applied to prove a strong form of eigenvalue universality for Wigner matrices (universality of the local correlation functions without any local averaging). Stronger and more general homogenization results were proved in \cite{fixed,bourgade2018extreme}. In particular, in the remarkable work \cite{bourgade2018extreme}, Bourgade developed a new approach to homogenization, allowing for optimal estimates in the bulk of the spectrum. As a consequence, the work \cite{bourgade2018extreme} obtained the universality of the minimal gap between consecutive eigenvalues for generalized Wigner matrices.

Despite these advances,  all of these  works nonetheless fall short of developing a homogenization theory at the spectral edges or any intermediate regimes. The main result of this paper is to - for the first time - develop a homogenization theory that holds throughout the entire spectrum, with optimal estimates. Our results apply to the case of two coupled DBMs started from generalized Wigner matrix initial data.


While interesting in its own right, the main purpose of developing the homogenization theory is for applications to eigenvalue universality. Let us now highlight some applications to an open problem. 
We will then return to the above discussion of homogenization in greater detail. In order to place the problem we have in mind in context, let us recall one form of the classical Berry-Esseen CLT: if $\{ X_i\}_i$ are i.i.d. centered, variance one random variables with finite third moment, then
\beq \label{eqn:be-CLT}
\sup_{ r \in \rr} \left| \pp\left[ \frac{1}{\sqrt{n}} \sum_{i=1}^n X_i \geq r \right] - \pp\left[ Z \geq r \right] \right| \leq \frac{C}{\sqrt{n}} 
\eeq
where $Z$ is a standard normal random variable. The quantity on the LHS is the Kolmogorov-Smirnov (KS) distance between the two probability distributions. 

Two important points about the Berry-Esseen CLT are:
\begin{enumerate}[label=(\roman*)]
\item Universality: the limiting Gaussian distribution is the same for every choice of the distribution of the $\{ X_i\}_i$
\item Rate of convergence: the error in terms of the KS distance is the square-root of the number of variables in the sum
\end{enumerate}

Random matrix statistics are widely believed to represent a new universality paradigm. They arise as the limiting distributions of large, complex, correlated systems (see, e.g., \cite{mehta,erdos2017dynamical} and the discussions therein), whereas the classical CLT and Poisson laws characterize systems with a  high degree of independence. 

One example of this phenomena is the \emph{edge universality of Wigner matrices}. Real symmetric and complex Hermitian Wigner matrices $W$ are $N \times N$ self-adjoint matrices whose upper triangular part consists of i.i.d. centered real or complex random variables with  variance $N^{-1}$.  If $\lambda_1$ is the largest eigenvalue of a Wigner matrix $W$ then,
\beq \label{eqn:TW}
\lim_{N \to \infty} N^{2/3} ( \lambda_1 - 2 ) \stackrel{d}{=} \mathrm{TW}_\beta
\eeq
where $\mathrm{TW}_\beta$ denotes a Tracy-Widom random variable of the same symmetry type as $W$ ($\beta=1$ for real symmetric, $\beta=2$ for complex Hermitian). The Tracy-Widom distributions are a new class of distributions first computed for random matrices in the case that the Wigner matrix consists of Gaussian random variables (the Gaussian Orthogonal and Unitary Ensembles, respectively (GOE/GUE)) \cite{tracy1994level,tracy1996orthogonal}; Gaussian matrices are exactly solvable, and so the limit may be characterized via exact formulas. As other ensembles are not exactly solvable, new approaches were developed by a variety of researchers proving the universality of the Tracy-Widom distributions for general Wigner matrices \cite{soshnikov1999universality,lee2014necessary,erdos2012rigidity,tao2010random}.

It is natural to ask whether an analog of the Berry-Esseen CLT holds for the convergence in \eqref{eqn:TW}; in particular does the convergence hold at rate $N^{-1}$, the square root of the number of variables in the Wigner matrix?

Progress on this problem has so far been limited. In order to understand possible approaches to this problem it is important to distinguish between two facets which are not present in the simpler case of the classical CLT. Recall above that the Tracy-Widom distributions were first computed from limits of Gaussian matrices (the GOE and GUE) via exact formulas.  There are several works that, via analysis of these formulas, have studied the rate of convergence of the GOE/GUE eigenvalues to the Tracy-Widom distribution, as well as other random matrix ensembles \cite{choup2008edgeworth,choup2006edgeworth,ma2012accuracy, el2006rate,johnstone2008multivariate}. 
For the ensembles considered here, the best results are $\O (N^{-2/3})$ for the GUE and $\O (N^{-2/3})$ for the GOE of odd\footnote{A rate of $\O(N^{-1/3+\eps})$ for the even dimension GOE can be deduced from the odd dimensional result via considering the eigenvalues of a minor.} dimension of Johnstone and Ma \cite{johnstone2012fast}.

The second aspect of this question (and the aspect with which we are concerned here) is one of \emph{quantitative universality}. At a high level, edge universality is proven by showing that the eigenvalue distribution of \emph{any} Wigner matrix is close to that of the GOE/GUE (of the same dimension), and so the convergence in \eqref{eqn:TW} is a consequence of the GOE/GUE convergence. The question then becomes whether or not one can provide good estimates of the KS distance between the largest eigenvalue of a Wigner matrix and that of the corresponding Gaussian ensemble.

The first contribution along this line is in Bourgade's work \cite{bourgade2018extreme}, whose results show that this KS distance is $\O(N^{-2/9+\eps})$ for general Wigner matrices. 
The important work of Schnelli and Xu \cite{schnelli2022convergence} (see also \cite{schnelli2023convergence,schnelli2023quantitative}) gave a new and comprehensive approach to quantitative edge universality and showed that the KS distance between Wigner matrices and the Gaussian ensembles is $\O (N^{-1/3})$, which is the sharpest available result. Combined with the results listed above for the Gaussian ensembles, they deduced  rate of convergence of $\O (N^{-1/3})$ for general Wigner matrices to the Tracy-Widom distribution.

Our contributions to this problem are as follows:
\begin{enumerate}[label=(\roman*)]
\item We will show that the result $\O (N^{-1/3} )$ of \cite{schnelli2022convergence} is optimal in the sense that we will find large classes of Wigner matrices with eigenvalues $\{ \lambda_i \}_i$ so that,
\beq
\sup_{r \in \rr} \left| \pp\left[ N^{2/3} ( \lambda_1 -2 ) \geq r \right] -\pp\left[ N^{2/3} ( \mu_1 -2 ) \geq r \right] \right| \geq \frac{c_\eps}{N^{1/3+\eps}}
\eeq
for any $\eps >0$, where $\mu_1$ is the largest eigenvalue of the corresponding Gaussian ensemble.
\item For a smooth class of Wigner matrices we have,
\beq
\sup_{r \in \rr} \left| \pp\left[ N^{2/3} ( \lambda_1 - \lambda_2 ) \geq r \right] -  \pp\left[ N^{2/3} ( \mu_1 - \mu_2 ) \geq r \right] \right| \leq C_\eps \frac{N^{\eps}}{N}
\eeq
for all $\eps >0$.
\end{enumerate}
The reason for the lower bound is the fact that if the fourth cumulant of a Wigner matrix is non-zero (recall the higher cumulants of the Gaussian distribution vanish) then the expectation $\ee[\lambda_1]$ differs from that of the GOE/GUE by a term of order $N^{-1}$. It is natural to ask whether or not one can improve the error in (i) to $\O (N^{-2/3})$ by subtracting this shift. We discuss pursuing this in greater detail in the next sections, but state for now that the result in (i) already  demonstrates the differences between the random matrix case and the classical CLT; in the classical CLT the KS distance between a sum of i.i.d. Gaussians and general i.i.d. random variables is $\O ( n^{-1/2})$, whereas in the random matrix case, even if the Gaussian ensembles converge very quickly, one would have to choose \emph{model-dependent} shifts and scalings to obtain fast convergence for other models (if possible - we expect that it is). 

The results in (ii) above show that the eigenvalue gaps are in general better behaved than the eigenvalue positions themselves. The individual eigenvalue positions are subject to correlated non-universal shifts; taking the difference in positions removes this effect.

We now return to the discussion of homogenization. 
Dyson Brownian motion is the solution $\{ x_i (t) \}_{i=1}^N$ to the system of stochastic differential equations given by,
\beq \label{eqn:def-dbm}
\d x_i (t) = \sqrt{ \frac{2}{N \beta}} \d B_i (t) + \frac{1}{N} \sum_{j=1, j \neq i}^N \frac{\d t}{x_i (t) - x_j(t)} - \frac{x_i(t)}{2} \d t.
\eeq
Above, $\beta \geq 1$ and the $\{ B_i (t) \}_{i=1}^N$ are a family of standard i.i.d. Brownian motions. We will denote by $\{ \lambda_i (t) \}_{i=1}^N$ and $\{ \mu_i(t) \}_{i=1}^N$ to be the solution to \eqref{eqn:def-dbm} with the \emph{same Brownian motion terms} with initial data the eigenvalues of two independent generalized Wigner matrices, $\{ \lambda_i (0) \}_{i=1}^N$ and $\{ \mu_i (0) \}_{i=1}^N$.  The eigenvalues will always be labeled in descending order.  

By considering the equation satisfied by the difference between the two flows, $u_i (t) := \e^{t/2} (\lambda_i (t) - \mu_i (t))$, one finds a parabolic differential equation,
\beq \label{eqn:intro-parabolic}
\del_t u_i (t) = \frac{1}{N} \sum_{j \neq i } \frac{u_j(t) - u_i(t)}{(\lambda_i (t) - \lambda_{i+1}(t) )( \mu_i (t) - \mu_{i+1} (t) )}.
\eeq
In this context,  homogenization amounts to proving that $u_i (t)$ is to leading order some deterministic function of the initial data, this function involving the fundamental solution of the expected limiting parabolic PDE. Homogenization results are important as they allow one to study general random matrix eigenvalues (the eigenvalues $\{ \lambda_i (t)\}_i$)  by comparing them directly to the invariant Gaussian ensembles in a high probability sense with effective error estimates. 

As stated above, the homogenization of DBM was initiated in the work \cite{fixed-wig}. This work established homogenization for long times $t \sim N^{-c}$ with a small polynomial error. Importantly, this allowed for the proof of fixed energy universality of Wigner matrices of all symmetry classes. Homogenization in a more general context was established in \cite{fixed}; this work proved homogenization for  general classes of initial data for optimal time scales and with a small polynomial error. 

In the breakthrough work \cite{bourgade2018extreme}, Bourgade gave a new approach to homogenization again in the case of Wigner matrix initial data. One of the key insights was finding a new observable satisfying  a stochastic advection equation. Through an iterative scheme leveraging the maximum principle for \eqref{eqn:intro-parabolic}, Bourgade was able to prove homogenization in the spectral bulk with optimal error estimates. In particular, this allowed him to prove the universality of the minimal eigenvalue gap  for smooth classes of Wigner matrices.

Our contribution is to prove a homogenization result for DBM valid throughout the entire spectrum with optimal estimates. Our approach is based on the iterative scheme and stochastic advection observable of \cite{bourgade2018extreme}. However, our work is not merely a simple extension of \cite{bourgade2018extreme} as there are significant obstacles to treating the edge (and all intermediate scales between the bulk and edge). In the next section we state our results on homogenization and then discuss the problems that arise and how our methods overcome them in Section \ref{sec:method}.

We remark here that the homogenization result in the bulk \cite{bourgade2018extreme} has been applied to solve several other bulk universality questions in random matrix theory \cite{bourgade2025optimal,landon2022almost}. Therefore there is significant motivation to establish homogenization throughout the spectrum,  as these results will likely play a role in future developments on universality results for problems that concern the entire spectrum (not just the bulk). We outline these and other possible future applications below.

\subsection{Homogenization at the edge of DBM}

The model of matrices we will consider is the following.

\bed[Generalized Wigner matrices]
We say that an $N \times N$ self-adjoint $H$ matrix is a generalized Wigner matrix if its upper triangular part consists of independent, centered random variables so that
\beq
\frac{1}{CN} \leq \ee[ |H_{ij}|^2] =: S_{ij} \leq \frac{C}{N}
\eeq
for some $C>0$ and all $1 \leq i, j\leq N$ and $\sum_{i=1}^N S_{ij} =1$ for all $j$. We say that $H$ is real symmetric if the matrix elements are real, and complex Hermitian if they are complex. In the complex Hermitian case we further require that the real and imaginary parts are independent and $\ee[ (H_{ij})^2] =0$ for all $i \neq j$.

We further assume that there is a $C>0$ so that
\beq \label{eqn:tail-def}
\pp\left[ \sqrt{N} | H_{ij}| > r \right] \leq C \e^{- r^{C^{-1} }}
\eeq
for all $r >0$.

A generalized Wigner matrix is called a Wigner matrix if, in the real symmetric case we have $S_{ij} = \frac{1+\delta_{ij}}{N}$ and $S_{ij} = \frac{1}{N}$ in the complex Hermitian case.
\eed

One could relax the  tail requirement \eqref{eqn:tail-def} somewhat, but we assume it for technical convenience. 

In order to state our result on homogenization we need to define several quantities from random matrix theory. 
Wigner's semicircle law is defined by,
\beq
\rhosc (E) := \1_{ \{ |E| \leq 2  \} } \frac{\sqrt{4-E^2}}{ 2 \pi} ,
\eeq
and is a good approximation to the empirical eigenvalue distribution of a generalized Wigner matrix for large $N$. 

The Stieltjes transform of $\rhosc(E)$ is,
\beq
\msc (z) := \int_\rr \frac{ \rhosc (x) \d x}{x - z } = \frac{-z+\sqrt{z^2-4}}{2}.
\eeq
We define the quantiles of the semicircle distribution by
\beq
\int_{ \gamma_i}^2 \rhosc (x) \d x = \frac{i-1/2}{N}
\eeq
for integer $1 \leq i \leq N$. 
We define
\beq
\rho_k := \sqrt{4 - \gamma_k^2}  \asymp \left( \frac{ \min\{ k, N+1-k \} }{N} \right)^{1/3} \asymp \rhosc ( \gamma_k).
\eeq
Above, $a_N \asymp b_N$ means that there is a $C>0$ so that $C^{-1} a_N \leq b_N \leq C a_N$. 
Let $\{ \lambda_k (t)\}_k$ and $\{\mu_k (t)\}_k$ satisfy DBM \eqref{eqn:def-dbm} with same Brownian terms. We will always assume that the initial data are the eigenvalues of two generalized Wigner matrices. 

Associated to DBM are the characteristics
\beq \label{eqn:char-def-1}
z_t := z_0 \cosh (t/2) + \sqrt{ z_0^2 -4} \sinh (t/2).
\eeq
Note that $z_t$ satisfies,
\beq \label{eqn:char-def-2}
\dot{z_t} = \msc (z_t) + \frac{z_t}{2}.
\eeq
To emphasize the initial data we will denote $z_t (z_0)$ the characteristic started at $z_0$. Then we introduce
\beq
\gamkt := \lim_{ \eta \dto 0} z_t ( \gamma_k + \i \eta ).
\eeq
We define,
\beq \label{eqn:bu-def}
\bu_k (t) := \frac{1}{N \Im[ \msc ( \gamkt) ]} \Im \left[ \sum_{j=1}^N \frac{ \lambda_j (0) - \mu_j (0)}{ \gamma_j - \gamkt } \right]
\eeq

Our main result on homogenization is the following. 

\bet \label{thm:main-homog} Let $\{ \lambda_i (t) \}_{i=1}^N$ and $\{ \mu_i (t) \}_{i=1}^N$ be two solutions of DBM with the same Brownian motion terms with initial data given by two  real symmetric or complex Hermitian generalized Wigner matrices in the case of $\beta=1, 2$ respectively. 

Let $\eps >0$. For any $D>0$ we have for all sufficiently large $N$ that,
\beq \label{eqn:main-homog-est}
\pp\left[ \left| \lambda_k (t) - \mu_k (t) - \e^{-t/2} \bu_k (t) \right| \leq \frac{ N^{\eps}}{N^2 t \rho_k (t+\rho_k)^2} \mbox{ } \forall \mbox{ } 1 \leq k \leq N {\normalfont \mbox{ and }} t \in (0, 1) \right] \geq 1  - N^{-D}.
\eeq
Above, $\bu_k (t)$ is defined as in \eqref{eqn:bu-def}. 
\eet

Bourgade's work \cite{bourgade2018extreme} proves the above result for indices $ \alpha N \leq k \leq (1- \alpha ) N$ for fixed $\alpha >0$. In this case  $\rho_k \geq c$ for some $c >0$ and so the error term simplifies to $N^{-2} t^{-1}$. 

We expect the error in \eqref{eqn:main-homog-est} to be optimal (up to the $N^{\eps}$ factor) for the following reason. Note that the only difference between $\bu_k (t)$ and, say, $\bu_{k+1} (t)$ is the starting location for the characteristic: $\gamma_k^t$  vs. $\gamma_{k+1}^t$. What is used about the quantiles is that $\gamma_k$ is a good approximation for $\lambda_k(t)$ and $\mu_k (t)$. However, $\gamma_{k+1}$ is also a good approximation for the locations of $\lambda_k, \mu_k$ as the difference $\gamma_{k+1} - \gamma_k$ is the same order as the fluctuation scale of $\lambda_k, \mu_k$ around $\gamma_k$. Therefore, one expects that $\bu_{k+1} (t)$ should be as good an approximation to $\lambda_k(t) - \mu_k (t)$ as $\bu_k (t)$; the error in \eqref{eqn:main-homog-est} is the same order as the difference $| \bu_k(t) - \bu_{k+1} (t) |$.

The fluctuation scale of $\lambda_k (t)$ is $(N \rho_k)^{-1}$. For $t$ close to $1$, the estimate in \eqref{eqn:main-homog-est} beats the natural fluctuation scale by a full factor of $1/N$.

\subsubsection{Discussion of methodology} \label{sec:method}

At a high level, our work proceeds similarly to Bourgade \cite{bourgade2018extreme} in that we prove Theorem \ref{thm:main-homog} through an iterative scheme that in particular, incorporates estimates arising from an observable $f_t (z)$ (see \eqref{def: f_u and s_u}) that satisfies a specific SDE (the \emph{stochastic advection equation}; see \eqref{DE: f_u(z_(t-u))}). 

However, significant obstacles arise that prevent this from being a simple generalization of Bourgade's method. The first such obstacle is that is that away from the bulk the natural spatial scale and time scales diverge; they are both $N^{-1}$ in the bulk whereas at the edge they are respectively $N^{-2/3}$ and $N^{-1/3}$. This introduces difficulties if, in particular, one wishes to find the optimal error in \eqref{eqn:main-homog-est}. For simplicity, let us assume that $t$ is of order $1$. The iterative procedure involves localizing the problem in a spatial window of size $\omega$ around a location $\gamma_k$ and studying the evolution over a time interval of length  $T$ ending at time $t$.\footnote{We use different notation in the actual proof, as the window is best understood in terms of eigenvalue indices and is denoted $R$; the time interval is of length $|u-t|$} In \cite{bourgade2018extreme}, the true solution of the eigenvalue difference equation is approximated in this space-time window by the constant solution $\bu_k (t)$. This naturally introduces errors of size
\beq \label{eqn:method-errors}
\omega \times \left| \del_k \bu_k (t) \right| \sim  \frac{\omega}{N} , \qquad T \times \left| \del_s \bu_k (s) \right| \sim \frac{T}{N}
\eeq
where $\del_k$ is somewhat formal and is interpreted as a spatial derivative. When $t \sim 1$, both of these derivatives are of the same size as $|\bu_k (t)|$ which is $\O (N^{-1})$. 

In the the last step of the iteration, one will choose $\omega$ and $T$ to be the natural length scales associated to DBM at the spectral index $k$. These are $(N \rho_k)^{-1}$ (the natural fluctuation of $\lambda_k (t)$) and $(N \rho_k^2)^{-1}$ (the time to local equilibrium of DBM near index $k$). In the bulk, the errors in \eqref{eqn:method-errors} are both $\O (N^{-2})$ which is acceptable. At the spectral edge (and at all intermediate scales), the error from the time derivative is much larger than the spatial error ($N^{-4/3}$ vs $N^{-5/3}$ at the edge) and is not acceptable; proceeding in this fashion would give an error of $\O (N^{-4/3})$ in \eqref{eqn:main-homog-est}, much worse than the expected optimal estimate. 

We overcome this problem by approximating the solution to \eqref{eqn:intro-parabolic} in this space-time window by $\bu_k (s)$; i.e., by a constant function in space but not in time. However, this approximating introduces extra terms in the main stochastic advection equation calculation in \cite{bourgade2018extreme}, coming from the time derivative $\del_s \bu_k (s)$. These terms are of larger order than  the optimal error in \eqref{eqn:main-homog-est} and so cannot simply be dropped. However, there is an additional approximation in \cite{bourgade2018extreme} that is lossy once one leaves the simpler bulk setting, where the solution to \eqref{eqn:intro-parabolic} is approximated by a ``short-range'' equation where only interactions between eigenvalue indices with $|i-j| \leq \ell$ are kept. The missing term that is dropped is a ``long-range'' interaction approximated by,
\beq
\frac{1}{N} \sum_{ j : |k-j| > \ell} \frac{ \bu_j (s) - \bu_k (s)}{ ( \lambda_i (s) - \lambda_k (s))( \mu_i (s) - \mu_k (s) ) }.
\eeq
We show that these two contributions cancel by a deterministic calculation (see Proposition \ref{prop:deltbu}) showing  that $\bu_k (s)$ solves the limit of \eqref{eqn:intro-parabolic},
\beq
\del_s g (x, s) = \int \frac{ g(y,s) - g(x, s)}{(x-y)^2} \rhosc (y) \d y
\eeq
with the integral interpreted in the principal value sense.

The second obstacle is that the iterative procedure based on spatial localization in \cite{bourgade2018extreme} is not guaranteed to work at the spectral edge, due to the inhomogeneity induced by the curvature of the semicircle distribution. Roughly, there are two parameters in the spatial localization procedure; one approximates \eqref{eqn:intro-parabolic} by a short-range equation involving interactions up to distance $\ell$, and studies the problem for eigenvalue indices up to range $R$ from a fixed index $k$ at which we are trying to prove homogenization; let us choose $k=1$ to illustrate. As in \cite{bourgade2018extreme}, this introduces an error of size $\frac{\ell}{R}$ which must be $o(1)$ in order to gain in the iterative procedure. Note that this constrains $\ell$ to be much smaller than $R$. However, there are additional errors which essentially come from the fact that $\Im[\msc]$ is much larger near $\gamma_R$ than near the spectral edge $\gamma_1 \approx 2$. This error gets worse if $R$ is larger, but is counteracted with a factor of $\ell$ in the denominator. By careful estimation we find that this error is, roughly, $\frac{R^{2/3}}{\ell}$. There is therefore room to choose $R^{2/3} \ll \ell \ll R$ and so the procedure works. Note that if instead of $\frac{R^{2/3}}{\ell}$ we obtained something like $\frac{R^{2/3}}{\ell^{2/3}}$, then there would be no room to gain in the iterative procedure and so sloppier analysis would prevent the method from working.

\subsection{Quantitative universality at the edge}

For our results on quantitative universality, we make additional smoothness assumptions on the generalized Wigner matrices: 
\bed [Smooth generalized Wigner matrices] \label{def:smooth}
 We say that a generalized Wigner matrix $H$ is smooth if the density real and imaginary parts of $\sqrt{N} H_{ij}$ are of the form $\e^{-V(x)}$ where $V$ is a smooth function whose derivatives grow at most polynomially,
 \beq
 |V^{(k)} (x) | \leq C_k (1+ |x|)^{C_k}
 \eeq
 where the $C_k$ are uniform in the matrix indices $i, j$.
\eed

While ideally one would like to remove the smoothness assumption above, it is a natural class to work in given our methods; we discuss this in detail below after stating our results. Our first result on quantitative universality is the following.

\bet \label{thm:gap-univ}
Let $\{ \lambda_i \}_{i=1}^N$ be the eigenvalues of a smooth generalized Wigner matrix, and let $\{ \mu_i \}_{i=1}^N$  be the eigenvalues of the GOE. For any $\eps >0$, there is a constant $C_\eps$ so that
\beq \label{eqn:main-gap-univ}
\sup_{s \in \rr} \left| \pp \left[ N^{2/3}(\lambda_1 - \lambda_2 ) \geq s \right] -  \pp \left[ N^{2/3}(\mu_1 - \mu_2 ) \geq s \right] \right| \leq C_\eps N^{\eps-1} .
\eeq
\eet

As discussed in the introduction, the above result states that the gap between the top two eigenvalues has a universal distribution with an error rate that is the analog of the Berry-Esseen CLT: the square root of the number of variables in the problem. We are not aware of any work studying the rate of convergence of the gap of the GOE/GUE to the gap between the top two particles in the Airy$_\beta$ point process (a small polynomial error rate could be deduced from \cite{landon2020edge} but this falls far short of the $\O (N^{-1})$ error above). It is possible that one could study this via the exact formulas available in the Gaussian cases but these methods are out of the scope of this paper, being vastly different from the tools we introduce.

Let us explain the quick proof of the above theorem given Theorem \ref{thm:main-homog}, which will also explain the smoothness restriction in Definition \ref{def:smooth}. Taking $t \approx 1$\footnote{We take $t= N^{-\eps}$ but ignore this for the sake of exposition} we will see that $| \bu_2 (t) - \bu_1 (t) | $ is of order $N^{-5/3}$ and so we see that
\beq
| ( \lambda_2 (t) - \lambda_1 (t) ) - (\mu_2 (t) - \mu_1 (t) ) | \leq N^{-5/3+\eps}
\eeq
with very high probability. Together with an elementary bound on the density of the random variable $N^{2/3} ( \mu_2(t) - \mu_1 (t) )$ (see Proposition \ref{prop:density}) this proves the theorem for $\lambda_i (t)$. Now the distribution of $\{ \lambda_i (t) \}_{i=1}^N$ is the same as the eigenvalue distribution of $H_t = \e^{-t/2} H_0 + \sqrt{1-\e^{-t}} G$ where $H_0$ is the matrix whose eigenvalues are the initial data $\{ \lambda_i (0) \}_{i=1}^N$, and $G$ is an independent GOE/GUE matrix. Broadly, there are two ways to remove the Gaussian component and/or approximate general matrices by ones that look like $H_t$. The first is based on moment matching and analyzing the Green function/resolvent $\frac{1}{H-z}$ and the other on the reverse heat flow \cite{local-relax,erdos2010bulk}. We rely on the reverse heat flow as it provides very strong estimates and no additional work; however, this introduces the smoothness assumption. 

We do not pursue the approach based on analyzing the resolvent as, without substantial new advances, it will not be able to produce an estimate as sharp as \eqref{eqn:main-gap-univ}. As can be seen from the passage from \eqref{eqn:main-homog-est} to \eqref{eqn:main-gap-univ}, one needs to access the eigenvalue gap $\lambda_2 - \lambda_1$ down to a scale $\eta = N^{-2/3-1} = N^{-5/3}$. This is roughly equivalent to controlling the resolvent $\frac{1}{H-z}$ for $z = 2 + \i \eta = 2 + \i N^{-5/3}$. 

The work of Schnelli and Xu \cite{schnelli2022convergence}, which proves a convergence rate of $N^{-1/3+\eps}$, is based on analyzing the resolvent down to the scale $\eta = N^{-1+\eps}$. Their work is further based on extensive and delicate expansions of the resolvent in the matrix entries; such expansions are only convergent above the critical threshold $\eta = N^{-1}$. Below this threshold, expansions seem only to make estimates worse, and the gap between the scale we need $N^{-5/3}$ and the scale $N^{-1}$ is quite substantial. As the methods of \cite{schnelli2022convergence} are already intricate and challenging we defer any possible extensions to shorter scales to future investigation.

It would be interesting to investigate whether resolvent methods combined with our homogenization results could achieve a convergence rate better than $\O (N^{-1/3})$, possibly under the assumptions allowing for the matching of many matrix entry moments, as in \cite{zhang2025quantitative}. We leave this for future work, as substantial resolvent computations would add to the length of the paper and detract from our main focus which is on optimal homogenization throughout the spectrum.

Having discussed the eigenvalue gap at the spectral edge, let us turn to our results for the largest eigenvalue itself.

\bet \label{thm:main-lb}
Let $H$ be a Wigner matrix such that the fourth cumulant of the entries is non-zero. Let $H_t = \e^{-t/2} H + \sqrt{1- \e^{-t} } G$ where $G$ is a GOE matrix, and let $t=\frac{1}{2}$. Then, if $\lambda_1$ denotes the largest eigenvalue of $H_t$ we have,
\beq  \label{eqn:main-lb-1}
\sup_{r \in \rr} \left| \pp\left[ N^{2/3} ( \lambda_1 - 2 ) \geq r \right] - \pp\left[ N^{2/3} ( \mu_1 - 2 ) \geq r \right] \right| \geq \frac{1}{N^{1/3+\eps}}
\eeq
for any $\eps >0$ and $N$ large enough.

 For any $b \in \rr$ there is a choice of $H$ so that if $\lambda_1$ is the largest eigenvalue of $H_t$ then,
 \beq  \label{eqn:main-lb-2}
\sup_{r \in \rr} \left| \pp\left[ N^{2/3} ( \lambda_1 - 2 ) \geq r \right] - \pp\left[ N^{2/3} ( \mu_1 - 2 ) -bN^{-1/3} \geq r \right] \right| \geq \frac{1}{N^{1/3+\eps}}
\eeq
\eet

Let us interpret the above results. The first estimate \eqref{eqn:main-lb-1} essentially says that the quantitative edge universality result of Schnelli and Xu \cite{schnelli2022convergence} is optimal if one asks for the worst case convergence amongst all possible Wigner matrices. As stated in the introduction, the work \cite{johnstone2012fast} showed that a convergence rate of $\O ( N^{-2/3})$ could be established for the odd-dimensional GOE to the Tracy-Widom distribution, if one allows for a shift of order $\O (N^{-1/3})$ and a similar scaling. It is natural to ask whether this could be implemented for other matrices. The estimate \eqref{eqn:main-lb-2} states that if this is true, then the shift must be model-dependent and not universal among all Wigner matrices. 

This then leads us to a natural question which we believe to be accessible with our work: given the estimates \eqref{eqn:main-homog-est} can one improve on the $\O (N^{-1/3})$ result of \cite{schnelli2022convergence}? We believe this to be the case and leave this for a follow-up paper. We do not include these results in the present work as the methods no longer have anything to do with homogenization and instead rely on analysis of the $\bu_k (t)$; this turns out to be to leading order a linear spectral statistic, an observable that can be analyzed via resolvent approaches.

\subsection{Related literature and applications}

There have been many works studying various forms of quantitative universality at the spectral edges of various random matrix ensembles.  By generalizing the approach of \cite{bourgade2018extreme} to the case of sample covariance matrices, the works  \cite{wang2024quantitative,wang2022optimal} prove quantitative universality results for the largest and smallest eigenvalues of matrices of the form $X X^T$. 

The work \cite{schnelli2023convergence} generalized the resolvent approach of \cite{schnelli2022convergence} to the case of sample covariance matrices; this was further developed in \cite{bucht2025quantitative} to prove estimates on the convergence rate for sparse random matrices to the Tracy-Widom distribution. 

For convergence to the Tracy-Widom distribution, all of these works prove rates that are at best $\O (N^{-1/3})$. This is due to the fact that none of these works establish homogenization near the spectral edge which seems to be a prerequisite for studying the problem at scales $\eta \ll N^{-1}$. 

In a different direction, the work \cite{zhang2025quantitative} studies quantitative universality of the gap in the spectral bulk, by generalizing the four moment method to match higher moments (under support assumptions on the distribution of the Wigner matrix elements).

On the dynamical side, convergence results for DBM sufficient for edge universality have been established in several different contexts. All of these works give weaker estimates than our homogenization results.  The work \cite{bourgade2014edge} proved convergence at the edge for general $\beta$-ensembles using an energy estimate and deduced edge universality of general $\beta$-ensembles and generalized Wigner matrices. The work \cite{landon2017edge} established a convergence result for Dyson Brownian motion for general initial data. 

Our work relies on a finite speed of propogation estimate for the short-range equation mentioned above. In the form we use, this estimate was first established in the work \cite{bourgade2017eigenvector} (see also \cite{erdHos2015gap} for an early form as well as \cite{fixed,landon2017edge} for related estimates). This work analyzed the convergence to equilibrium of the eigenvectors of generalized Wigner matrices and used this to prove the universality of the eigenvector distribution. 

\vspace{5 pt}

\noindent{\bf Applications of homogenization.} As mentioned above, the bulk homogenization result \cite{bourgade2018extreme} has been applied to solve several challenging universality problems in random matrix theory. This is due in part to the fact that the homogenization result allows one to access the eigenvalues on scales a full factor of $N^{-1}$ beyond their natural fluctuation scale (which is $(N \rho_i)^{-1}$). In \cite{bourgade2018extreme}, the result was applied to prove universality of the smallest gap between consecutive eigenvalues; in the case of real symmetric matrices, this quantity is of size $\O (N^{-3/2})$, requiring a very fine homogenization estimate. The homogenization result played an important role in work of the first author together with Sosoe \cite{landon2022almost} which proved an almost-optimal result for the regularity conditions in the CLT for linear spectral statistics of Wigner matrices; here the fact that the result allowed access to scales as small as $\O (N^{-2})$ in the bulk was particularly important. The homogenization result was also used in the work \cite{bourgade2025optimal} to prove results on the maximum of the characteristic polynomial and an optimal rigidity result for general Wigner matrices. 

Given the above results (which are all limited to the spectral bulk) we expect that - in addition to the important quantitative edge universality results mentioned above - our work will have numerous future applications. For example, our results are the first step in extending the results of \cite{landon2022almost} to the spectral edges, which is a delicate problem. Moreover, we have learned of the forthcoming work \cite{bouci} which independently develops relaxation estimates throughout the spectrum for the purposes of estimating the spectral form factor of Wigner matrices (see also \cite{cipolloni2023spectral}).

\section{Preliminaries}

\subsection{Notation}

Let us define the control parameter,
\beq
\varphi := \e^{ ( \log \log N)^3}.
\eeq
We will use $C, c>0$ to denote large and small constants whose value may vary from line to line. If $a_N$ and $b_N$ are two sequences of positive numbers we say that
$
a_N \asymp b_N
$
if there are $c,C>0$ so that $c a_N \leq b_N \leq C a_N$ for all sufficiently large $N$. We use $a_N \lesssim b_N$ if there is a $C>0$ so that $a_N \leq b_N$ for all sufficiently large $N$.

\bed[Overwhelming probability] If $\I_N$ is some abstract index set and $\F^{(N)}_i$ for $i \in \I_N$ is a family of events, we say they hold with overwhelming probability if for all $D>0$ we have that
\beq
\sup_{ i \in \I_N} \pp\left[ (\F^{(N)}_i)^c \right] \leq N^{-D}
\eeq
for all sufficiently large $N$.
\eed
We will typically use the notation,
\beq
z = E + \i \eta, \qquad E \in \rr, \quad \eta >0 .
\eeq
We also recall,
\beq
\rho_k := \sqrt{4 - \gamma_k^2}  \asymp \left( \frac{ \min\{ k, N+1-k \} }{N} \right)^{1/3} \asymp \rhosc ( \gamma_k)
\eeq
We will define,
\beq
\kappa(z) := \mathrm{dist} (z, \{-2, 2\} ), \qquad \kappa_j = \kappa ( \gamma_j) \asymp \rho_j^2 .
\eeq
We define
\beq
[\![A, B]\!] := \{ n \in \zz : A \leq n \leq B\} ,
\eeq
and for $i \in [\![1, N]\!]$,
\beq
\hati := \min\{ i, (N+1-i ) \}
\eeq
We denote the upper half plane $\hh_+ := \{ z \in \cc : \Im[z] > 0 \}$.

\subsection{Interpolating DBM, rigidity and local laws}

Let $\{ \lambda_i (0) \}_{i=1}^N$ and $\{ \mu_i (0) \}_{i=1}^N$ be the eigenvalues of two generalized Wigner matrices. For $\nu \in [0, 1]$ we let $\{ \xnu_i (t) \}_{i=1}^N$ be the solution to \eqref{eqn:def-dbm} with initial data $\xnu_i (0) = \nu \lambda_i (0) + (1- \nu) \mu_i (0)$. Note that $x^{(1)}_i (t) = \lambda_i (t)$ and $x^{(0)}_i (t) = \mu_i (t)$. Following \cite{bourgade2018extreme} we define,
\beq
\fu^{(\nu)}_i (t) := \e^{t/2} \del_\nu \xnu_i (t),
\eeq
which solves the equation
\beq
\del_t \fu^{(\nu)}_i (t) = ( \sB \fu^{(\nu)} )_i := \frac{1}{N} \sum_{j \neq i} \frac{ \fu^{(\nu)}_j (t) - \fu^{(\nu)}_i (t) }{(\xnu_i (t) - \xnu_j (t) )^2}.
\eeq
Throughout the paper we will judiciously suppress the dependence of various quantities on $\nu$ and denote $x_i (t) = \xnu_i (t)$ and $\fu_i (t) = \fu^{(\nu)}_i (t)$, etc. We introduce the empirical Stieltjes transform,
\beq
s_t (z) = s_t^{(\nu)} (z) := \frac{1}{N} \sum_{i=1}^N \frac{1}{ \xnu_i(t) -z } .
\eeq

\bet[Rigidity and local laws] With overwhelming probability we have for all $ \nu \in [0, 1]$, $t \in [0, 1]$ and $k \in [\![ 1, N ]\!]$,
\beq \label{eqn:rig}
| x_i (t) - \gamma_i | \leq \frac{ \varphi^{\frac{1}{100}} }{N^{2/3} \hati^{1/3} }.
\eeq
Whenever the above estimates hold we have for all 
$z \in \hh_+$ with $|z| \leq 10$,
\beq \label{eqn:ll}
| \msc (z) - s_t (z) | \leq \frac{ \varphi^{\frac{1}{100}}}{N(\eta + \1_{ \{ |E| > 2 + \varphi^{\frac{1}{100}} N^{-2/3} \} } \kappa(E) )} 
\eeq
\eet
\proof The rigidity estimates are proven similarly to \cite[Lemma 2.3]{bourgade2018extreme}. One can deduce the estimate for $s_t (z)$ in the same manner as the proof of \cite[Theorem 10.3]{semi-lec}. \qed

\subsection{Properties of characteristics}

We record here that by direct calculation using \eqref{eqn:char-def-1}, \eqref{eqn:char-def-2} one sees that
\beq \label{eqn:msc-evolution}
\del_t \msc(z_t) = - \frac{ \msc(z_t)}{2}
\eeq

The following is \cite[Lemma 2.2]{bourgade2018extreme}. 

\bel Let $z_t$ be a characteristic starting from some $z_0 = z$ with $|z-2| < \frac{1}{10}$ and $\Im[z] >0$. 
Uniformly in $0 < t < 1$ we have,
\beq
\Re[z_t - z_0] \asymp \frac{ t a(z)}{ k(z)^{1/2}} + t^2, \qquad \Im[z_t -z_0] \asymp \frac{ b(z)}{ k(z)^{1/2}} t
\eeq
where 
\beq
a(z) = \mathrm{dist} (z, [-2, 2] ), \quad b(z) = \mathrm{dist} (z, [-2, 2]^c), \quad \kappa(z) = \mathrm{dist} (z, \{ -2, 2\} ).
\eeq
\eel

Define the region
\beq
\D := \{ z : \eta \in (0, 1), E \in [-2, 2]\} \cup \{ z : \eta \in [ \varphi /N^{2/3}, 1], |E| \in [2,2+  \varphi / N^{2/3} ] \}.
\eeq
It is clear that if $z \in \D$ then $a(z) \asymp \eta$ and $b(z) \asymp \kappa(E) + \eta$ so that uniformly for $z \in \D \cap \{ z : |z-2| < 0.1 \}$ we have,
\beq \label{eqn:char-bourgade}
\Re[z_t - z_0] = t \frac{ \eta}{ ( \kappa(E) + \eta)^{1/2}} + t^2, \qquad \Im[z_t-z_0] \asymp ( \kappa(E) + \eta)^{1/2} t
\eeq

\bel
Let $z_t $ be a characteristic starting at $z_0 = z = E + \i \eta \in \D$ and $\gamma \in [-2, 2]$. Then,
\beq \label{eqn:char-est}
|z_t - \gamma| \asymp |E - \gamma| + t^2 + t \sqrt{ \kappa(E) + \eta} + \eta
\eeq
\eel
\proof If $\kappa(z) >c$ then the RHS of \eqref{eqn:char-est} is $\asymp |E-\gamma| + \eta + t$. In this case the estimate is easy to prove directly from \eqref{eqn:char-def-1} and \eqref{eqn:char-def-2}. We do the case $|z-2| < 0.1$, the case where $z$ is near $-2$ following from symmetry.  

For such $z$, the upper bound is an immediate consequence of \eqref{eqn:char-bourgade}. It remains to prove the lower bound. If $E \geq \gamma$, then this lower bound is again an easy consequence of \eqref{eqn:char-bourgade}. 

It therefore remains to prove the lower bound when $ E < \gamma$. Note that in this case we have $|E-\gamma| \leq \kappa(E)$. Moreover, $\Im[z_t] \asymp t \sqrt{ \kappa(E) + \eta} + \eta$, this following from \eqref{eqn:char-bourgade}. Assume first that $t \leq C \sqrt{ \kappa(E) + \eta}$. It suffices to prove that $| \Re[z_t - \gamma] | + \Im[z_t] \geq c |E-\gamma|$ for some $c>0$ (possibly depending on $C>0$). If $|\Re[z_t-z_0] | \leq \frac{|E-\gamma|}{2}$ then we have $| \Re[z_t-\gamma] | \geq \frac{|E-\gamma|}{2}$ by the reverse triangle inequality. On the other hand, if $| \Re[z_t-z_0] | \geq \frac{|E-\gamma|}{2}$ then by \eqref{eqn:char-bourgade} we get that
\beq
t \frac{ \eta}{\sqrt{ \kappa(E) + \eta}} + t^2 \geq c_1 |E-\gamma|
\eeq
for some $c>0$. One of the terms on the LHS must be greater than $\frac{c_1}{2} |E - \gamma|$. If it is $t^2$ then 
\beq
\Im[ z_t] \geq c_2 t \sqrt{ \kappa(E) + \eta } \geq c_3 |E-\gamma|^{1/2} \kappa(E)^{1/2} \geq c_3 |E-\gamma|
\eeq
and so the desired lower bound holds. On the other hand, if $t \frac{ \eta}{ \sqrt{ \kappa(E) + \eta}} \geq \frac{c_1}{2} |E-\gamma|$ then,
\beq
t \sqrt{ \kappa(E) + \eta} \geq \frac{c_1}{2} \frac{ \kappa(E) + \eta}{\eta} |E-\gamma| \geq \frac{c_1}{2} |E-\gamma|
\eeq
and so $\Im[z_t] \geq c |E-\gamma|$. This completes the argument in the case that $t \leq C \sqrt{ \kappa(E) + \eta}$, for any large $C>0$.

So we now assume that $t \geq C \sqrt{ \kappa(E) + \eta}$. By taking $C$ sufficiently large we can assume that
\beq
|\Re[z_t-z_0]| \geq c t^2 \geq 10 (\kappa(E) + \eta) \geq 10|E-\gamma|.
\eeq
Then,
\beq
| \Re[z_t-\gamma]| \geq |\Re[z_t-z_0] | - |E-\gamma| \geq 5 |E-\gamma| - |E-\gamma| + \frac{1}{2} \Re[z_t-z_0] | \geq |E - \gamma| + c t^2
\eeq
for some $c>0$. This completes the proof. \qed

\bel Let $z_t$ be a characteristic starting at $z_0 = z \in \D$. Then uniformly  in $t \in [0, 1]$ we have,
\beq \label{eqn: msc(z_t)}
\Im [ \msc (z_t) ] \asymp \sqrt{ \kappa(E) + \eta}
\eeq
\eel
\proof From the equation \eqref{eqn:msc-evolution} 
we see that $\Im[\msc (z_t) ] \asymp \Im[ \msc (z)]$. For $z \in \D$ we have $\Im[\msc(z)] \asymp \sqrt{ \kappa(E) + \eta}$ using,
\beq
\Im[ \msc(a+ \i b) ] \asymp \begin{cases} \sqrt{\kappa(a) + b }, & a \in [-2, 2]\\ \frac{b}{ \sqrt{\kappa(a) + b } } , & a \notin [-2, 2] \end{cases}
\eeq
which follow from \cite[(3.3)]{semi-lec}. \qed

\bel
Let $z_t$ be a characteristic starting from some $z_0 = z \in \D$.  Assume that 
\beq \label{eqn:outside-assump}
\eta + t \sqrt{ \kappa(E) + \eta} > \frac{ \varphi^2}{N (\rhosc(E) + N^{-1/3} )}.
\eeq
  On the event that the rigidity estimates \eqref{eqn:rig} we have that, 
\beq  \label{eqn: |s_u(z_t)-msc(z_t)|}
\left| s_u (z_t) - \msc (z_t) \right| \leq \varphi^{\frac{1}{99}} \frac{1}{N ( t^2 + t \sqrt{ \kappa(E) + \eta} + \eta)} .
\eeq
\eel
\proof We see from \eqref{eqn:char-bourgade} that $\Im[z_t] \geq c ( \eta + t \sqrt{ \kappa(E) + \eta} )$. Fix some constant $C_1 >0$. If $t^2 \leq C_1 ( t \sqrt{ \kappa(E) + \eta} + \eta)$ then we conclude the estimate from \eqref{eqn:ll}. On the other hand, if $t^2 \geq  C_1 ( t \sqrt{\kappa(E) + \eta} +\eta)$ then $t^2 \geq C_1 ( \kappa(E) +\eta)$ and so by taking $C_1$ large enough we conclude from the first estimate in \eqref{eqn:char-bourgade} that 
\beq
\Re[z_t] \geq 2 + c t^2 + 10 \kappa(E).
\eeq
If $\kappa(E) \geq \varphi^{1/2} N^{-2/3}$ then the RHS of the above is greater than $2 + \varphi^{1/2} N^{-2/3}$. If $\kappa(E) \leq \varphi^{1/2} N^{-2/3}$ then \eqref{eqn:outside-assump}  implies that $t^2 \geq \varphi N^{-2/3}$. In either case $\Re[z_t] \geq 2 + \varphi^{1/2} N^{-2/3}$ and $\kappa(\Re[z_t]) \geq c t^2$ and so the estimate follows from \eqref{eqn:ll}. \qed


\section{Properties of homogenized function}

We note here for future reference an alternative formula for $\bu_k (t)$. Since along characteristics we have by \eqref{eqn:msc-evolution}, 
\beq
\msc (z_t) = \e^{ - t/2} \msc (z_0)
\eeq
as well as by definition \eqref{eqn:char-def-1},
\beq
\gamma_k^t = \gamma_k \cosh (t/2) + 2 \Im[ \msc ( \gamma_k ) ] \sinh (t/2) = \gamma_k \cosh (t/2) + \rho_k \sinh(t/2)
\eeq
we see that
\begin{align} 
\bu_k (t) &:= \frac{1}{N \Im[ \msc ( \gamkt) ]} \Im \left[ \sum_{j=1}^N \frac{ \lambda_j (0) - \mu_j (0)}{ \gamma_j - \gamkt } \right] \notag\\
&= \frac{ 2 \e^{t/2} \sinh(t/2)}{N} \sum_{j=1}^N \frac{ \lambda_j(0) - \mu_j (0) }{ | \gamma_j - \gamma_k^t|^2}  \label{eqn:buk-alt-formula} \\
&= \frac{ 2 \e^{t/2} \sinh(t/2)}{N} \sum_{j=1}^N \frac{ \lambda_j(0) - \mu_j (0) }{ ( \cosh(t/2) \gamma_k - \gamma_j)^2 + (  \sinh(t/2) \rho_k)^2} \label{eqn:buk-formula}
\end{align}

\subsection{Regularity}

We will use the formula \eqref{eqn:buk-formula} in this section. In order to analyze it, it is helpful to recall the facts,
\beq
\sinh(t/2) \asymp t, \quad  \rho_k \asymp (k/N)^{1/3}, \quad \kappa(\gamma_j) \asymp (j/N)^{2/3}
\eeq
 which will be used repeatedly in this section.  We first require the following simple lemma, allowing us to handle the denominator in \eqref{eqn:buk-formula}. 
 
\bel \label{lem:buk-denom}
We have uniformly for $t \in [0, 1]$,
\beq 
| \cosh (t/2) \gamma_k - \gamma_j | + t \rho_k \asymp |\gamma_j -\gamma_k | + t^2 +t \rho_k 
\eeq
\eel
\proof This follows from \eqref{eqn:char-est} and the fact that
\beq
| \cosh (t/2) \gamma_k - \gamma_j | + t \rho_k  \asymp | \gamma_k^t - \gamma_j | .
\eeq
\qed

\bel \label{lem:t^2rho}
There is a $c>0$ so that the following holds. For any $A \geq 1$, we have that if $t^2 + t \rho_k \geq A (N \rho_k)^{-1}$ then,
\beq
t \rho_k \geq c A^{1/2} (N \rho_k)^{-1}
\eeq
\eel
\proof Assume that $t \rho_k \leq \frac{A^{1/2}}{2 N \rho_k}$. Then 
\beq
t^2 \geq \frac{A}{N \rho_k} - t \rho_k \geq \frac{A}{2 N \rho_k} \geq \frac{A^{1/2}t  \rho_k}{2}
\eeq
so $t \geq c A^{1/2} \rho_k$. But then $t \rho_k \geq c A^{1/2} \rho_k^2 \geq c A^{1/2} (N \rho_k)^{-1}$ since $\rho_k \geq c N^{-1/3}$ for all $k$. 
\qed

\bel \label{lem:buk-bound} Assume that $ t^2 + t \rho_k \geq (N \rho_k)^{-1}$. 
We have that
\beq \label{eqn:buk-bound}
| \bu_k (t) | \leq C \frac{ \varphi^{1/100}}{N (t + \rho_k ) }
\eeq
Otherwise,
\beq \label{eqn:buk-bound-2}
| \bu_k (t) | \leq C \varphi^{1/100} \left( \frac{ 1}{N (t + \rho_k ) } + \frac{1}{N^2 \rho_k (t+\rho_k)^2 t} \right)
\eeq
\eel
\proof By \eqref{eqn:rig}, \eqref{eqn:buk-formula} and Lemma \ref{lem:buk-denom} we have
\beq \label{eqn:buk-bound-a1}
| \bu_k (t) | \leq C \varphi^{1/100} \frac{t}{N} \sum_{j=1}^N \frac{1}{ ( \gamma_j -\gamma_k)^2 + (t^2 + t \rho_k)^2} \frac{1}{N^{2/3} \hat{j}^{1/3}}.
\eeq
We start with the first estimate. The assumption on $t^2 +t \rho_k$ implies that the sum can be compared to the following integral, 
\begin{align}
 &\sum_{j=1}^N \frac{1}{ ( \gamma_j -\gamma_k)^2 + (t^2 + t \rho_k)^2} \frac{1}{N^{2/3} {\hat{j}}^{1/3}} \leq C \frac{1}{N} \sum_{j=1}^N \frac{1}{ ( \gamma_j -\gamma_k)^2 + (t^2 + t \rho_k)^2} \frac{1}{2-|\gamma_j|} \notag \\
 \lesssim & \int_{-2}^2 \frac{1}{ ( \gamma_k - x)^2 + (t^2 + t \rho_k )^2} \frac{ \rhosc(x)}{ (4-x^2)^{1/2}} \d x \lesssim \frac{1}{ t^2 + t \rho_k}
\end{align}
which yields the first estimate. For the second estimate, every term on the RHS of\eqref{eqn:buk-bound-a1} except for the term $j=k$ may still be estimated about by the integral above. The $j=k$ term is bounded above by,
\beq
\frac{t}{N} \frac{1}{ (\gamma_k - \gamma_k)^2 + (t^2 + t \rho_k)^2}  \frac{1}{N^{2/3} \hat{k}^{1/3}} \leq \frac{C}{N^2 \rho_k ( \rho_k + t)^2 t}
\eeq
yielding the second estimate.  \qed

\bel
Assume that $t^2+ t \rho_k \geq (N \rho_k)^{-1}$. Let $t/2 < s <t $. Then,
\beq \label{eqn:buk-time-bd}
| \bu_k (s) - \bu_k (t) | \leq C \varphi^{1/2} \frac{|t-s|}{Nt (t+ \rho_k ) }
\eeq
\eel
\proof We first note that since
\beq
| \sinh (t/2) - \sinh(s/2) | \leq C |t-s|
\eeq
we can replace the prefactor $\sinh(s/2)$ in the definition of $\bu_k (s)$ in \eqref{eqn:buk-formula} by $\sinh(t/2)$ at a cost of $\O ( |t-s| |\bu_k (s) | /t)$, and by \eqref{eqn:buk-bound} we can absorb this into the RHS of \eqref{eqn:buk-time-bd}. By using the inequalities,
\begin{align}
 & \left| ( \cosh(t/2) \gamma_k - \gamma_j )^2 - ( \cosh(s/2) \gamma_k - \gamma_j )^2 \right| \\
\leq & Ct (t-s) \left( | \cosh(t/2) \gamma_k - \gamma_j | + | \cosh(s/2) \gamma_k - \gamma_j | \right)
\end{align}
and
\beq
\left| \sinh(t/2)^2 - \sinh(s/2)^2 \right| \leq C t (t-s) 
\eeq
we find,
\begin{align}
& \left| \frac{t}{N} \sum_{j=1}^N \left( \frac{\lambda_j (0) - \mu_j (0) }{ ( \cosh(t/2) \gamma_k - \gamma_j )^2 + ( \sinh(t/2) \rho_k )^2} - \frac{\lambda_j (0) - \mu_j (0) }{ ( \cosh(s/2) \gamma_k - \gamma_j )^2 + ( \sinh(s/2) \rho_k )^2} \right)  \right|  \notag \\
 & \leq C t^2 (t-s) \frac{\varphi^{1/2}}{N} \sum_{j=1}^N \frac{ \rho_k^2 + | \gamma_k -\gamma_j | +t^2}{(( | \gamma_j - \gamma_k|^2 + (t^2 + t \rho_k )^2)^2} \frac{1}{N^{2/3} \hat{j}^{1/3}} \notag\\
& \leq C \varphi \frac{t-s}{t} \left[t \sum_{j=1}^N \frac{1}{ ( \gamma_j - \gamma_k)^2 + (t^2 + t \rho_k )^2} \frac{1}{N^{2/3} \hat{j}^{1/3}} \right]
\end{align}
In the proof of Lemma \ref{lem:buk-bound} we showed that the term in the square brackets is $\O ( N^{-1} (t + \rho_k)^{-1} )$. This completes the proof. \qed

\bel
Let $a, b$ be indices and suppose that $| \gamma_a - \gamma_b | \leq C_1 ( t^2 + t ( \rho_a + \rho_b ) )$ for some $C_1 >0$. Then,
\beq \label{eqn:b1}
t^2 + t \rho_a \asymp t^2 + t \rho_b
\eeq
and there is a $c_1 >0$ so that if $ t \rho_b \geq A  (N \rho_b)^{-1}$ then $ t \rho_a \geq c_1 A^{1/3} (N \rho_a)^{-1}$, for any $A \geq 1$. 
\eel
\proof For the first part, assume that $\rho_b \geq \rho_a$ and by symmetry that $a \leq N/2$. Then,
\begin{align}
t^2 + t \rho_b & \leq t^2 + t C|2 - | \gamma_b | |^{1/2} \leq t^2 + Ct ( |2-|\gamma_a| )^{1/2} +  | \gamma_a - \gamma_b|^{1/2} ) \notag\\
&\leq C(t^2 + t |2 - |\gamma_a ||^{1/2} ) + \frac{t \rho_b}{2}
\end{align}
with the last inequality following from Cauchy-Schwarz. This completes the proof of the first claim. For the second part we can assume that $\rho_b \geq \rho_a$. From \eqref{eqn:b1} we see that there is a constant $B_1 >0$ so that if $\rho_b \geq B_1 t$ then $\rho_a \asymp \rho_b$. So assume that $\rho_b \leq B_1 t$. Then $t \geq c A^{1/3} N^{-1/3}$, and since $\rho_a \geq N^{-1/3}$ we see that $t \rho_a^2 \geq c A^{1/3} N^{-1}$. \qed

\bep \label{prop:bu-dif}
Assume that $a$ and $b$ are indices and $\rho_a \leq \rho_b$. Assume $t^2+t \rho_a \geq (N \rho_a)^{-1}$.  Then,
\beq 
| \bu_a (t) - \bu_b (t) | \lesssim \varphi^{1/100} \frac{ | \gamma_a - \gamma_b|}{N t(t+ \rho_a)^2}
\eeq
\eep
\proof We have the two estimates,
\begin{align}
 & \left| ( \cosh(t/2) \gamma_a - \gamma_j)^2 - ( \cosh(t/2) \gamma_a - \gamma_j)^2 \right| \notag\\
\leq & C| \gamma_a - \gamma_b | \left( | \cosh(t/2) \gamma_b - \gamma_j| + | \cosh(t/2) \gamma_a - \gamma_j | \right)
\end{align}
and
\beq
\sinh(t/2)^2 | \rho_a^2 - \rho_b^2| \leq C t^2 | \gamma_a - \gamma_b |.
\eeq
From the above two estimates, \eqref{eqn:buk-formula} and Lemma \ref{lem:buk-denom} we see that
\begin{align}
&| \bu_a (t) - \bu_b (t) | \leq \frac{ C \varphi^{1/100} | \gamma_a - \gamma_b| t}{N} \sum_{j=1}^N \frac{ (t^2 + | \gamma_a - \gamma_j| + | \gamma_b - \gamma_j| )N^{-2/3} ( \hat{j})^{-1/3}}{(( \gamma_a - \gamma_j)^2 + (t^2 +t \rho_a)^2)(( \gamma_b - \gamma_j)^2 + ( t^2 + t \rho_b)^2)}  \notag\\
& \leq C \varphi^{1/100} \frac{ |\gamma_a - \gamma_b|}{ (t^2 + t \rho_a )} \left[ \frac{t}{N} \sum_{j=1}^N \frac{1}{ ( \gamma_a - \gamma_j)^2 + (t^2 + t \rho_a)^2} \frac{1}{ N^{2/3} \hat{j}^{1/3}} \right].
\end{align}
The term in the square brackets was estimated by $\O ( N^{-1} (t + \rho_a)^{-1})$ in the Proof of Lemma \ref{lem:buk-bound}, completing the proof. \qed

\bec
Assume $\rho_a \leq \rho_b$. Then, 
\beq \label{eqn:bua-dif}
| \bu_a (t) - \bu_b (t) | \lesssim \varphi^{1/100} \frac{ | \gamma_a - \gamma_b|}{N t(t+ \rho_a)^2}
\eeq
\eec
\proof We may assume $a \neq b$. First if 
\beq
| \gamma_a - \gamma_b | \leq t^2 + t \rho_a
\eeq
holds, then since $| \gamma_a - \gamma_b| \geq |\gamma_a - \gamma_{a+1} | \asymp (N \rho_a)^{-1}$ we get the estimate by Proposition \ref{prop:bu-dif}. 
 Otherwise, assuming $| \gamma_a - \gamma_b | \geq t^2 + t \rho_a$ we see that,
\beq
|\gamma_a - \gamma_b | \gtrsim t^2 + t \rho_a + \frac{1}{N \rho_a}.
\eeq
Therefore, the RHS of \eqref{eqn:bua-dif} is larger than
\beq
\frac{ | \gamma_a - \gamma_b|}{N t(t+ \rho_a)^2} \geq c \left( \frac{1}{N(t + \rho_a)} + \frac{1}{N^2 t \rho_a (t + \rho_a)^2} \right)
\eeq
which is larger than the unconditional $\ell^\infty$ bound \eqref{eqn:buk-bound-2}. \qed

\subsection{Evolution}

We promote the function $\bu_k (t)$ to a function on the upper half plane $\hh_+$ given by,
\beq
\bu (z, t) := \frac{1}{N \Im[ \msc (z_t) ] } \Im \left[ \sum_{j=1}^N \frac{ f_j}{ \gamma_j - z_t} \right], \qquad f_j := \lambda_j (0) - \mu_j (0).
\eeq
Above, $z_t$ denotes the characteristic \eqref{eqn:char-def-1} starting at $z_0 = z$. We define $\bu (x, t) := \bu ( x + \i 0, t)$. With this notation we note that $\bu_k (t) = \bu ( \gamma_k, t)$.  We now compute the time derivative of $\bu (x, t)$. 

\bep \label{prop:deltbu}
For $t >0$ and $x \in (-2, 2)$ we have
\beq \label{eqn:deltbut}
\del_t \bu (x, t) = \int_{-2}^2 \frac{ \bu (y, t) - \bu (x, t)}{ (x-y)^2} \rhosc (y) \d y
\eeq
\eep
\remark The integral on the RHS of \eqref{eqn:deltbut} is interpreted in the principal value sense. \qed

We first compute the RHS of \eqref{eqn:deltbut}.
\bel \label{lem:deltbu-1}
For $x \in (-2, 2)$ and $t>0$ we have,
\beq \label{eqn:deltbu-a}
\int_{-2}^2 \frac{ \bu (y, t) - \bu (x, t)}{ (y-x)^2} \rhosc (y) \d y = \frac{ \bu (x, t)}{2} + \frac{1}{ 2 \Im [ \msc (x_t) ] } \sum_{j=1}^N \frac{ f_j}{N} \Im\left[ \frac{ \cosh (t/2) \sqrt{x^2-4} + \sinh(t/2) x }{ ( \gamma_j - x_t)^2} \right]
\eeq
\eel
\proof For simplicity of notation we let $\bu (x) := \bu (x, t)$. We may compute the integral in question via,
\begin{align} \label{eqn:deltbu-1}
\int_{-2}^2 \frac{ \bu (y) - \bu (x) }{ (y-x)^2} \rhosc (y) \d y = \lim_{ \eta \dto 0} \int_{-2}^2 \frac{ \bu (y) - \bu (x) }{ (y-x)^2} \rhosc (y) \d y = \lim_{ \eta \dto 0} \Im \left[ \frac{1}{ \eta} \int_{-2}^2 \frac{ \bu (y) - \bu (x) }{ y - x - \i \eta } \rhosc (y) \d y \right].
\end{align}
The integral on the RHS consists of two terms, with $\bu (x)$ and one with $\bu (y)$. For the former we have by the fact that
\begin{align}
\int_{-2}^2 \frac{1}{y -x - \i \eta } \rhosc (y) \d y = \msc ( x+ \i \eta) = \msc (x) + \i \eta \msc' (x) + \O ( \eta^2)
\end{align}
that
\beq \label{eqn:deltbu-2}
\Im\left[ \frac{ \bu (x)}{ \eta} \int_{-2}^2 \frac{ \rhosc (y) \d y }{ y - x - \i \eta} \right] = \frac{ \bu (x)}{2 \eta} \sqrt{ 4 -x^2} - \frac{ \bu(x)}{2} + \O ( \eta).
\eeq
Above and in the remainder of the proof the big $\O$ notation is uniform in $\eta$ but not necessarily in $N$ or $t>0$. This is immaterial as we are taking $\eta \to 0$ in the proof. 

Let $ z= x + \i \eta$ for notational simplicity. We have by the definition of $\bu (y)$,
\begin{align}
 \int_{-2}^2 \frac{ \bu (y) \rhosc (y) \d y}{ y-z} = & \sum_{j=1}^N \frac{ f_j}{N} \int_{-2}^2 \frac{ \rhosc (y)}{ \Im[ \msc (y_t) ]} \Im\left[ \frac{1}{ \gamma_j - y_t} \right] \frac{1}{y-z} \d y \notag \\
= & \sum_{j=1}^N \frac{ f_j \e^{t/2}}{N \pi} \int_{-2}^2 \Im \left[ \frac{1}{ \gamma_j -y_t} \right] \frac{1}{y-z} \d y \notag \\
= & \sum_{j=1}^N \frac{ f_j \e^{t/2}}{N \pi} \int_{\rr} \Im \left[ \frac{1}{ \gamma_j -y_t} \right] \frac{1}{y-z} \d y \notag \\
= & \sum_{j=1}^N \frac{ f_j \e^{t/2}}{N} \frac{1}{ \gamma_j - z_t}
\end{align}
In the second line we used that $\Im[ \msc (y_t) ]  = \e^{-t/2} \Im[ \msc (y+ \i 0) ] = \e^{-t/2} \pi \rhosc (y)$ for $y \in (-2, 2)$. The third line follows from the fact that $\Im[y_t] = 0$ for $y \notin (-2, 2)$. The fourth line is a consequence of the residue theorem. Recall now that $z_t$ denotes the characteristic starting at $z_0 = z = x + \i \eta$. We further compute,
\begin{align}
\frac{ 1}{\gamma_j -z_t} &= \frac{1}{ \gamma_j -x_t}+ \frac{1}{ ( \gamma_j -x_t)^2} (z_t -x_t) + \O ( \eta^2) \notag \\
&= \frac{1}{ \gamma_j -x_t}+ \frac{\i \eta}{ ( \gamma_j -x_t)^2} \left( \cosh(t/2) + \sinh(t/2) \frac{ x}{ \sqrt{x^2-4}} \right) + \O ( \eta^2).
\end{align}
Therefore,
\begin{align} \label{eqn:deltbu-3}
& \Im\left[ \frac{1}{ \eta} \int_{-2}^2 \frac{ \bu (y) \rhosc (y) \d y }{ y-z} \right] \notag\\
= & \frac{1}{ \eta} \sum_{j=1}^N \frac{ f_j \e^{t/2}}{N} \Im\left[ \frac{1}{ \gamma_j - x_t } \right] + \sum_{j=1}^N \frac{ f_j \e^{t/2}}{N} \Im\left[ \i \frac{ \cosh(t/2) + \sinh(t/2) \frac{x}{\sqrt{x^2-4}}}{(\gamma_j -x_t)^2} \right] + \O ( \eta).
\end{align}
By \eqref{eqn:deltbu-1}, it remains to show that the first term on the RHS of \eqref{eqn:deltbu-3} cancels with the first term on the RHS of \eqref{eqn:deltbu-2} and that the second term on the RHS of \eqref{eqn:deltbu-3} gives the second term on the RHS of \eqref{eqn:deltbu-a}. First,
\beq
\sum_{j=1}^N \frac{ f_j \e^{t/2}}{N} \Im\left[ \frac{1}{ \gamma_j -x_t} \right] = \frac{ \sqrt{4-x^2}}{2} \bu (x)
\eeq
since $\frac{ \sqrt{4-x^2}}{2} = \e^{t/2} \Im[ \msc (x_t)]$. Second, since $ \sqrt{x^2-4} = \i \sqrt{4-x^2}$,
\begin{align}
 & \sum_{j=1}^N \frac{ f_j \e^{t/2}}{N} \Im\left[ \i \frac{ \cosh(t/2) + \sinh(t/2) \frac{x}{\sqrt{x^2-4}}}{(\gamma_j -x_t)^2} \right] \notag\\
= & \sum_{j=1}^N \frac{ f_j}{ N \e^{-t/2} \sqrt{4-x^2}} \Im \left[ \frac{ \cosh(t/2) \sqrt{x^2-4} + \sinh(t/2) x }{ (\gamma_j -x_t)^2} \right] \notag\\
= & \frac{1}{ 2 \Im[ \msc (x_t) ]} \sum_{j=1}^N \frac{f_j}{N} \Im \left[ \frac{ \cosh(t/2) \sqrt{x^2-4} + \sinh(t/2) x }{ (\gamma_j -x_t)^2} \right].
\end{align}
This completes the proof. \qed

\vspace{5 pt}

\noindent{\bf Proof of Proposition \ref{prop:deltbu}}. In light of Lemma \ref{lem:deltbu-1} it remains to compute the LHS of \eqref{eqn:deltbut}. By direct computation,
\begin{align}
\del_t \bu (x, t) &= \del_t \left( \frac{1}{N \Im[ \msc (x_t) ] } \sum_{j=1}^N \Im \left[ \frac{1}{ \gamma_j -x_t} \right] \right) \notag\\
&= - \frac{ \del_t \Im[ \msc (x_t)]}{ \Im[ \msc(x_t)]} \bu(x, t) + \frac{1}{ N \Im[ \msc(x_t)]} \sum_{j=1}^N f_j \Im \left[ \frac{ \del_t x_t}{ ( \gamma_j -x_t)^2} \right] \notag\\
&= \frac{ \bu (x, t)}{2} + \frac{1}{2 N \Im[ \msc (x_t)] } \sum_{j=1}^N \Im \left[ \frac{ \sinh(t/2) x + \cosh(t/2) \sqrt{ x^2-4} }{ ( \gamma_j -x_t)^2} \right]
\end{align}
which completes the proof. \qed

Finally we record here a generalization of Proposition \ref{prop:bu-dif}. The proof is the same, so is omitted.
\bel
Let $E \in [\gamma_N, \gamma_1]$ and assume $t^2 +t \rhosc (E) \geq (N \rhosc(E))^{-1}$. Then,
\beq \label{eqn:bu-deriv}
| \del_y \bu (y, t) | \vert_{y = E} \lesssim \varphi^{1/100} \frac{1}{N t  (t+ \rhosc(E))^2} .
\eeq
Furthermore, for any index $a$ we have,
\beq \label{eqn:bua-dif-sup}
\sup_{ y \in [\gamma_a, \gamma_{a+1} ] } | \bu (y, t) - \bu_a (t) | \lesssim \varphi^{1/100} \frac{ (N \rho_a)^{-1}}{N t (t + \rho_a)^2}
\eeq
\eel

\section{Homogenization by maximum principle}

We introduce,
\beq
B(t, k) := N t \rho_k (t+\rho_k), \qquad Q_a (t, k) := N(t + \rho_k)[ Nt \rho_k (t+\rho_k)]^{1-a} = N(t+\rho_k ) B(t, k)^{1-a}
\eeq
for $a \in [0, 1]$. 
\bed
We say that property $\pP_a$ holds if for any $D>0$ there exists $C_1 \geq 1$ and $N_0$ so that for all $N \geq N_0$ we have,
\beq
\pp\left[ \forall t \in (0, 1], k \in [\![1, N ]\!] : N t \rho_k (t + \rho_k)  \geq \varphi^{C_1} , |\fu_k (t) - \bu_k (t) | \leq \frac{ \varphi^{C_1}}{Q_a (t, k) } \right] \geq 1 - N^{-D}
\eeq
\eed

\bet \label{thm:induction-step} There is a $c \in (0, 1)$ so that if Property $\pP_a$ holds then $\pP_{c a }$ holds.
\eet

\bel \label{lem:initial}
Property $\pP_1$ holds.
\eel
\proof From \cite[(2.18)]{bourgade2018extreme} we have
\beq \label{eqn:u-inf-bd}
| \fu _k (t) | \leq \frac{ \varphi^C}{N (t + \rho_k)} 
\eeq
with overwhelming probability. This, combined with \eqref{eqn:buk-bound} yields the claim. \qed

Property $\pP_a$ makes an assumption on $t$. The following shows that in practice this is not necessary, which will simplify estimation of various error terms. 
\bel \label{lem:uncond-Pa}
Suppose Property $\pP_a$ holds. Then for any $D>0$ there is a $C_1 >0$ so that with probability at least $1 - N^{-D}$ we have, for all $ t\in [0, 1]$ and $k \in [\![1, N ]\!]$
\beq
| \bu_k (t) - \fu_k (t) | \leq \varphi^{C_1} \left(  \frac{1}{Q_a (t, k)} + \frac{1}{ Q_0 (t, k)} \right)
\eeq
\eel
\proof If $N t \rho_k (t + \rho_k) \leq \varphi^{C_1}$ then the second term on the RHS of \eqref{eqn:buk-bound-2} is larger than $|\bu_k (t)| + | \fu_k (t) |$, and this term is $1 / Q_0 (t, k)$, up to factors of $\varphi$. \qed

The function $\fu_k (t)$ satisfies the equation,
\beq
\del_t \fu_k (t) = (\bB_t \fu_k (t) )_k, \qquad ( \bB_t v)_k := \frac{1}{N} \sum_{i \neq k} \frac{ v_i -v_k}{ (x_i (t) - x_k(t) )^2} 
\eeq
We denote its semigroup by $\UB (u, t)$ i.e., 
\beq
\del_t \UB (u, t)_{ij} = \sum_k (\sB_t)_{ik} \UB (u, t)_{kj}, \qquad \UB (u, u)_{ij} = \delta_{ij} .
\eeq
We require the short-range operator defined by,
\beq
( \sS u)_i := \frac{1}{N} \sum_{j : |j-i| < \ell} \frac{u_j - u_i}{ (x_i -x_j)^2}
\eeq
where $\ell$ is a parameter. We will always assume $\varphi^{100} \leq \ell \leq \varphi^{-100} N$. We define its semigroup by $\US$. We also introduce the long-range operator by
\beq
\bB = \sS + \lL .
\eeq
The following estimate will be useful. It is proven at the end of Section \ref{sec:fs}, in Subsection \ref{sec:fs-proof}

\bet \label{thm:fs} With overwhelming probability the following holds. For any index $a$,  define
\beq
T_a := \varphi^{-12} \frac{ \ell}{N^{1/3} ( \ell + a)^{2/3}}.
\eeq
For all $j$ such that $|j-a| > \varphi^5 \ell$ and $|s-t| \leq T_a$ we have
\beq
\US (s, t)_{ja}  + \US (s, t)_{aj} \leq \e^{ -c \varphi}.
\eeq
\eet

\subsection{Induction argument}

The purpose of this section is to give the proof of Theorem \ref{thm:induction-step}. 

\vspace{5 pt}

\noindent{\bf Notational convention.} \emph{Throughout the proof we will make the following simplifying assumption. We will always assume that $k < \frac{N}{200}$. Bourgade's work \cite{bourgade2018extreme} proves Theorem \ref{thm:induction-step} for the case that $\alpha N \leq k \leq (1- \alpha )N$ for any fixed $ \alpha >0$. It therefore suffices to consider $k$ near $1$ or $N$. By symmetry it suffices to consider the case that $k$ is near $1$. Therefore, we will assume that $k < \frac{N}{200}$ throughout the rest of the paper (with heavier notation  what is written here could accommodate the case of bulk $k$, but there does not seem to be much gained by introducing further complexity). }

\emph{Furthermore, various estimates involve another index $j$ which will always be close to $k$: $|j-k| \leq \varphi^{-1} N$. We therefore always assume that $j < \frac{N}{100}$. Furthermore, various estimates in the paper will involve sums over $i \in [\![1, N]\!]$. In general, the contribution from terms $i > N/2$ will always be smaller than the contribution from indices $i \leq N/2$, because the objects being summed over will usually have contributions like $|\gamma_i - \gamma_j|$ in the denominator which will be smaller when $i$ is close to $j$ (for example, many of the estimates in Section \ref{sec:flow}). When it is safe to do so, we will therefore tacitly ignore the contribution from the indices $i \geq N/2$ as it is dominated by the contribution from terms $i \leq N/2$. }

\vspace{5 pt}

The proof requires the construction of an approximating flow which involves several parameters (in addition to the choice of $\ell$ above in the definition of $\sS$). These choices depend on the index $k$ and time $t$ at which we are trying to prove the desired estimate. We therefore fix $(t, k)$ for the remainder of the section until the end of the proof of Theorem \ref{thm:induction-step}. However, as can be observed from the proof, the events on which we prove the desired estimates hold simultaneously for all $(t, k)$ satisfying $B(t, k) \geq \varphi^{C_1}$ for some large $C_1 >0$.

We define an intermediate flow $m_j (s)$ for $s \in [u, t]$. Here, $u$ is a time satisfying,
\beq
\frac{t}{2} < u < t.
\eeq 
We construct the initial data at time $u$ to equal to $\fu$ nearby index $k$ and equals to $\bu_k$ for indices away from $k$. More precisely, choose a parameter $R >0$ s.t. $R \geq \varphi^{100} \ell$ as a localization scale, then define, for $b\in [\![R, 2R]\!]$
\begin{align}\label{def: average operator}
(\Flat_b(s) w)_j = \begin{cases}
     w_j,  |j-k|<b\\ \bu_k(s), |j-k|\geq b
\end{cases},
\quad (\Av(s) w)_j = \frac{1}{R}\sum_{R\leq b \leq 2R}(\Flat_b(s) w)_j.
\end{align}
Note our average operator $\Av(s)$ is time dependent. Define the intermediate flow $m_j (s)$ for $s \in [u, t]$ by,
\begin{align}\label{def: intermediate m}
\begin{cases}
\partial_t m_j(s) = (\sS_s m(s))_j + \del_s \bu_k (s)\\
m_j(u) = (\Av(u) \fu (u))_j.
\end{cases}
\end{align}
Note that in the first equation, the inhomogeneous term $ \del_s \bu_k (s)$ does not depend on the index $j$ but is time dependent. The point is that by \eqref{eqn:deltbut}, we hope that it is a good approximation of $( \lL \fu)_j$ for $j$ near $k$. The proof the following is given in Section \ref{sec:|u-m|-proof}.

\bel \label{lem:|u-m|}
 Assume Property $(\pP_a)$ holds. The following event holds with overwhelming probability. Suppose $\ell \geq \varphi^{20}  N|u-t|\kappa_{k\vee \ell}$. Assume also that for all $j : |j-k| \leq \varphi^{20} \ell$ we have $(N \rho_j)^{-1} \leq t^2 +t \rho_j$. 
 Then we have
\begin{align}\label{eqn: |u-m|}
\abs{\fu_k(t) - m_k(t)} \leq \varphi^{C_2} |u-t|\sup_{j: |j-k|\leq \varphi^{5} \ell }\pare{\frac{N\kappa_{j\vee \ell}}{Q_0(t,j\wedge k)} + \frac{N\kappa_{j\vee \ell}}{\ell Q_a(t,j)}}.
\end{align}
\eel

\subsubsection{Maximum principle}

The second step is to compare the intermediate flow $m(t)$ with $\bu_k (t)$ by a maximum principle.

Define, for $s\in [u,t]$, $i\in [\![1,N]\!]$, 
\beq
M_i(s):=m_i(s)-\bu_k (s) .
\eeq
Then, by the definition of $m(s)$ in \eqref{def: intermediate m} and the fact that constant vectors are in the kernel of $\sS$ we have,
\beq \label{eqn:M-evolution}
\partial_sM_i(s) = (\sS m(s))_i  = (\sS(m(s)-\bu_k(s)))_i = \pare{\sS M(s)}_i.
\eeq

We wish to bound $|M_i(t)|$ by a maximum principle. We first bound $\max_i M_i (t)$ above. Since $M(s)$ satisfies \eqref{eqn:M-evolution}, the maximum $\max_i M_i(s)$ is non-increasing. So, without loss of generality, we assume $\max_i M_i(s)>0$ for all $s\in [u,t]$.

Let
\begin{align}\label{def: maximum index set J}
J(s) = \brac{j: M_j(s) = \max_i M_i(s)}.
\end{align}
be the set of maximum indices at time $s$. Pick an index $j\in J(s)$ and choose a parameter $\eta = \eta(j)>0$. We will always assume that
\beq \label{eqn:eta-j-assump}
\eta(j) \geq \frac{ \varphi}{N \rho_j} .
\eeq
We remark here that the choice $\eta(j)$ can depend on the index $j$ at which the maximum is attained, but the other parameters such as $R, u$ and $\ell$ cannot be chosen to depend on $j$, as  $M_i(s)$ is itself a function of these parameters. Set $z=\gamma_j+\i\eta$. We have, with overwhelming probability, 
\begin{align} \label{eqn:M-max-principle-1}
\partial_s \pare{\max_iM_i(s)} &= (\sS M(s))_j =\frac{1}{N}\sum_{i:|i-j|<\ell, i \neq j}\frac{(m_i(s)-\bu_k(s))-(m_j(s)-\bu_k(s))}{(x_i(s)-x_j(s))^2}\notag\\
&\leq \frac{1}{N\eta}\im\sum_{i:|i-j|<\ell}\frac{(m_i(s)-\bu _k(s))-(m_j(s)-\bu_k(s))}{x_i(s)-(x_j(s) + \i \eta)}\notag \\
&\lesssim  \frac{1}{N\eta}\im\sum_{i:|i-j|<\ell}\frac{(m_i(s)-\bu_k(s))-(m_j(s)-\bu_k(s))}{x_i(s)-z } \notag\\
&= -\frac{1}{N\eta}\pare{\im\sum_{i:|j-i|<\ell}\frac{1}{x_i(s)-z}}M_j(s) + \frac{1}{N\eta}\im\sum_{i:|i-j|<\ell}\frac{m_i(s)-\bu_k(s)}{x_i(s)-z}.
\end{align}
In the first inequality, we used that $M_j(s) \geq M_i(s)$ and
$
\frac{1}{(x_i(s)-x_j(s))^2}>\frac{1}{\eta}\im \frac{1}{x_i(s)-z}.
$
In the second inequality we used the fact that \eqref{eqn:eta-j-assump} implies that $| x_i (s) - z | \asymp |x_i(s) - (x_j (s) + \i \eta) |$ whenever the rigidity estimates \eqref{eqn:rig} hold.

The following lemma determines the size of the prefactor in the equation \eqref{eqn:M-max-principle-1}. 
\bel \label{lem:im-m-ell}
Assume that $\eta(j) \geq \varphi (N \rho_j)^{-1}$ and that
\beq
\ell \geq \varphi N \eta \sqrt{ \eta + \kappa_j}. 
\eeq
Then,
\beq
\frac{1}{N} \sum_{ i : |i-j| < \ell} \frac{\eta}{ (x_i - \gamma_j)^2 + \eta^2} \asymp \sqrt{ \kappa_j + \eta}
\eeq
\eel
The proof is deferred to Appendix \ref{sec:im-m-ell-proof}. From Lemma \ref{lem:im-m-ell} and \eqref{eqn:M-max-principle-1} we see that an important role will be played by the error term,
\beq
X_j (s) := \frac{ 1}{N \sqrt{ \kappa_j + \eta}} \Im \sum_{ i :  |i-j| < \ell} \frac{ m_i (s) - \bu_k (s) }{x_i - z}
\eeq

The following collects the maximum principle argument.
\bep
Suppose that $R \geq \varphi^{100} \ell$. Assume that
\beq
|u - t | \leq \varphi^{-100} \frac{ \ell}{ N^{1/3} (R+ k)^{2/3}}
\eeq
and that for every choice of $\eta = \eta(j)$ for $|j-k| \leq 3 R$ we have
\beq
\eta \geq \frac{ \varphi}{N \rho_j} , \qquad \ell \geq \varphi^{10} N \eta \sqrt{ \kappa_{j} + \eta }
\eeq
as well as 
\beq \label{eqn:param-assump-1}
|u - t | \geq \varphi \frac{\eta}{ \sqrt{\kappa_j + \eta} }.
\eeq
Then,
\beq \label{eqn: |m-u_bar|}
|m_k(t)-\bu_k(t)|  \leq \e^{ - c \varphi} + \sup_{j : |j-k| \leq 3 R , s \in [u, t] } |X_j (s) |
\eeq
\eep
\proof We first claim that if the maximum of $M_i (s)$ occurs at any $j$ such that $|j-k| \geq 3 R$ then $\max_i M_i (s) \leq \e^{ -c \varphi}$. Indeed, this follows from the fact that $M_i (u) = 0$ for $|i-k| > 2 R$ and the fact that by Theorem \ref{thm:fs}
\beq
\US (u, s)_{ji} \leq \e^{ - c \varphi}
\eeq
for any $|i-k| \leq 2 R$ and $|j-k| \geq 3 R \geq 2 R + \varphi^{10} \ell$ since $|u-t| \leq \varphi^{-50} \frac{\ell}{ N^{1/3} ( \ell + i )^{2/3} }$ by our assumption on $|u-t|$. Since the maximum of $M_i (s)$ is non-increasing in time, we may dispense with this case and assume that for all $s \in [u, t]$ the index $j$ at which the maximum occurs satisfies $|j-k| \leq 3 R$. 

Then we have by integrating \eqref{eqn:M-max-principle-1}, (note that $j, \eta$ depend on $s$), 
\begin{align}
& \max_i M_i(t) \lesssim \e^{-c \int_u^t \frac{\sqrt{\kappa_j+\eta}}{\eta} \d v}\max_i M_i(u)+\int_u^t\e^{-c\int_s^t\frac{\sqrt{\kappa_j+\eta}}{\eta} \d v} \frac{ \sqrt{\kappa_j+\eta}}{\eta} X_j(s) \d s \notag\\
&\leq \e^{-c \varphi} + \int_u^t\e^{-\int_s^t\frac{\sqrt{\kappa_j+\eta}}{\eta}\d v}\frac{\sqrt{\kappa_j+\eta}}{\eta}\d s \times \sup_{s\in [u,t]}|X_j(s)| \notag \\
&\lesssim \e^{-c \varphi}+ \int_u^t \frac{\d}{\d s}\left(  \e^{-c\int_s^t \frac{\sqrt{\kappa_j+\eta}}{\eta}\d v} \right) \d s \times \sup_{s\in [u,t]}|X_j (s)| \notag \\
&\leq \e^{-c \varphi} + \max_{j:|j-k|<3 R }\sup_{s\in [u,t]} \abs{X_j(s)} .
\end{align}
This bounds $m_k(t) - \bu_k (t) =M_k (t) \leq \max_i M_i (t)$ above by the RHS of \eqref{eqn: |m-u_bar|}. By a similar argument we bound $\min_i M_i (t)$ by the same quantity, completing the proof. \qed

It remains to estimate $X_j(s)$. 
The derivation of the dummy index $j$ and designated index $k$ above gives extra complexity to estimate the RHS of \eqref{eqn: |m-u_bar|}, in contrast to the bulk case in \cite{bourgade2018extreme} where $j\asymp k$ always holds. As $j\asymp k$ if $R\ll \hat{k}$, it is natural to divide into bulk and edge cases according to the relative size of $R$ and $k$, and discuss differently on two cases. We recall our assumption that,
\beq
B(t, k) := N t (t + \rho_k) \rho_k \geq \varphi^{C_1} =: \varphi^{A_1 /a }
\eeq
for some large $A_1 >0$ to be determined. We will be choosing
\beq
R = [ Nt (t + \rho_k ) \rho_k]^{a \mfc_1}
\eeq
for some $\mfc_1 \in (0, 1)$ to be determined.  The arguments are somewhat different depending on whether or not $R \ll k$. To facilitate this, let us introduce the following sets of indices.
\begin{align}
\Da &:= \{ k : B(t, k) \geq \varphi^{A_1/a} \} \notag \\
\Db ( \cb ) &:= \{ k \in \Da : B(t, k)^{ \cb a} \leq \varphi^{-100} k \} \notag \\
\De ( \ce ) &:= \{ k \in \Da : B(t, k)^{ \ce a/100 } \geq k \} 
\end{align}

\bel
Assume that $\cb, \ce \in (10^{-9}, 0.5 )$ and that $\cb \leq \ce/200$. Then if $A_1 \geq 10^{20}$ we have 
\beq
\Da = \Db \cup \De .
\eeq
If $k \in \Db$ then 
\beq
\sqrt{ \frac{ k^{(\mfc_b a)^{-1}}}{N \rho_k } } \geq c \varphi^{50/ (\mfc_b a ) } t
\eeq
and if $k \in \De$ then $t > \rho_k$ and 
\beq
t \geq \sqrt{ \frac{ k^{( \mfc_e a)^{-1}}}{N \rho_k}}
\eeq
Furthermore, 
\beq
t \geq  c \rho_{B(t, k)^{\mfc_e a }}
\eeq
\eel
\proof The claim about the union is obvious. If $k \in \Db$ and if $t \geq \rho_k$ then the definition of $\Db$ gives
\beq
Nt^2 \rho_k \asymp B(t, k) \leq \varphi^{-100  / ( \mfc_b a)} k^{ \frac{1}{ \mfc_b a}}
\eeq
which gives the desired inequality. On the other hand, if $t \leq \rho_k$ then since $k \geq \varphi^{100}$,
\beq
\sqrt{ \frac{ k^{(\mfc_b a)^{-1}}}{N \rho_k } }  \geq \sqrt{ \frac{ k}{N \rho_k} } k^{1/(4 \mfc_b a)} \geq c \rho_k \varphi^{50 / (\mfc_b a) }.
\eeq
Now assume that $k \in \De$. If $t \leq \rho_k$ then,
\beq
B(t, k)^{ \mfc_e a /2} \leq \varphi^{-1} B(t, k) \leq \varphi^{-1} C N \rho_k^3 \leq \varphi^{-1} C k
\eeq
and so $t > \rho_k$. Then,
\beq
N t^2 \rho_k \asymp B(t, k) \geq k^{\frac{2}{ \mfc_e a } } \geq k^{ \frac{1}{ \mfc_e a } } .
\eeq
For the last claim, since we know $t > \rho_k$
\begin{align}
t \asymp \sqrt{ \frac{ B(t, k)}{N \rho_k} } \asymp \frac{ B(t, k)^{1/2}}{ N^{1/3} k^{1/6}} \geq \frac{ B(t, k)^{1/2-\mfc_e a/6 }}{N^{1/3}} \geq \frac{B(t, k)^{\mfc_e a /3}}{N^{1/3}} \asymp \rho_{ (B (t, k))^{\mfc_e a } } 
\end{align}
\qed

In the remainder of the proof if $k \in \Db$ we take $R = B(t, k)^{a \cb}$ and if $k \in \De$ we take $R = B(t, k)^{a \ce}$. Note that if $k \in \De$ then
\beq
k \leq R^{1/100}, 
\eeq
and if $k \in \Db$ then
\beq
R \leq \varphi^{-100} k.
\eeq
In particular if $|j-k| \leq 3 R$ then $j \asymp k$.

\bel \label{lem:index-properties}
Suppose that $B(t, k) \geq \varphi^{C_1}$. If $k \in \Db ( \cb)$ and $R = B(t, k)^{a \cb}$ and if $|j -k| \lesssim  R$ then $B(t, j) \gtrsim \varphi^{C_1}$. If $k \in \De ( \ce)$, then $t \gtrsim \varphi^{C_1/2} N^{-1/3}$ and so $B(t, j) \gtrsim \varphi^{C_1}$ for all $j$. 
\eel
\proof We discussed the $\Db$ case above.  If $k  \in \De $ then it follows from previous lemma that $t \gtrsim \rho_k$ and so $B(t, k) \geq \varphi^{A_1/a}$ implies that $t^2 \gtrsim \varphi^{C_1} (N \rho_k)^{-1} \gtrsim \varphi^{C_1} N^{-2/3}$ and so $t \gtrsim \varphi^{C_1/2} N^{-1/3}$. Therefore $B(t, j) \geq N \rho_j t^2 \gtrsim \varphi^{C_1}$.  \qed 

The proof of the following is given in Section \ref{sec:ave-proof}
\bel\label{lem: ave interpolated dyn}
Assume that Property $\pP_a$ holds, let $j$ be an index satisfying $|j-k| \leq 3 R$. Let $z = \gamma_j + \i \eta$. Assume that,
\beq \label{eqn:param-assump-2}
\frac{\varphi}{N \rho_j} \leq \eta \leq \varphi^{-1} (t^2 + t \rho_j), \quad R \varphi^{-100} \geq \ell \geq \varphi^{100} N |u-t|  \kappa_{R \vee k} , \quad \ell \geq \varphi N \eta \sqrt{ \kappa_{k \vee R } + \eta }
\eeq
Then we have that the term,
\beq
\frac{1}{N \sqrt{ \kappa_j + \eta} } \Im \sum_{ i : |i-j| < \ell} \frac{ m_i (s) - \bu_k (s) }{ x_i (s) -z } 
\eeq
is bounded by the following error terms. If $k \in \Db$ by (up an overall $\varphi^C$ factor),
\begin{align}
& \frac{1}{Q_a (t, k)} \left( \frac{1}{N \eta \rho_k} + \frac{N \eta \rho_k}{\ell} + \frac{\ell}{R} + \frac{ |u-t| N \kappa_k}{\ell} \right) +  \frac{1}{Q_0 (t, k) } \left(R 
+ \frac{\ell^2}{N \eta \rho_k } \right)
\end{align}
and if $k \in \De$ by  (up an overall $\varphi^C$ factor),
\begin{align}
 & \frac{1}{Q_a (t, j) } \left( \1_{ j > \varphi^{12} \ell} \frac{ |u-t| N \kappa_j}{\ell} + \frac{1}{N \eta \sqrt{\kappa_j + \eta}} + \frac{N \eta \kappa_{j \vee \ell}}{\ell \sqrt{ \kappa_j+\eta}} \right) + \frac{\ell}{R} \frac{1}{Q_a (t, R)} \notag\\
+ & \frac{R + |u-t| N \kappa_R k^{1/3}}{Q_0 (t, k) } + \frac{1}{Q_0 (t, j \wedge k) } \left( \frac{R \rho_{j \vee \ell}}{\sqrt{\kappa_j + \eta}} + \frac{ \ell^2}{N \eta \sqrt{\kappa_j + \eta}} \right)  \notag\\
+ & \1_{j \leq \varphi^{12} \ell} \frac{ |u-t| N \kappa_\ell}{Q_a (t, 1) \ell}
\end{align}
\eel

\subsubsection{Proof of Theorem \ref{thm:induction-step}}

We recall that we approximate $\fu_k$ by $m_k$ using \eqref{eqn: |u-m|} and then use \eqref{eqn: |m-u_bar|} to approximate $m_k$ by $\bu_k$, and bound error that results in \eqref{eqn: |m-u_bar|} with Lemma \ref{lem: ave interpolated dyn}.

The argument is a little different in the bulk and edge cases, so we split them accordingly, beginning with the bulk. 

\vspace{5 pt}
\noindent{\bf Case $k \in \Db$}. In this case we obtain from \eqref{eqn: |u-m|}, \eqref{eqn: |m-u_bar|} and Lemma \ref{lem: ave interpolated dyn} that,
\begin{align}
& |\fu_k (t) - \bu_k (t) | \leq  \frac{\varphi^C}{Q_0 (t, k) } \left( N \kappa_k |u-t| + R + \frac{ \ell^2}{N \eta \rho_k} \right) \notag\\
+ & \frac{\varphi^C}{Q_a (t, k) } \left( \frac{ |u-t| N \kappa_k}{\ell} + \frac{1}{N \eta \rho_k}  + \frac{ \ell}{R} \right).
\end{align}
The assumptions on the parameters we require in \eqref{eqn:param-assump-1} and \eqref{eqn:param-assump-2} become,
\begin{align} \label{eqn:bulk-ass}
|u-t| \rho_k \geq \varphi \eta, \quad \frac{ \varphi}{N \rho_k} \leq \eta \leq \varphi^{-1} (t^2 + t\rho_k ) , \quad & R \varphi^{-100} \geq \ell \geq \varphi^{100} N |u-t| \kappa_k , \quad \ell \geq \varphi N \eta \rho_k \notag\\
|u-t| & \leq \frac{t}{2}.
\end{align}
Recall that $R = B(t, k)^{\cb a} \geq \varphi^{A_1 \cb}$. Choose,
\beq
|u -t| N \rho_k^2 = R^{1/3} , \quad \eta N \rho_k = \frac{R^{1/3}}{\varphi^{10}}, \quad \ell = R^{2/3}.
\eeq
It is easy to see that if $A_1$ is large enough (depending on $\cb$) then the  first line of the assumptions \eqref{eqn:bulk-ass} hold. For the last one, if $t < \rho_k$, if $A_1$ is large enough and $\cb < \frac{1}{2}$,
\beq
|u-t| = \frac{R^{1/3}}{N \rho_k^2}  \leq \varphi^{-1} \frac{ N t \rho_k^2}{N \rho_k^2} \leq \varphi^{-1} t
\eeq
and if $ \rho_k < t$ then
\beq
|u -t | \leq \varphi^{-1} \frac{ (N t^2 \rho_k)^{1/2}}{N \rho_k^2} = \varphi^{-1} t \frac{1}{(N \rho_k^3)^{1/2} } \leq \varphi^{-1} t.
\eeq
 Then, with these choices we see that
\begin{align}
| \bu_k (t) - \fu_k (t) | \leq&  \varphi^C \left(  \frac{R}{Q_0 (t, k) } + \frac{1}{Q_a (t, k) R^{1/3} } \right) \notag\\
\lesssim&   \varphi^C \left( \frac{1}{Q_{\cb a} (t, k)} + \frac{1}{Q_{a(1-\cb/3)} (t, k) }  \right) \leq \frac{\varphi^C}{Q_a{(1-\cb/3)} (t, k)}
\end{align}
with the last inequality holding if $\cb < \frac{1}{2}$.

\vspace{5 pt}

\noindent{\bf Case $k \in \De$.} In this case we parameterize
\beq
\eta = \nu^{1/2}\sqrt{ \nu + \kappa_j}
\eeq
for some $\nu >0$. With this choice,
\beq
\sqrt{\kappa_j + \eta} \asymp \sqrt{ \kappa_j + \nu}, \quad \frac{\eta}{\sqrt{\kappa_j + \eta}} \asymp \nu^{1/2}.
\eeq
We make the choice
\beq
|u-t| = \varphi \nu^{1/2}.
\eeq
With these choices, we obtain from \eqref{eqn: |u-m|}, \eqref{eqn: |m-u_bar|} and Lemma \ref{lem: ave interpolated dyn} that
\begin{align}
& \varphi^{-C} | \fu_k (t) - \bu_k (t) |  \leq   \sup_{j | j-k|< 3 R}\left\{  \frac{N \nu^{1/2} \kappa_{j \vee \ell}}{\ell Q_a (t, j) } + \frac{ N \nu^{1/2} \kappa_{j \vee \ell}}{Q_0 (t, j \wedge k ) } \right\} \notag\\
+ & \frac{R + \nu^{1/2} N \kappa_R  k^{1/3}}{Q_0 (t, k) } + \sup_{ j : |j-k | \leq 3 R } \frac{1}{Q_0 (t, j \wedge k ) } \left( \frac{R \rho_{j \vee \ell}}{\sqrt{\kappa_j + \nu}} + \frac{ \ell^2}{N \nu^{1/2} ( \kappa_j + \nu)}  \right)  + \frac{ \ell}{R} \frac{1}{Q_a (t, R) } \notag\\
+ & \sup_{ j : |j-k| \leq 3 R }  \frac{1}{Q_a (t, j) } \left(  \frac{1}{N \nu^{1/2} ( \kappa_j + \nu )} + \frac{N \nu^{1/2} \kappa_{j \vee \ell}}{\ell} \right)  
+  \frac{  \nu^{1/2} N \kappa_\ell}{Q_a (t, 1) \ell}
\end{align}
under the assumptions that for all $j \leq 4 R$,
\beq
\frac{ \varphi}{N \rho_j} \leq \nu^{1/2} ( \nu^{1/2} + \rho_j ) \leq \varphi^{-1} (t^2 + t \rho_j ), \quad \ell \geq \varphi N \nu^{1/2} ( \nu^{1/2} + \rho_j) \sqrt{ \kappa_R + \nu }
\eeq
and
\beq
R \varphi^{-100} \geq \ell \geq  \varphi^{100} N \nu^{1/2} \kappa_R, \qquad \nu^{1/2} \leq \varphi^{-10} t .
\eeq
If we choose $\nu = R^\alpha / N^{2/3}$ and $\ell = R^\beta$ for $\alpha$ and $\beta$ satisfying
\beq
0 < \alpha < \frac{1}{5}, \qquad \frac{2}{3} + \frac{\alpha}{2} < \beta < 1
\eeq
then it is easy to see that all of the assumptions on the parameters are satisfied. The only non-trivial one is $\nu^{1/2} \leq \varphi^{-10} t$ (which implies $ \nu^{1/2} ( \nu^{1/2} + t \rho)j) \leq \varphi^{-1} (t^2 + t \rho_j)$) which we see from the calculation (using that $t > \rho_k$ for $k \in \De$)
\beq
\nu^{1/2} = \frac{ R^{\alpha/2}}{N^{1/3}} \lesssim \varphi^{-100} \frac{ (N t^2 \rho_k)^{1/3}}{N^{1/3}} = \varphi^{-100} t \left( \frac{\rho_k}{t} \right)^{1/3} \leq \varphi^{-100} t. 
\eeq
Using that $ t> \rho_{R}$ and that $Q_a (t, j)^{-1} \leq \frac{R^{1/100}}{Q_a (t, k)}$ because $k \leq R^{1/100}$ on $k \in \De$ we get,
\begin{align}
& | \bu_k (t) - \fu_k (t) | \leq  \varphi^C \left( \frac{R^3}{Q_0 (t, k)} 
+  \frac{1}{Q_a (t, k)} \left( R^{\alpha/2 +0.01 +2/3 - \beta} + R^{\beta-1} + R^{-3 \alpha /2 } \right) \right)
\end{align}
If we choose, for example $\beta= \frac{2}{3} + \frac{1}{10}$, $\alpha = \frac{1}{10}$ then we obtain, 
\beq
| \bu_k (t) - \fu_k (t) |  \leq \varphi^C \left( \frac{R^3}{Q_0 (t, k) } + \frac{ R^{-1/30}}{Q_a (t, k)} \leq \frac{1}{Q_{a(1-\ce/30)} (t, k) } \right)
\eeq

\vspace{5 pt}

\noindent{\bf Conclusion.} We can take $\ce = \frac{1}{10}$ and $\cb = \frac{1}{3000}$. Then, in either case we have
\beq
| \bu_k (t) - \fu_k (t) | \leq \frac{ \varphi^{C_1}}{Q_{a c_1 } (t, k) }
\eeq
for $c_1 = 1 - 10^{-5}$. 
This completes the proof. \qed

\vspace{5 pt}

\noindent{\bf Proof of Theorem \ref{thm:main-homog}} Given Theorem \ref{thm:induction-step} and Lemma \ref{lem:initial}, we see that Property $\pP_\eps$ holds for any $\eps >0$. The conclusion of the estimate in Theorem \ref{thm:main-homog} is then the same as in \cite{bourgade2018extreme}, and so we omit it. \qed

\section{Short range approximations}

\subsection{Finite speed of propogation} \label{sec:fs}

We begin by proving the following estimate.
\bep \label{prop:fs-1} Fix an index $a$ satisfying $a \geq \varphi^3 \ell$.
With overwhelming probability, there is an event so that the following holds. For all times $t$ satisfying
\beq
0 \leq t \leq \varphi^{-10} \frac{ \ell}{ N^{1/3} (  a )^{2/3}} =: T
\eeq
we have that,
\beq
\US (0, t)_{ij} + \US (0, t)_{ji} \leq \e^{ - c \varphi}
\eeq
for any indices $i, j$ satisfying
\beq
i \leq a - \varphi^3 \ell, \qquad a + \varphi^3 \ell \leq j .
\eeq
\eep
\proof We first bound $\US(0, t)_{ji}$, detailing changes to bound $\US (0, t)_{ij}$ at the end of the proof. Let $w_t$ satisfy the equation $\dot{w}_t = \sS w_t$ with initial data
\beq \label{eqn:es-a2}
w_0 (b) = \1_{ \{ b \leq a - \varphi^3 \ell \} }.
\eeq
Let $\psi : \rr \to \rr$ be a function satisfying the following:
\begin{enumerate}[label=(\roman*)]
\item $\psi (x) = 0$ for $x \geq \gamma_{a - \varphi^2 \ell}$
\item $\psi'(x) = -\frac{1}{2}$ for $x \in (\gamma_{a+\varphi^2\ell}, \gamma_{a} )$
\item $\psi'(x) =0$ for $x \leq \gamma_{a + 2 \varphi^2 \ell}$ \label{it:psi-const}
\item $\psi'(x) \leq 0$ for all $ x\in \rr$
\item $| \psi''(x)| \leq C \frac{ N^{2/3} a^{1/3}}{\ell \varphi^2}$ for all $ x \in \rr$ \label{it:psi-second-der}
\item $| \psi'(x) | \leq 1$ for all $x \in \rr$ \label{it:psi-first-der}
\item $\psi(x) \asymp  \gamma_{a } - \gamma_{a + \varphi^2 \ell} \asymp \frac{ \ell \varphi^2}{N^{2/3} a^{1/3}}$ for $x \geq \gamma_{a +2 \varphi^2 \ell}$. \label{it:psi-growth}
\end{enumerate}
Define now
\beq
F(t) := \sum_{k=1}^N w_t (k)^2 \e^{ \nu \psi ( x_k (t) ) }  =: \sum_{k=1}^N w_t (k)^2 \phi_k^2 =: \sum_{k=1}^N v_k^2 .
\eeq
Let $\tau$ be the stopping time,
\begin{align}
\tau &:=  \tau_1 \wedge \tau_2, \notag\\
\tau_1 &:= \inf \left\{ 1> t >0 : \exists i : |x_i(t) - \gamma_i (t) | > \varphi^{1/10} (i^{1/3} N^{2/3} )^{-1} \right\}  , \quad \tau_2 := \inf \{ 1> t >0 : F(t) > 10 a \}.
\end{align}
Note that $\tau_1 = 1$ with overwhelming probability. 
By the Ito lemma we have,
\begin{align}
    \d F &= - \frac{1}{N} \sum_{|i-j| < \ell} \frac{ (v_i -v_j)^2}{ (x_i -x_j)^2} \label{eqn:es-1} \\
    &+ \frac{1}{N} \sum_{|i-j| < \ell} \frac{1}{ (x_i -x_j )^2} v_i v_j \left( \frac{ \phi_i}{\phi_j} + \frac{\phi_j}{\phi_i} -2 \right) \label{eqn:es-2} \\
    &+ \nu \sum_i v_i^2 \psi'(x_i ) \d x_i  \label{eqn:es-3} \\
    &+ \sum_i v_i^2 \left( \frac{\nu^2}{N} (\psi'(x_i  ))^2 + \frac{\nu}{N} \psi''(x_i ) \right) \d t \label{eqn:es-4}
\end{align}
We first consider the term $\eqref{eqn:es-2}$. If a summand in \eqref{eqn:es-2} is non-zero then we must have that either $i$ or $j$ is less than $a + 2 \varphi^2 \ell + \varphi \ell$, as otherwise we have $\phi_i = \phi_j$ (by property \ref{it:psi-const} of $\psi$), as long as $t < \tau_1$.  Then the restriction that $|i-j| \leq \ell$ implies that $i$ and $j$ are both smaller than $a + 3 \varphi^2 \ell$.  By a similar argument we also have both $i$ and $j$ are greater than $a- 2 \varphi^2 \ell$.  For such $i$, $j$ such that $|i-j | \leq \ell$ we have that
\beq
\nu | \psi ( x_i  ) - \psi (x_j ) | \leq \nu |x_i -x_j| \leq \frac{ \nu |i-j| }{ N^{2/3} a^{1/3}} ,
\eeq
for $t < \tau_1$. 
Consequently, as long as 
\beq
\frac{ \nu \ell}{ N^{2/3} a^{1/3} } \leq 1
\eeq
we can Taylor expand the term in the parentheses in \eqref{eqn:es-2} and estimate
\begin{align}
& \frac{1}{N} \sum_{|i-j| < \ell} \frac{1}{ (x_i -x_j )^2} v_i v_j \left( \frac{ \phi_i}{\phi_j} + \frac{\phi_j}{\phi_i} -2 \right) \leq  \frac{C}{N} \sum_{ |i-j| < \ell} \frac{|v_i v_j|}{ (x_i -x_j)^2}  v^2 (x_i -x_j)^2 \notag\\
 \leq & C  \frac{ \nu^2 \ell}{N} F,
\end{align}
with the last line following from Cauchy Schwarz, for $t < \tau_1$. 

The Ito terms in \eqref{eqn:es-4} may be bounded by
\beq
 \sum_i v_i^2 \left( \frac{\nu^2}{N} (\psi'(x_i  ))^2 + \frac{\nu}{N} \psi''(x_i  ) \right)  \leq \frac{ \nu^2}{N} F + \frac{ \nu N^{2/3} a^{1/3}}{ \ell \varphi^2 N} F
\eeq
using properties  \ref{it:psi-second-der} and \ref{it:psi-first-der} of $\psi$. 

We write the term \eqref{eqn:es-3} as 
\begin{align}
\sum_i v_i^2 \psi'(x_i  ) \d x_i  &= \sum_i v_i^2 \psi'(x_i ) \frac{ \d B_i}{\sqrt{N} } \label{eqn:es-3-5}\\
&+ \sum_i v_i^2 \psi'(x_i ) \left( \frac{1}{N} \sum_{ |i - j | < \varphi } \frac{1}{x_i -x_j} \right) \label{eqn:es-3-1} \\
&+\sum_i v_i^2 \psi'(x_i  ) \left( \frac{1}{N} \sum_{ |i - j | > \varphi} \frac{1}{x_i -x_j} - \frac{x_i}{2} \right) \label{eqn:es-3-2}
\end{align}
The quadratic variation of the Martingale \eqref{eqn:es-3-5} is bounded by
\beq
\left \langle \sum_i v_i^2 \psi'(x_i ) \frac{ \d B_i}{\sqrt{N} }  \right \rangle \leq \frac{1}{N} \sum_{i=1}^N v_i^4 \leq \frac{F^2}{N}.
\eeq
Therefore by the BDG inequality,
\beq
\sup_{ 0 < s < \tau \wedge T } \left| \nu \sum_{i=1}^N v_i^2 \psi' (x_i ) \frac{ \d B_i}{ \sqrt{N}} \right| \leq ( \log N)^2 \frac{ \nu a T^{1/2}}{N^{1/2}}
\eeq
with overwhelming probability. 

The term \eqref{eqn:es-3-1} can be written as
\begin{align}
\sum_i v_i^2 \psi'(x_i - \gamma_a ) \left( \frac{1}{N} \sum_{ |i - j | < \varphi } \frac{1}{x_i -x_j} \right) &= \frac{1}{N} \sum_{|i-j| < \varphi } v_i^2 \frac{ \psi' ( x_i ) - \psi' (x_j )}{ x_i -x_j} \label{eqn:es-3-3} \\
&+ \frac{1}{N} \sum_{ |i-j| < \varphi  } \frac{ \psi' (x_i) (v_i^2 - v_j^2)}{ x_i -x_j} \label{eqn:es-3-4}
\end{align}
Using property \ref{it:psi-second-der} of $\psi$ we have for the conribution of \eqref{eqn:es-3-3},
\beq
\left| \nu \frac{1}{N} \sum_{|i-j| < \varphi } v_i^2 \frac{ \psi' ( x_i ) - \psi' (x_j )}{ x_i -x_j} \right| \leq \nu \frac{ F N^{2/3} a^{1/3}}{  N \ell}
\eeq
for $t < \tau_1$. 

 Writing $v_i^2 -v_j^2 = (v_i - v_j)(v_i + v_j)$ and applying Cauchy-Schwarz, the line \eqref{eqn:es-3-4} can be bounded as
\beq
\frac{\nu}{N} \sum_{ |i-j| < \varphi  } \frac{ \psi' (x_i) (v_i^2 - v_j^2}{ x_i -x_j} \leq \frac{1}{ 10 N} \sum_{|i-j| < \varphi} \frac{ (v_i -v_j)^2}{(x_i -x_j)^2} + \frac{ \nu^2}{N} \varphi^2 F
\eeq
the first term being absorbed into the negative term \eqref{eqn:es-1}.

 It remains to estimate the term \eqref{eqn:es-3-2}. First we note that if   the summand is non-zero, then $i \asymp a$. We have that by rigidity and \eqref{eqn:misc-3},
\begin{align}
\frac{1}{N} \sum_{ j :|i-j| > \varphi} \frac{1}{x_i -x_j} - \frac{x_i}{2} &= \Re[ \msc ( \gamma_i ) ] - \frac{\gamma_i}{2} \notag \\
& - \int_{\gamma_{i - \varphi}}^{ \gamma_{i+\varphi}} \frac{ \rhosc (x)}{ x - \gamma_i } + \O \left( \varphi \frac{ a^{1/3}}{N^{1/3}} \right) \label{eqn:es-a1}
\end{align}
for $t < \tau_1$. 
Note that the RHS of the first line in \eqref{eqn:es-a1} is $0$. We turn to estimating the deterministic term on the second line of \eqref{eqn:es-a1}. We have, 
\begin{align}
\int_{\gamma_{i - \varphi}}^{ \gamma_{i+\varphi}} \frac{\rhosc (x)}{ x- \gamma_i } \d x = \int_{\gamma_{i - \varphi}}^{ \gamma_{i+\varphi}} \frac{\rhosc (x)- \rhosc( \gamma_i )}{ x- \gamma_i } + \rhosc(\gamma_i) \int_{\gamma_{i - \varphi}}^{ \gamma_{i+\varphi}} \frac{1}{ x- \gamma_i } 
\end{align}
Using $\rhosc(x) - \rhosc(y) \asymp |x-y| / (x^{1/2} + y^{1/2} )$ the first term is
\beq
\int_{\gamma_{i - \varphi}}^{ \gamma_{i+\varphi}} \frac{\rhosc (x)- \rhosc( \gamma_i )}{ x- \gamma_i }  \leq C \sqrt{ \kappa_i } \leq C \frac{ a^{1/3}}{N^{1/3}}
\eeq
and the second term can be bounded by $\rhosc( \gamma_i ) \log(N) \leq C ( \log N) \frac{a^{1/3}}{N^{1/3}}$.

Collecting all of the above estimates we see that for any $t < \tau \wedge T$ we have,
\beq
F(  t ) \leq a + C a \varphi^2 \left( \frac{ T^{1/2} \nu}{N^{1/2}} + \frac{ \nu^2 \ell T }{N} +\frac{ \nu T N^{2/3} a^{1/3}}{ \ell N} + \frac{ \nu T a^{1/3}}{N^{1/3}} \right)
\eeq
as long as $\nu \ell / ( N^{2/3} a^{1/3} ) \leq 1$. Choose $\nu = N^{2/3} a^{1/3} / ( \ell \varphi)$ and use $T \leq \ell / (N^{1/3} a^{2/3} \varphi^{10} )$. Then the above inequality simplifies to,
\beq
F(  t ) \leq a + C a \varphi^{-1}
\eeq
If $\tau < T$ then $\lim_{ t \to \tau^+} F(t) \geq 10 a$. The above inequality rules this out, and we therefore see that $\tau = T$ with overwhelming probability. On the other hand, for $i \geq a + 3 \varphi^2 \ell$ we see from property \ref{it:psi-growth} that
\beq
\nu \psi (x_i) \geq v \varphi.
\eeq
Therefore, for $t < \tau$ we have for $j > a + 3 \varphi^2 \ell$ and $i < a - 3 \varphi^2 \ell$ that
\beq
(\US (0, t)_{ji} )^2 \e^{ \nu \psi (x_j (t) ) } \leq F(t) \leq 10 N
\eeq
and so $\US_{ji} (0, t) \leq \e^{ - ( \log N)^2}$. This completes the estimate for $\US_{ji}$. The argument for $\US_{ij}$ is exactly the same, except that one takes for initial data in \eqref{eqn:es-a2} a step function of indices $b > a +3 \varphi^2 \ell$,  and one takes a function $\psi$ that is $0$ for $x$ near $\gamma_{a + 3 \varphi^2 \ell}$ and has slope $\psi'(x) = \frac{1}{2}$ for $x \in ( \gamma_{a + 2 \varphi^2}, \gamma_{a + \varphi^2 \ell} )$. \qed

We extend the above estimate. 
\bep \label{prop:fs-box}
Fix an index $b$ satisfying $b \geq \varphi^5 \ell$. Define,
\beq
T := \varphi^{-11} \frac{ \ell}{ N^{1/3} b^{2/3}}.
\eeq
For any times $0 \leq s \leq t \leq T$ we have that
\beq
\US (s, t)_{ij} + \US (s, t)_{ji} \leq \e^{ - c \varphi }
\eeq
for any indices $i, j$ satisfying
\beq
i \leq b - \varphi^4 \ell, \qquad b + \varphi^4 \ell \leq j.
\eeq
\eep
\proof Due to Proposition \ref{prop:fs-1} we may work on the event $\E$ on which the following holds. For all $a$ satisfying
\beq
\frac{b}{2} \leq a \leq \frac{3}{2} b
\eeq
we have that for any $n, m$ satisfying
\beq
n \leq a - \varphi^3 \ell, \qquad a + \varphi^3 \ell \leq m
\eeq
we have
\beq
\US (0, r)_{nm}  + \US (0, r)_{mn} \leq \e^{ - c \varphi} ,
\eeq
for all $0 < r < T$. In the remainder of the proof we will show that the desired estimates hold on the event $\E$.

Let $i$ and $j$, $s$ and $t$ be as in the statement of the Proposition we are proving. As $\US$ is a semigroup and its matrix elements are non-negative we have for any $n$
\begin{align}
\US (0, t)_{jm} &= \sum_{k} \US(s, t)_{jk} \US_{km}(0, s) \notag\\
&\geq \US (s, t)_{ji} \US (0, s)_{im} .
\end{align}
Moreover, $\US(0, r)$ is a doubly stochastic matrix (since the constant vector is an eigenvector with eigenvalue $0$ of $\sS$). Therefore, there exists an $m$ so that
\beq
\US (0, s)_{im} \geq N^{-1}.
\eeq
Consider now the index $a_1 := b - \varphi^{3.5} \ell$. Then,
\beq
i \leq a_1 -  \varphi^3 \ell
\eeq
and so if $m > a_1 + \varphi^3 \ell$ then $\US (0, s)_{im} \leq N^{-2}$  on the event $\E$. Therefore $m \leq a_1 + \varphi^3 \ell$. Consider the index $a_2 := b + \varphi^{3.5} \ell$. Then $m \leq a_2 - \varphi^3 \ell$ and $j \geq a_2 + \varphi^{3.5} \ell$. Therefore, on the event $\E$ we have that $\US(0, t)_{jm} \leq \e^{ -c \varphi}$ and so
\beq
\US (s, t)_{ji} \leq N \US(0, t)_{jm} \leq \e^{ -c/2 \varphi}.
\eeq 
This completes the estimate of $\US (s, t)_{ji}$. The estimate of $\US(s, t)_{ij}$ is similar. We have,
\beq
\US(0, t)_{im} \geq \US(s, t)_{ij} \US(0, s)_{jm}
\eeq
for some $m$ so that $\US (0, s)_{jm} \geq N^{-1}$. As before, we see that $m \geq a_2$ and then that $\US(0, t)_{im} \leq \e^{ - c \varphi}$ as $i \leq a_1 - \varphi^3 \ell$ and $a_1 + \varphi^3 \ell \leq m$. This completes the proof. \qed

\subsubsection{Proof of Theorem \ref{thm:fs}} \label{sec:fs-proof}

This now follows easily from Proposition \ref{prop:fs-box}. Indeed, defining $T_b := \varphi^{-11} \frac{ \ell}{N^{1/3}(\ell+b)^{2/3}}$ we see that for any set of the form $S_u := \{ (s, t) \in [0, 1]^2 : u \leq s \leq t \leq u +T_b\}$ in the square $[0, 1]^2$ we have 
\beq
\US (s, t)_{ib} + \US (s, t)_{bi} \leq \e^{ - \varphi}
\eeq
with probability at least $ 1- N^{-D}$  for all $(s, t) \in S_u$ and $|i-b| \geq \varphi^5 \ell$. But then we can cover the set of all $s \leq t \in [0, 1]^2$ s.t. $|s-t| \leq \varphi^{-1} T_b$ by at most $N^2$ of the sets $S_u$ since $T_b \geq N^{-3/2}$. \qed

\subsection{Short range approximation}

\bel \label{lem: short-range approximation}
Assume that Property $\pP_a$ holds. With overwhelming probability we have for all $u, t, k$ such that $u \in [t/2, t]$  and
\beq \label{eqn:sr-assump}
|u-t| \leq \varphi^{-12} \frac{ \ell}{N^{1/3} ( \ell + k)^{2/3}} 
\eeq
we have for some $C_1 >0$
\begin{align}
& \left| [\UB (u, t) \fu (u)]_k - [\US (u, t) \fu (u) ]_k - ( \bu_k (t) - \bu_k (u) ) \right|  \notag\\
\leq & \varphi^{C_1}  |u-t| \sup_{ j : |k-j| \leq \varphi^{5} \ell} \left( \frac{ N \kappa_{j \vee \ell}}{Q_0 (t, j \wedge k)} + \frac{ N \kappa_{j \vee \ell}}{Q_a (t, j) \ell} \right)
\end{align}
under the assumption that for all $j$ s.t. $|j-k| \leq \varphi^{20} \ell $ we have,
\beq \label{eqn:sra-assump-2}
(N \rho_j)^{-1} \leq t^2 + t \rho_j
\eeq
\eel
\proof Let $\be$ be the constant vector $\be = (1, \dots, 1) \in \rr^N$ and
\beq
v_s := \be \times  \del_s \bu (\gamma_k, s).
\eeq
Recall \eqref{eqn:deltbut}, motivating the following. 
  By the Duhamel formula,
\begin{align}
&(\UB (u,t)\fu (u))_k \notag \\
&=(\US (u,t) \fu (u))_k +  \int_u^t ( \US (s,t)\lL \fu (s))_kds \notag \\
&=(\US (u,t) \fu (u))_k + \int_u^t ( \US (s,t) v_s )_k \d s + \int_u^t ( \US (s,t)(\lL \fu (s)-v_s))_k \d s \notag \\
&=(\US (u,t)\fu (u))_k + \int_u^t v_s \d s + \int_u^t ( \US (s,t)(\lL \fu (s) - v_s ))_k \d s \notag \\
&=(\US (u,t) \fu (u))_k + (\bu_k(t) - \bu_k(u)) +  \int_u^t ( \US (s,t)(\lL \fu (s) - v_s ))_k \d s .
\end{align}
In the third equality we used the fact that $\US \be = \be$. In order to prove the lemma, it suffices to estimate the last integral. By Theorem \ref{thm:fs}, the assumption \eqref{eqn:sr-assump}, and the fact that $\| \US \|_{\ell^\infty \to \ell^\infty} \leq 1$ we have that
\begin{align}
  \left|  \int_u^t ( \US (s,t)(\lL \fu (s) - v_s ))_k \d s \right| 
\leq |u-t| \times \sup_{s \in [u, t]} \sup_{j : |j-k| \leq \varphi^{5} \ell} | (\lL \fu (s) )_j - v_s | +|u-t| \e^{ - c \varphi} .
\end{align}
We estimate the quantity in the supremum on the RHS for some $j : |j-k| \leq \varphi^5 \ell$ and some $s \in [u, t] \subseteq [t/2, t]$. We start with, 
\beq \label{eqn:Lu-dtu-a1}
| (\lL \fu (s) )_j - v_s | \leq | (\lL \fu (s) )_j - \del_s \bu  ( \gamma_j, s) | + | \del_s \bu ( \gamma_j , s ) - \del_s \bu ( \gamma_k, s) | .
\eeq
We start with the first term on the RHS. 
Recalling \eqref{eqn:deltbut} we have, with overwhelming probability,
\begin{align}\label{eqn: |Lu-d_t(u_bar)|}
& |( \lL \fu) _j(s)-\del_s \fu (\gamma_j, s)| \notag \\
&\leq\sum_{i:|i-j|> \ell} \left|  \frac{\fu_i-\fu_j}{N(x_i-x_j)^2}-\int_{\gamma_i}^{\gamma_{i+1}}\frac{\bu (y,s)-\bu (\gamma_j, s)}{(y-\gamma_j)^2}\rhosc(y) \d y \right| + \int_{\gamma_{j+\ell}}^{\gamma_{j-\ell}}\frac{|\bu(y, s)-\bu (\gamma_j, s) |}{(y-\gamma_j)^2}\rhosc (y)\d y\notag \\
&\leq\frac{\varphi}{N}\sum_{i:|i-j|> l}\left( \frac{|\fu_i-\bu_i|+|\fu_j-\bu_j|}{(\gamma_i-\gamma_j)^2}+\left( \frac{1}{N\rho_j}+\frac{1}{N\rho_i} \right) \frac{|\bu_i-\bu_j|}{(\gamma_i-\gamma_j)^3} + \sup_{ y \in [\gamma_i, \gamma_{i+1}] } \frac{ | \bu (y) - \bu_i| }{(\gamma_i - \gamma_j)^2} \right) \notag  \\
&\quad +  \int_{\gamma_{j+l}}^{\gamma_{j-l}}\frac{|\bu (y)-\bu_j|}{(y-\gamma_j)^2}\rhosc (y)\d y
\end{align}
For notational simplicity we suppressed above and in the rest of the proof the time argument $s$ of $\fu, \bu$. 
Note that we used the rigidity estimates \eqref{eqn:rig} in the second estimate. We now estimate each of the  terms on the RHS. Starting with the first term we have, using Property $\pP_a$ and Lemma \ref{lem:uncond-Pa},
\begin{align} \label{est:|Lu-d_t(u_bar)|-a1}
&\frac{1}{N}\sum_{i:|i-j|> \ell}\frac{|\fu_i-\bu_i|+|\fu_j-\bu_j|}{(\gamma_i-\gamma_j)^2} \lesssim \frac{\varphi^{C_1}}{N}\sum_{i:|i-j|> l}\frac{Q_a(s,i\wedge j)^{-1}+Q_0(s,i\wedge j)^{-1}}{(\gamma_i-\gamma_j)^2}\notag \\
&\lesssim\frac{\varphi^{C_1}}{N}\pare{\sum_{i>j+ \ell }\frac{1}{(\gamma_i-\gamma_j)^2Q_a(s,j)}+\mathbbm{1}_{j> \ell}\sum_{i<j-\ell}\frac{Q_a (s, i)^{-1} + Q_0 (s, i)^{-1} }{(\gamma_i-\gamma_j)^2 }    } \notag \\
&\lesssim \varphi^{C_1} \bigg\{ \int_{-2} ^{\gamma_{j+ l}}\frac{\rhosc (y)dy}{(y-\gamma_j)^2Q_a(s,j)} + \mathbbm{1}_{j> l}\int_{\gamma_{j- l}}^2\frac{(\rho(y))^ady}{(y-\gamma_j)^2N^{2-a}(s+\rho(y))^{2-a}s^{1-a}} \notag \\
&+ \mathbbm{1}_{j> l}\int_{\gamma_{j- l}}^2\frac{\d y}{(y-\gamma_j)^2N^{2}(s+\rhosc (y))^{2}s^{1}} \bigg\} \notag\\
& \lesssim \varphi^{C_1} \left( \frac{ \rho_{j+\ell}^2}{\rho_\ell^3 Q_a (s, j)} + \frac{N\kappa_j}{lQ_a(s,j)} \right) 
\lesssim \varphi^{C_1} \frac{N\kappa_{j\vee l}}{\ell Q_a(t,j)}
\end{align}
Above we used \eqref{eqn:misc-4} and \eqref{eqn:misc-5}. For the second and third terms on the RHS of \eqref{eqn: |Lu-d_t(u_bar)|} we have, using \eqref{eqn:bua-dif} and \eqref{eqn:bua-dif-sup} 
\begin{equation}\label{est:|Lu-d_t(u_bar)|,part1}
\begin{aligned}
&\frac{1}{N}\sum_{i:|i-j|> \ell}\left\{ \pare{\frac{1}{N\rho_j}+\frac{1}{N\rho_i}}\frac{|\bu_i-\bu_j|}{(\gamma_i-\gamma_j)^3} + \sup_{ y \in [\gamma_i, \gamma_{i+1}] } \frac{ | \bu (y) - \bu_i| }{(\gamma_i - \gamma_j)^2} \right\} \\
&\lesssim \frac{\varphi}{N}\pare{\sum_{i>j+\ell}\frac{1}{N\rho_j(\gamma_i-\gamma_j)^2Ns(s+\rho_j)^2}+\mathbbm{1}_{j> \ell}\sum_{i<j- \ell}\frac{1}{N\rho_i(\gamma_i-\gamma_j)^2Ns(s+\rho_i)^2}}\\
& \lesssim \frac{\varphi}{Q_0 (t, j)} \int_{-2}^{\gamma_{j+\ell}} \frac{ \rhosc (y)}{ (y- \gamma_j)^2} \d y + \varphi \mathbbm{1}_{j> \ell} \int_{\gamma_{j-\ell}}^2  \frac{1}{(y-\gamma_j)^2 N^2(s+\rhosc(y))^2 s} \d y \\
&\lesssim \varphi \frac{N\kappa_{j\vee \ell}}{\ell Q_0(s,j)}.
\end{aligned}
\end{equation}
again using \eqref{eqn:misc-4} and \eqref{eqn:misc-5}. For the last term, we use regularity of $\bu$ in \eqref{eqn:bu-deriv} (note that our assumption \eqref{eqn:sra-assump-2} implies that this bound is applicable) as well as the fact $\partial_x \bu (x, s), \partial_x^2\bu (x, s)$ are $O(N^{10})$ (by a straightforward estimate) to find 
\begin{equation}\label{est:|Lu-d_t(u_bar)|,part2}
\begin{aligned}
&\varphi^{-1}\int_{\gamma_{j+\ell}}^{\gamma_{j-\ell}}\frac{|\bu (y)-\bu_j|}{(y-\gamma_j)^2}\rhosc(y)\d y\\
&\lesssim \int_{\gamma_{j+l}}^{\gamma_j-N^{-100}}\frac{\rhosc(y)\d y}{Ns(s+\rho_j)^2|y-\gamma_j|} + \int_{\gamma_j+N^{-100}}^{\gamma_{j-l}}\frac{\rhosc(y)\d y}{Ns(s+\rhosc(y))^2|y-\gamma_j|}\\
&\quad + \int_{|y-\gamma_j|<N^{-100}}\left( \frac{( \partial_x \bu ) (\gamma_j )(\rhosc (y)}{y-\gamma_j} + \O ( N^{10} ) \rhosc (y) \right) \d y\\
&\lesssim \varphi \frac{\rho_{j\vee \ell}}{Ns(s+\rho_j)^2}+\O(N^{-50}) \lesssim \varphi \frac{N\kappa_{j\vee\ell}}{Q_0(s,j)} ,
\end{aligned}
\end{equation}
where we used \eqref{eqn:misc-6}. 
In summary, from \eqref{eqn: |Lu-d_t(u_bar)|}, \eqref{est:|Lu-d_t(u_bar)|-a1} \eqref{est:|Lu-d_t(u_bar)|,part1} \eqref{est:|Lu-d_t(u_bar)|,part2} we find 
\beq \label{eqn:sra-est-b4}
|( \lL \fu) _j(s)-\del_s \fu (\gamma_j, s)|  \lesssim \varphi^{C_1} \left( \frac{N\kappa_{j\vee \ell }}{\ell Q_a(s,j)}+\frac{N\kappa_{j\vee\ell}}{Q_0(s,j)} \right). 
\eeq
We now turn to the second term on the RHS of \eqref{eqn:Lu-dtu-a1}. We start with,
\begin{align}
&\abs{\bar{\bB}\bu_j(s)-\bar{\bB}\bu_k(s)} = \abs{\int \pare{\frac{\bu(y)-\bu_j}{(y-\gamma_j)^2}-\frac{\bu(y)-\bu_k}{(y-\gamma_k)^2}}\rhosc(y)\d y} \notag \\
&\leq \abs{\int_{\gamma_{j+\varphi^6 \ell}}^{\gamma_{j-\varphi^6 \ell}}\frac{\bu (y)-\bu_j}{(y-\gamma_j)^2}\rhosc(y)\d y}+\abs{\int_{\gamma_{j+\varphi^6 \ell}}^{\gamma_{j-\varphi^6 \ell}}\frac{\bu (y)-\bu_k}{(y-\gamma_k)^2}\rhosc(y) \d y} \notag \\
&+\int_{[-2,\gamma_{j+\varphi^6 \ell}]\cup [\gamma_{j-\varphi^6 \ell},2]} \left\{ \frac{|\bu_j-\bu_k|}{(y-\gamma_j)^2}+\pare{\frac{1}{(y-\gamma_j)^2}-\frac{1}{(y-\gamma_k)^2}}|\bu (y)-\bu_k| \right\} \rhosc(y) \d y \label{eqn:sra-est-b1}
\end{align}
The two terms in the second to last line is similar to the last term in \eqref{eqn: |Lu-d_t(u_bar)|}. Repeating the argument in \eqref{est:|Lu-d_t(u_bar)|,part2} we obtain (recall also $|k-j| \leq \varphi^5 \ell$), 
\beq \label{eqn:sra-est-b2}
\abs{\int_{\gamma_{j+\varphi^6 \ell}}^{\gamma_{j-\varphi^6 \ell}}\frac{\bu (y)-\bu _j}{(y-\gamma_j)^2}\rho(y)dy} + \abs{\int_{\gamma_{j+\varphi^6 \ell}}^{\gamma_{j-\varphi^6 \ell}}\frac{\bu (y)-\bu_k}{(y-\gamma_k)^2}\rhosc(y)dy} \lesssim \varphi^{10} \left(  \frac{N\kappa_{j \vee\ell}}{Q_0(s,j)} + \frac{N\kappa_{k \vee\ell}}{Q_0(s,k)} \right)
\eeq
By regularity of $\bu $ in \eqref{eqn:bua-dif}
\begin{equation}\label{est:|B_bar u_bar_j - B_bar u_bar_j|}
\begin{aligned}
&\varphi^{-1} \int_{[0,\gamma_{j+\varphi^6 \ell}]\cup [\gamma_{j-\varphi^6 \ell},2]}\frac{|\bu_j-\bu_k|}{(y-\gamma_j)^2}\rhosc (y)dy \lesssim \int_{[0,\gamma_{j+\varphi^6 \ell}]\cup [\gamma_{j-\varphi^6 \ell},2]}\frac{|\gamma_k-\gamma_j|\rhosc (y)}{Ns(s+\rho_{j\wedge k})^2(y-\gamma_j)^2} \d y\\
&\lesssim \frac{|\gamma_j-\gamma_k|}{Ns(s+\rho_{j\wedge k})^2} \sum_{i:|i-j|>\varphi^6 \ell}\frac{1}{N(\gamma_i-\gamma_j)^2} \lesssim\frac{|\gamma_j-\gamma_k|}{Ns(s+\rho_{j\wedge k})^2}\sum_{i:|i-j|>\varphi^6 \ell}\frac{N(\kappa_i+\kappa_j)}{(i-j)^2}\\
&\lesssim \varphi^{10} \frac{\kappa_j+\kappa_l}{Ns(s+\rho_{j\wedge k})^2\rho_{j\wedge k}} \lesssim \varphi^{10} \frac{N\kappa_{j\vee \ell}}{Q_0(s,j\wedge k)}.
\end{aligned}
\end{equation}

For $y > \gamma_{j-\varphi^6 \ell}$ or $y<\gamma_{j+\varphi^6 \ell}$, we have $|y-\gamma_j|\asymp|y-\gamma_k|$. Then, a similar computation to \eqref{est:|B_bar u_bar_j - B_bar u_bar_j|} gives,
\begin{align}
&\int_{[0,\gamma_{j+\varphi^6 \ell}]\cup [\gamma_{j-\varphi^6 \ell},2]}\pare{\frac{1}{(y-\gamma_j)^2}-\frac{1}{(y-\gamma_k)^2}}|\bu (y)-\bu_k|\rhosc (y) \d y \notag \\
&\lesssim \int_{[0,\gamma_{j+\varphi^6 \ell}]\cup [\gamma_{j-\varphi^6 \ell},2]} \frac{|\gamma_k-\gamma_j||\bu (y)-\bu_k|}{|y-\gamma_k|^3}\rhosc (y) \d y \notag \\
&\lesssim \varphi \int_{[0,\gamma_{j+\varphi^6 \ell}]\cup [\gamma_{j-\varphi^6 \ell},2]} \frac{|\gamma_k-\gamma_j|}{|y-\gamma_k|^2Ns(s+\rho_k\wedge \rho(y))^2}\rhosc (y) \d y \notag \\
&\lesssim \varphi \frac{\kappa_j+\kappa_l}{Ns(s+\rho_{j\wedge k})^2\rho_{j\wedge k}} \lesssim \varphi^{C_1} \frac{N\kappa_{j\vee \ell}}{Q_0(s,j\wedge k)}, \label{eqn:sra-est-b3}
\end{align} 
where we used \eqref{eqn:misc-7}. 

In summary, from \eqref{eqn:sra-est-b1}, \eqref{eqn:sra-est-b2}, \eqref{est:|B_bar u_bar_j - B_bar u_bar_j|} and \eqref{eqn:sra-est-b3} we have, 
\beq
\abs{\bar{\bB}\bu_j(s)-\bar{\bB}\bu_k(s)} \lesssim \varphi^{C_1} \frac{N\kappa_{j\vee l}}{Q_0(s,j\wedge k)} .
\eeq
The above estimate together with \eqref{eqn:sra-est-b4} complete the proof. 
\qed

\subsubsection{Proof of Lemma \ref{lem:|u-m|}} \label{sec:|u-m|-proof}

Duhamel formula implies,
\beq
m_k(t)  = \pare{\US (u,t)\Av(u) \fu (u)}_k + \int_u^t \left( \US (s,t) \del_s \bu_k (s) \right)_k \d s.
\eeq
 Note that the integrand in the second term is the $k$th entry of the vector obtained from the operator $\US (s, t)$ acting on the constant vector that has value $\del_s \bu_k (s)$. Therefore
\beq
 \int_u^t \left( \US (s,t) \del_s \bu_k (s) \right)_k \d s= \bu_k (t) -\bu_k (u) .
\eeq
Recalling that that $\fu_k(t)  = \pare{\UB (u,t) \fu (u)}_k$, we get
\begin{align}
\fu_k(t) - m_k(t)  &= ((\UB (u,t) - \US (u,t)) \fu(u))_k - (\bu_k(t)-\bu_k(u)) \notag \\
&\quad +\US (u,t)( \fu (u) - \Av(u) \fu (u))_k\label{eqn:u-m-est-1}
\end{align}
By finite speed of propagation in Theorem  \ref{thm:fs}, and the fact that $[ \Av u]_j = u_j$ for $|j-k| \leq \varphi^{100} \ell$ (due to the assumption that $R > \varphi^{100} \ell$) we have, 
\beq
\US(u,t)(\fu(u) - \Av(u) \fu (u))_k =\O ( \e^{ -c \varphi} )
\eeq
with overwhelming probability. 
The claim then follows from using short-range approximation in Lemma \ref{lem: short-range approximation} to estimate the first line of \eqref{eqn:u-m-est-1}. \qed 

\subsection{Estimate for averages}

\bel \label{lem: ave mod dyn, part 1}
Let $|j - k | \leq 3R$ and let $ z= \gamma_j + \i \eta$. Assume that
\beq
\eta \geq \varphi \frac{1}{N \rho_j}, \qquad |u-t| \leq \varphi^{-100} \frac{ \ell}{N^{1/3} ( k + R)^{2/3} }
\eeq
and that $k \in \Da$ and that Property $\pP_a$ holds. Consider,
\beq \label{eqn:ave-part-1}
\frac{1}{ N \sqrt{\kappa_j+\eta}}\Im \sum_{i : |j-i| < \ell } \frac{\pare{\US (u,s)\Av(u) \fu (u)}_i-\pare{\Av(s)\UB (u,s) \fu (u)}_i + (\bu_k(s)-\bu_k(u))}{x_i-z}.
\eeq
For $k \in \Db$ the term in \eqref{eqn:ave-part-1} is bounded by (up to factors of $\varphi$)
\beq
\left( \frac{\ell}{R}+ \frac{ |u-t| N \kappa_k}{\ell} \right) \frac{1}{Q_a (t, k)}+ \frac{R + |u-t| N \kappa_k}{Q_0 (t, k) }
\eeq
and for $k \in \De$ the term in \eqref{eqn:ave-part-1} is bounded by  (up to factors of $\varphi$)
\beq
\left( \1_{ j \leq \varphi^{12} \ell} \frac{ |u-t| N \kappa_\ell}{Q_a(t, 1) \ell} +\1_{ j > \varphi^{12} \ell} \frac{ |u-t| N \kappa_j}{Q_a(t, j) \ell}  + \frac{R+ |u-t| N \kappa_R N^{1/3}\rho_k}{Q_0 (t, k)}\right)+ \frac{\ell}{R} \frac{1}{Q_a (t, R)}
\eeq
\eel
\proof 
Let us focus, for fixed $i$, on  the numerator in \eqref{eqn:ave-part-1}. Note that because $|i - j | < \ell$ and our assumptions on $|u-t$ we may use Theorem \ref{thm:fs} for $\US$ with one index being $i$. 
From definition \eqref{def: average operator}, we have that 
\begin{align}\label{quant: ave mod dyn, part 1, numerator}
&\pare{\US (u,s)\Av(u) \fu (u)}_i-\pare{\Av(s)\UB (u,s) \fu (u)}_i + (\bu_k(s)-\bu_k(u)) \notag\\
= &\frac{1}{R} \sum_{ R \leq b \leq 2R }\left\{  \pare{\US (u,s)\Flat_b(u) \fu (u)}_i-\pare{\Flat_b(s)\UB(u,s) \fu (u)}_i+(\bu _k(s)-\bu_k(u)) \right\}
\end{align}
Let us now consider the summand on the RHS for fixed $b$.
If $|i-k|>b+\varphi^{10} \ell$, then
\beq
\pare{\Flat_b(s)\UB (u,s) \fu (u)}_i = \bu_k(s),
\eeq
and by finite speed of propagation Theorem \ref{thm:fs}
\beq
\pare{\US (u, s) \Flat_b(u) \fu (u)}_i = \pare{\US(u,s) \bu_k(u)}_i  + O(N^{-D}) = \bu_k(u)+\O ( \e^{ - c \varphi} )
\eeq
So the summand in \eqref{quant: ave mod dyn, part 1, numerator} is $\O ( \e^{ -c \varphi} )$ if $|i-k| > b + \varphi^{10} \ell$.

If $b<|i-k|<b+\varphi^{10} \ell$, the summand in \eqref{quant: ave mod dyn, part 1, numerator} becomes (note that $(\Flat_b \UB(u, s) \fu (u) )_i = \bu_k (s)$ for such $i$)
\begin{align}
(\US (u,s)\Flat_b(u) \fu (u))_i - \bu_k(u) &= \pare{\US (u,s)\pare{\Flat_b(u)\fu (u) - \bu_k(u)}}_i \notag \\
&=\pare{\US(u,s)\mathbbm{1}_{\brac{|\cdot-k|<b}}\pare{\fu (u)-\bu_k(u)}}_i \notag \\
&\lesssim \sup_{m:b-\varphi^{10} \ell<|m-k|<b}|\fu_m(u)-\bu_k(u)| + \O( \e^{-c \varphi} ) .
\end{align}
In the first line we used that $\US$ is an identity on constant function. In the second line, the function $\mathbbm{1}$ is a cutoff on the spatial index. In the third line, we used that $\US$ is a contraction on $\ell^\infty$ as well as Theorem \ref{thm:fs}. 

 If $b-\varphi^{10}\ell<|i-k|\leq b$, using a similar argument, \eqref{quant: ave mod dyn, part 1, numerator} becomes
\begin{align*}
&(\US \Flat_b(u) \fu (u))_i - (\UB \fu (u))_i + (\bu_k(s)-\bu_k(u))\\
& = \pare{\US \pare{\Flat_b(u)\fu (u)-\fu (u)}}_i + \sqbrac{(\US \fu (u))_i - (\UB \fu (u))_i + (\bu_k(s)-\bu_k(u))}\\
& =  \sup_{m:b<|m-k|<b+\varphi l}|\fu_m(u)-\bu_k(u)|+ \sqbrac{(\US \fu (u))_i - (\UB \fu (u))_i + (\bu_k(s)-\bu_k(u))}.
\end{align*}
If $|i-k|<b-\varphi^{10} \ell $, \eqref{quant: ave mod dyn, part 1, numerator} is (again using Theorem \ref{thm:fs}) 
\beq
(\US \fu (u))_i - (\UB \fu (u))_i + (\bu_k(s)-\bu_k(u)) + \O( \e^{ - c \varphi} ).
\eeq

Therefore, we have for the second line of \eqref{quant: ave mod dyn, part 1, numerator}, 
\begin{align}
 & \left| \frac{1}{R} \sum_{ R \leq b \leq 2R }\left\{  \pare{\US (u,s)\Flat_b(u) \fu (u)}_i-\pare{\Flat_b(s)\UB(u,s) \fu (u)}_i+(\bu _k(s)-\bu_k(u)) \right\} \right| \notag\\
\leq &  
\frac{1}{R}\sum_{b:| |i-k|-b|  \leq \varphi^{10} \ell}\sup_{m:| |m-k|-b| \leq \varphi^{10} \ell}|\fu_m(u)-\bu_k(u)| \notag  \\
+ & \label{eqn: ave mod dyn, part 1-2}
\frac{1}{R}\sum_{b:|i-k|<b} \left| (\US \fu (u))_i - (\UB \fu (u))_i + (\bu_k(s)-\bu_k(u)) \right| .
\end{align}
We have,
\beq
\label{eqn: ave mod dyn, part 1-1}
\frac{1}{R}\sum_{b:| |i-k|-b|  \leq \varphi^{10} \ell}\sup_{m:| |m-k|-b| \leq \varphi^{10} \ell}|\fu_m(u)-\bu_k(u)|  \leq \varphi^{11} \mathbbm{1}_{\brac{|j-k|\asymp R}}\frac{\ell}{R}\sup_{m:|m-k|\asymp R}|\fu_m(u)-\bu_k(u)| .
\eeq
Here we used the fact that in the LHS of \eqref{eqn: ave mod dyn, part 1-1}, since $|i-j|<\varphi^{10} \ell$ the sum over $b$ is empty unless $|j-k|\asymp R$. By decomposing
\beq
|\fu_m(u)-\bu_k(u)|\leq |\fu_m(u)-\bu_m(u)| + |\bu_m(u)-\bu_k(u)|,
\eeq
and using Property $(\pP_a)$ (recall Lemma \ref{lem:index-properties}) and the spatial regularity of $\bu$ in \eqref{eqn:bua-dif}, we derive that \eqref{eqn: ave mod dyn, part 1-1} is bounded by, for $(t,k) \in \Db$, 
\beq
\frac{\ell}{R}\sup_{m:|m-k|\asymp R}|\fu_m(u)-\bu_k(u)| \leq \varphi^{C_1} \left( \frac{\ell}{R}\frac{1}{Q_a(t,k)} + \frac{\ell}{N\rho_k}\frac{1}{Nt(t+\rho_k)^2}\right).
\eeq
On the other hand for $k \in \De$ we have that since $k \leq R^{1/2}$ and $|m-k| \asymp R$ that $m \gtrsim R$. Additionally, $t \gtrsim \rho_k$ and so, 
\beq
\mathbbm{1}_{\brac{|j-k|\asymp R}}\frac{\ell}{R}\sup_{m:|m-k|\asymp R}|\fu_m(u)-\bu_k(u)| \leq \varphi^{C_1} \left( \frac{\ell}{R} \frac{1}{Q_a (t, R)} + \frac{\ell}{N \rho_k} \frac{1}{N t (t+\rho_k)^2} \right).
\eeq
This completes estimation of the second line of \eqref{eqn: ave mod dyn, part 1-2}. 

Using short-range approximation (Lemma \ref{lem: short-range approximation}), as well as \eqref{eqn:bua-dif} for $|\bu_i-\bu_k|$, the third line of \eqref{eqn: ave mod dyn, part 1-2} is bounded by, for $k \in \Db$,
\begin{align}
& \frac{1}{R}\sum_{b:|i-k|<b} \left| (\US \fu (u))_i - (\UB \fu (u))_i + (\bu_k(s)-\bu_k(u)) \right| \notag\\
\leq & \varphi^{C_1} \left(  \frac{|u-t|N\kappa_k}{Q_a(t,k)\ell}+\frac{|u-t|N\kappa_k+ R}{Q_0(t,k)} \right) ,
\end{align}
and for $k \in \De$ by,  (recall $|i-j| \leq \ell$; also $|i-k| \lesssim R$ so our assumption on $|u-t|$ implies that the $|u-t|$ assumption of Lemma \ref{lem: short-range approximation} is satisfied) 
\begin{align}
& \frac{1}{R}\sum_{b:|i-k|<b} \left| (\US \fu (u))_i - (\UB \fu (u))_i + (\bu_k(s)-\bu_k(u)) \right| \notag\\ 
\leq & \varphi^{C_1} \sup_{i:|i-j|\leq \varphi^{11} \ell} \left(  \frac{|u-t|N\kappa_{i\vee l}}{Q_a(t,i)\ell }+\frac{|u-t|N\kappa_{i\vee \ell}}{Q_0(t,1 \vee (i-\ell))} + \frac{|\gamma_i-\gamma_k|}{Nt(t+\rho_{i\wedge k})^2} \right) \notag\\
\leq & \varphi^{C_1} \left( \1_{ j \leq \varphi^{12} \ell} \frac{ |u-t| N \kappa_\ell}{Q_a(t, 1) \ell} +\1_{ j > \varphi^{12} \ell} \frac{ |u-t| N \kappa_j}{Q_a(t, j) \ell}  + \frac{R+ |u-t| N \kappa_R N^{1/3}\rho_k}{Q_0 (t, k)}\right)
\end{align}
The claim now follows from summing over $i$ and using Lemma \ref{lem:im-m-ell}. \qed

\section{Flow lemmas}
\label{sec:flow} 

In this section important roles will be played by
\begin{align}\label{def: f_u and s_u}
f_u(z) = \e^{-u/2}\sum_{i=1}^N \frac{\fu_i(u)}{x_i(u)-z}, \quad s_u(z) = \frac{1}{N}\sum_{i=1}^N \frac{1}{x_i(u)-z}.
\end{align}

\subsection{Estimate for $f_0$}

%
%
%
%
 
\bel
Let $z_t$ be a characteristic started from $z_0 = z = \gamma_j + \i \eta$. Assume that
\beq
\frac{\varphi}{N \rho_j} \leq \eta \leq \varphi^{-1} ( t^2 + t \rho_j ),
\eeq
and let $t/2 < s < t$. Then,
\beq \label{eqn:f0-est}
\left| \Im\left[ \frac{ \e^{s/2}}{N} f_0 (z_s) - \msc (z) \bu_j (s) \right] - \eta \del_s \bu_j (s) \right| \leq \varphi^3 \frac{ N \eta \rho_j \sqrt{ \kappa_j + \eta}}{Q_0 (t, j) }
\eeq
for $s \in (t/2, t)$ 
\eel
\proof For notational simplicity let us write $E = \gamma_j$ and $E_s = ( \gamma_j)^s$. Note that $z_s$ and $E_s$ are characteristics started from $E+ \i \eta$ and $E$, respectively. We begin by computing the time derivative  $ \del_s \bu_j (s)$ to leading order. We first write it as (using \eqref{eqn:buk-alt-formula})
\begin{align}
\del_s \bu_j (s) &= \frac{1}{N} \sum_{i=1}^N \fu_i (0)  \del_s \left( \frac{ 2 \e^{s/2} \sinh (s/2)}{| \gamma_i - E_s |^2}  \right) \notag\\
&= \frac{1}{N} \sum_i \fu_i (0) \left( \frac{\e^s}{| \gamma_i - E_s|^2} - 2 \e^{s/2} \sinh(s/2) \frac{ \del_s | \gamma_i - E_s|^2}{ | \gamma_i - E_s|^4} \right) .  \label{eqn:f0-est-a1}
\end{align}
For the second term on the last line we have,
\begin{align}
\del_s | \gamma_i - E_s|^2 &= \del_s \left( ( \gamma_i - \cosh(s/2) E)^2 + ( \sinh (s/2) \sqrt{ 4 - E^2} )^2 \right) \notag\\
&= - ( \gamma_i - \cosh(s/2) E) \sinh(s/2) E  \sinh(s/2) \cosh(s/2) (4 - E^2) \notag\\
&= - ( \gamma_i - \cosh(s/2) E) \sinh(s/2) E + \O ( t \kappa_j)
\end{align}
for $s \leq t$. The contribution of the $\O ( t \kappa_j)$ term to the sum on the last line of \eqref{eqn:f0-est-a1} is bounded by, 
\begin{align}
 & \left| \frac{1}{N} \sum_{i=1}^N |\fu_i (0) | 2 \e^{s/2} \sinh(s/2) \frac{ t \kappa_j}{ | \gamma_i - E_s|^4} \right| \leq \varphi \frac{1}{N} \sum_{i=1}^N \frac{1}{N \rho_i} \frac{ t^2 \kappa_j}{ | \gamma_i - E|^4 + (t^2 + t \rho_j )^4 }  \notag\\
\leq & C \varphi \frac{ t^2 \kappa_j}{N ( t^2 + t \rho_j )^3} \leq C \varphi \frac{ N \kappa_j }{ Q_0 (t, j)}.
\end{align}
In the first inequality we used the rigidity estimates and Lemma \ref{lem:buk-denom}. The second inequality uses \eqref{eqn:misc-1}. Combining the above estimate with the last line of \eqref{eqn:f0-est-a1} we have
\beq \label{eqn:f0-est-a3}
\eta \del_s \bu_j (s) = \frac{1}{N} \sum_{i=1}^N \fu_i (0) \left( \frac{ \eta \e^s}{| \gamma_i - E_s|^2} + \frac{ 2 \e^{s/2} \sinh^2(s/2) (\gamma_i - \cosh(s/2) E)E \eta }{ | \gamma_i - E_s|^4} \right) + \O \left( \varphi \frac{ N \eta \kappa_j}{ Q_0 (t, j) } \right).
\eeq
We will now compute the other terms appearing in \eqref{eqn:f0-est} to leading order, showing that they cancel with the terms in the above estimate. By definition,
\beq \label{eqn:f0-est-a2}
\Im \left[ \frac{ \e^{s/2}}{N} f_0 (z_s) - \msc (z) \bu_j (s) \right] = \frac{ \e^{s/2}}{N} \sum_i \fu_i (0) \Im \left[  \frac{1}{ x_i - z_s} - \frac{1}{ \gamma_i - E_s} \frac{ \Im  \msc (z)  }{ \Im \msc (E)  } \right].
\eeq
We first replace $x_i$ in the denominator of the first term in \eqref{eqn:f0-est-a2} by $\gamma_i$. We have,
\begin{align}
 & \left| \frac{1}{N} \sum_{i=1}^N \fu_i (0) \Im\left[ \frac{1}{ x_i - z_s} - \frac{1}{\gamma_i - z_s } \right]  \right| \leq \frac{\varphi}{N} \sum_{i=1}^N \frac{\Im[z_s]}{ N \rho_i } \frac{ \left|  |\gamma_i - z_s |^2 - | x_i - z_s |^2 \right|}{ | \gamma_i - z_s|^4} \notag\\
 \leq & \frac{  \varphi^2}{N} \sum_{i=1}^N \frac{t \sqrt{ \kappa_j + \eta}}{N^2 \rho_i^2} \frac{ 1}{ | \gamma_i - z_s|^3} \leq \frac{ t \sqrt{ \kappa_j + \eta}}{N^2} \frac{1}{ \rho_j (t^2 + t \sqrt{ \kappa_j + \eta} + \eta)^2} \leq \frac{ \sqrt{ \kappa_j  + \eta}}{Q_0 (t, j) } \label{eqn:f0-est-a7}
\end{align}
In the first inequality we used the fact that $| \gamma_i -z_s | \asymp | x_i - z_s|$ since $\eta \gg (N \rho_j)^{-1}$ by assumption. In the second inequality we used $\Im[z_s] \leq C( \eta + t \sqrt{ \kappa_j + \eta} ) \leq C ( t \sqrt{ \kappa_j + \eta} )$ by the assumption $\eta \ll t^2 + t \rho_j$, and that
\beq
\left| | \gamma_i -z_s| - | x_i - z_s| \right| \leq C |x_i - \gamma_i| |\gamma_i -z_s|
\eeq
by difference of square factoring.  The third inequality uses \eqref{eqn:misc-2}. 

We now turn to the second term in \eqref{eqn:f0-est-a2}. We first note the identity,
\beq
\frac{ \Im z_s}{ \Im \msc (z) } = 2 \frac{ \cosh(s/2) \eta + \sinh(s/2) \Im[ \sqrt{ z^2-4} ] }{ - \eta + \Im[ \sqrt{z^2-4} ] } = 2 \sinh(s/2) + \frac{ \e^{s/2} \eta}{ \Im \msc (z) }.
\eeq
Using the above identity in the second equality below we have,
\begin{align}
 & \frac{ \e^{s/2}}{N} \sum_{i=1}^N \fu_i (0) \Im \left[ \frac{1}{ \gamma_i -z_s} - \frac{1}{ \gamma_i - E_s} \frac{ \Im \msc(z) }{ \Im \msc (E)} \right] \notag\\
= & \frac{ \e^{s/2} \Im \msc (z) }{ N} \sum_{i=1}^N \fu_i (0) \left( \frac{ \Im z_s}{ | \gamma_i -z_s|^2 \Im \msc (z) } - \frac{ 2 \sinh (s/2)}{| \gamma_i - E_s|^2} \right) \notag\\
= & \frac{ \e^{s/2} \Im \msc (z) }{ N} \sum_{i=1}^N \fu_i (0) \left[ \frac{ \e^{s/2} \eta}{ \Im \msc (z) } \frac{1}{ | \gamma_i -z_s|^2} + 2 \sinh (s/2) \left( \frac{1}{ | \gamma_i -z_s|^2} - \frac{1}{ | \gamma_i - E_s|^2} \right) \right] \notag\\
= & \frac{1}{N} \sum_{i=1}^N \fu_i (0) \left[ \frac{ \e^s \eta}{ | \gamma_i - z_s|^2} + \left( \e^s \eta + 2 \e^{s/2} \sinh(s/2) \Im \msc (z) \right) \left( \frac{1}{ | \gamma_i - z_s|^2} - \frac{1}{ | \gamma_i - E_s|^2} \right) \right] \label{eqn:f0-est-a4}
\end{align}
Note that the first term on the last line cancels with the first term on the RHS of \eqref{eqn:f0-est-a3}. We will further compute the second term on the last line of \eqref{eqn:f0-est-a4}. In order to do so we start with
\begin{align}
& | \gamma_i - E_s|^2 - | \gamma_i - z_s|^2 = 2 \Re [ \gamma_i - E_s] \Re[ z_s - E_s] + ( \Re [ z_s - E_s] )^2 + ( \Im[ E_s] )^2 - ( \Im[z_s] )^2 \notag\\
= & 2 ( \gamma_i - \cosh (s/2) E) \sinh (s/2) \Re[ \sqrt{z^2-4} ] \notag\\
+ & \O \left( \left( \frac{ s \eta}{ \sqrt{ \kappa_j + \eta} } \right)^2 + ( \eta + s \sqrt{ \eta + \kappa_j} )\left( \eta + s \frac{ \eta}{ \sqrt{\kappa_j + \eta}} \right) \right) \notag\\
= &2 ( \gamma_i - \cosh (s/2) E) \sinh (s/2) \Re[ \sqrt{z^2-4} ] + \O \left( t \eta \sqrt{ \kappa_j + \eta} + t^2 \eta \right) \label{eqn:f0-est-a5}
\end{align}
In the second line we used the estimates,
\begin{align}
\Re[z_s -E_s] &= \O \left( s \frac{ \eta}{ \sqrt{ \kappa_j + \eta}} \right), \qquad \Im [z_s] - \Im[ E_s] = \O \left( \eta + s  \frac{ \eta}{ \sqrt{ \kappa_j + \eta}} \right) \notag\\
\Im[z_s]+ \Im[E_s] &= \O \left( \eta + t \sqrt{ \kappa_j + \eta} \right)
\end{align}
all established by direct computation.  Note that in the last line of \eqref{eqn:f0-est-a5} we used the assumption that $\eta \leq t^2 + t \sqrt{ \kappa_j}$.  Using the fact that $\Im[ \sqrt{z^2-4} ] \Re[ \sqrt{z^2-4} ] = E \eta$ we see that
\beq \label{eqn:f0-est-a6}
2 \Im[ \msc (z) ] \Re[ \sqrt{z^2-4} ] = E \eta + \O\left( \frac{ \eta^2}{ \sqrt{ \kappa_j+\eta} } \right).
\eeq
Using \eqref{eqn:f0-est-a5} and \eqref{eqn:f0-est-a6} we see for the $\sinh(s/2)$ term in on the last line of \eqref{eqn:f0-est-a4} that we have,
\begin{align}
 & \frac{1}{N} \sum_{i=1}^N \fu_i (0) 2 \e^{s/2} \sinh(s/2) \Im \msc (z) \left( \frac{1}{ | \gamma_i -z_s|^2} - \frac{1}{ | \gamma_i - E_s|^2} \right) \notag\\
= & \frac{1}{N} \sum_{i=1}^N \fu_i (0) \frac{ 2 \e^{s/2} \sinh^2 (s/2) ( \gamma_i - \cosh(s/2) E) E \eta + \O \left( ( | \gamma_i - E| + t^2) t^2 \frac{ \eta^2}{ \sqrt{ \kappa_j + \eta} } + t^2 \eta \sqrt{ \kappa_j + \eta} ( t + \sqrt{ \kappa_j + \eta} ) \right) }{| \gamma_i -z_s|^2 | \gamma_i - E_s|^2}  \notag\\
= &\frac{1}{N} \sum_{i=1}^N \fu_i (0) \frac{ 2 \e^{s/2} \sinh^2 (s/2) ( \gamma_i - \cosh(s/2) E) E \eta }{| \gamma_i -z_s|^2 | \gamma_i - E_s|^2} + \O \left( \frac{ N \eta \rho_j \sqrt{ \kappa_j + \eta}}{Q_0 (t, j) } \right) \notag \\
= & \frac{1}{N} \sum_{i=1}^N \fu_i (0) \frac{ 2 \e^{s/2} \sinh^2 (s/2) ( \gamma_i - \cosh(s/2) E) E \eta }{ | \gamma_i - E_s|^4} + \O \left( \frac{ N \eta \rho_j \sqrt{ \kappa_j + \eta}}{Q_0 (t, j) } \right)  \label{eqn:f0-est-a8}
\end{align}
The second estimate uses \eqref{eqn:misc-1}. In the last estimate we used that \eqref{eqn:f0-est-a5} again, and \eqref{eqn:misc-1}. 

We now examine the $\e^s \eta$ term on the last line of \eqref{eqn:f0-est-a4}. We have, using \eqref{eqn:f0-est-a5},
\begin{align}
 & \left| \frac{1}{N} \sum_{i=1}^N \fu_i (0) \e^s \eta \left( \frac{1}{ | \gamma_i -z_s|^2} - \frac{1}{ | \gamma_i -E_s|^2} \right) \right| \notag\\
& \leq \varphi \sum_{i=1}^N \frac{ \eta}{N \rho_i} \frac{ ( | \gamma_i - E | + s^2 ) \frac{ s \eta}{ \sqrt{ \kappa_j + \eta}} +  (t \sqrt{ \kappa_j + \eta } )^2}{ | \gamma_i - z_s|^4} \notag\\
&\leq  \varphi^2 \frac{ N \eta \rho_j \sqrt{ \kappa_j + \eta}}{Q_0 (t, j) } \label{eqn:f0-est-a9}
\end{align}
again using \eqref{eqn:misc-1}. Using \eqref{eqn:f0-est-a7}, \eqref{eqn:f0-est-a4}, \eqref{eqn:f0-est-a8} and \eqref{eqn:f0-est-a9} in \eqref{eqn:f0-est-a2}, we find,
\begin{align}
& \Im \left[ \frac{ \e^{s/2}}{N} f_0 (z_s) - \msc (z) \bu_j (s)\right] \notag\\ 
= & \frac{1}{N} \sum_{i=1}^N \fu_i (0) \left( \frac{ \e^s \eta}{ | \gamma_i - E_s|^2} + \frac{ 2 \e^{s/2} \sinh(s/2)^2 ( \gamma_i - \cosh (s/2) E ) E \eta}{| \gamma_i - E_s|^4} \right)  + \O \left( \frac{ N \eta \rho_j \sqrt{ \kappa_j + \eta}}{Q_0 (t, j) } \right).
\end{align}
We conclude the proof by comparing the above with \eqref{eqn:f0-est-a3}. \qed

\subsection{Main flow lemma}

In this section it will be helpful to denote
\begin{align}\label{def: eta_u}
\eta_u:=u^2+ u\sqrt{\kappa_j+\eta}+\eta.
\end{align}
We recall also the notation 
\begin{align}\label{def:z_t, msc}
\gamma_k^t = \lim_{\eta \to 0+}z_t, \quad \msc(z) = \int \frac{\rho(x)dx}{x-z} = \frac{-z+\sqrt{z^2-4}}{2}.
\end{align}

\bet[Characteristic flow]\label{thm: character flow}Assume Property $(\pP_a)$ holds. The following event holds with overwhelming probability. Let $t \in [0,1]$. Let $ z= \gamma_j + \i \eta$ and $z_t$ the characteristic started at $z$. Assume that,
\beq
\frac{\varphi}{N \rho_j} \leq \eta \leq \varphi^{-1} (t^2 + t \rho_j )
\eeq
Then we have,
\begin{align}
&\abs{f_t(z)-f_0(z) - \e^{-\frac{t}{2}}N\pare{s_t(z)-\msc(z)}\bu_j(t)} \lesssim \varphi^{C_1} \left(  \frac{1}{\eta Q_a (t, j)}+N\rho_j\frac{1}{Q_0(t,j)} \right)
\end{align}
\eet
\begin{remark}
The assumption $\eta \leq \varphi^{-1} ( t^2+t\rho_j)$ is essentially equivalent to $t\sqrt{\eta+\kappa_j} \geq \varphi \eta$. When $t<\rho_j$, this assumption will force $\eta < \kappa_j$.
\end{remark}

We first estimate the short-time increment $f_{\frac{t}{2}}(z_{\frac{t}{2}}) - f_0(z_t)$.
\bel[Short-time increment]\label{lem: short-time increment}
The following event holds with overwhelming probability. Let $t \in [\frac{1}{N},1]$. Let $z = E + \i \eta$ and $z_s$ be the characteristic started from $z$. 
Suppose,
\beq
\frac{\varphi}{N \rho_j} \leq \eta \leq \varphi^{-1} (t^2 + t \rho_j ) .
\eeq
Then, 
\begin{align}
\abs{f_0(z_t)-f_{\frac{t}{2}}(z_{\frac{t}{2}})} \leq \varphi \frac{N \rho_j }{Q_0 (t, j)} .
\end{align}
\eel
\proof We will just prove the lemma for fixed $t$; we omit a straightforward discretization argument, given in detail in \cite{bourgade2018extreme}, to extend the result to all $t$ simultaneously.

Using the definitions in \eqref{def:z_t, msc} and \eqref{def: f_u and s_u}, Ito's formula gives
\begin{equation}\label{DE: f_u(z_(t-u))}
\begin{aligned}
df_u(z_{t-u}) &= \pare{s_u(z_{t-u})-\msc(z_{t-u})}\partial_zf_u(z_{t-u})+\frac{1}{N}(\frac{2}{\beta}-\frac{1}{2})\partial_{zz}f_u(z_{t-u}) du\\
&-\sqrt{\frac{2}{N\beta}}e^{-u/2}\sum_i\frac{\fu_i(u)}{(x_i(u)-z_{t-u})^2}dB_i(u)
\end{aligned}
\end{equation}
Therefore,
\begin{align}
& \abs{f_0(z_t)-f_{\frac{t}{2}}(z_{\frac{t}{2}})} \lesssim  \abs{\int_0^{t/2}\pare{s_u(z_{t-u})-\msc(z_{t-u})}\e^{-u} \sum_i \frac{\fu_i(u)}{(x_i-z_{t-u})^2} \d u} \label{eqn: short-time increment 1} \\
+ & \abs{\frac{1}{N}\int_0^{t/2}\e^{-u/2}\sum_i \frac{\fu_i(u)}{(x_i-z_{t-u})^3}\d u }\label{eqn: short-time increment 2} \\
+ & \abs{\int_0^{t/2}\sqrt{\frac{2}{N\beta}}\e^{-u/2}\sum_i\frac{\fu_i(u)}{(x_i(u)-z_{t-u})^2}\d B_i}. \label{eqn: short-time increment 3}
\end{align}
We now bound each of the three terms above.

 By \eqref{eqn:char-est}, \eqref{eqn: |s_u(z_t)-msc(z_t)|}, and \eqref{eqn:u-inf-bd} we have that the term \eqref{eqn: short-time increment 1} is bounded by (with overwhelming probability),
\begin{align}
& \varphi^{-C} \abs{\int_0^{t/2}\pare{s_u(z_{t-u})-\msc(z_{t-u})}\e^{-u} \sum_i \frac{\fu_i(u)}{(x_i-z_{t-u})^2} \d u} \notag\\
 \leq &  \frac{1}{N\eta_t}\sum_i\frac{1}{|\gamma_i-\gamma_j|^2+\eta_t^2}\int_0^{t/2}\frac{du}{N|u+\rho_i|} \lesssim   \frac{1}{N^2\eta_t}\sum_i\frac{1}{|\gamma_i-\gamma_j|^2+\eta_t^2} \frac{t}{\rho_i} \notag\\
 \lesssim & \frac{t}{N \eta_t} \int_{\rr} \frac{1}{ (x - \gamma_j)^2 + \eta_t^2} \lesssim \frac{t}{N \eta_t^2} \lesssim \frac{1}{N t (t+\rho_j)^2} = \frac{N \rho_j}{Q_0 (t, j)} \label{eqn:short-time-a1}
\end{align} 
Similarly, \eqref{eqn: short-time increment 2} is bounded by
\begin{align}
& \varphi^{-C} \abs{\frac{1}{N}\int_0^{t/2}\e^{-u/2}\sum_i \frac{\fu_i(u)}{(x_i-z_{t-u})^3}\d u } \lesssim \frac{1}{N^2}\sum_i\frac{1}{|\gamma_i-\gamma_j|^3+\eta_t^3}\frac{t}{ \rho_i} \notag\\
\lesssim &  \frac{1}{N^2\eta_t}\sum_i\frac{1}{|\gamma_i-\gamma_j|^2+\eta_t^2} \frac{t}{\rho_i} \lesssim \frac{ N \rho_i}{Q_0 (t, j)} ,
\end{align}
as the sum in the second line was bounded in \eqref{eqn:short-time-a1}

For \eqref{eqn: short-time increment 3}, introduce the stopping time 
\beq \label{eqn:tau-1-def}
\tau_1 := \inf \left\{ s \in (0, 1) : \exists j : |x_j (s) - \gamma_j | > \frac{ \varphi^{1/100}}{N^{2/3} \hat{i}^{1/3} } \mbox{ or } | \fu_i (s) | > \frac{ \varphi^C}{N ( \rho_i +s ) } \right\}
\eeq
Note that for $t < \tau_1$ we have that \eqref{eqn: |s_u(z_t)-msc(z_t)|} holds. 
With overwhelming probability we have $\tau =1$. The quadratic variation of the martingale inside \eqref{eqn: short-time increment 3} stopped at $\tau_1$ is bounded by,
\begin{align}
& \varphi^{-C} \frac{1}{N}\int_0^{\tau \wedge (t/2) }\sum_i \pare{\frac{|\fu_i(u)|}{|x_i-z_{t-u}|^2}}^2 \d u  \lesssim \frac{1}{N}\sum_i\frac{1}{|\gamma_i-\gamma_j|^4+(\eta_t)^4}\int_0^{t/2}\frac{du}{N^2(u+\rho_i)^2} \d u \notag\\
\lesssim & \frac{t}{N^3 \eta_t^2} \sum_i \frac{1}{ (\gamma_i - \gamma_j)^2 + \eta_t^2} \frac{1}{ \rho_i^2} \lesssim \frac{t}{N^2 \eta_t^3 \sqrt{ \kappa_j + \eta_t}} \lesssim \left( \frac{N \rho_j}{Q_0 (t, j) } \right)^2 .
\end{align}
Above we used \eqref{eqn:misc-10}. Therefore, by Lemma \ref{lem:martingale} we have for \eqref{eqn: short-time increment 3} that with overwhelming probability,
\beq
\abs{\int_0^{t/2}\sqrt{\frac{2}{N\beta}}\e^{-u/2}\sum_i\frac{\fu_i(u)}{(x_i(u)-z_{t-u})^2}\d B_i} \leq \varphi^C \frac{ N \rho_j}{Q_0 (t, j)}
\eeq
This completes the proof. \qed

\subsubsection{Proof of Theorem \ref{thm: character flow}}

The short-time increment is dealt by Lemma \ref{lem: short-time increment}. To deal with the long-time increment $f_t(z)-f_{t/2}(z_{t/2})$, we need a finer estimate. Define a reference function $g_u(z):=e^{-\frac{u}{2}}N s_u(z)$. Ito's formula gives
\begin{align*}
dg_u(z_{t-u}) &= \pare{s_u(z_{t-u})-\msc(z_{t-u})}\partial_zg_u(z_{t-u}) + \frac{1}{N}(\frac{2}{\beta}-\frac{1}{2})\partial_{zz}g_u(z_{t-u})du \\
&-\sqrt{\frac{2}{N\beta}}e^{-u/2}\sum_i\frac{1}{(x_i(u)-z_{t-u})^2}dB_i,
\end{align*}

which is similar to $df_u(z_{t-u})$. It suggests us to consider the difference flow
\begin{align}
\label{eqn: character diff flow 0}&\d f_u(z_{t-u})-\bu_j(u)\d (e^{-u/2}Ns_u(z_{t-u}))\\
\label{eqn: character diff flow 1}&=\pare{s_u(z_{t-u})-\msc(z_{t-u})}\e^{-u/2}\sum_i\frac{u_i-\bu_j}{(x_i-z_{t-u})^2} \d u\\
\label{eqn: character diff flow 2}&+ \frac{1}{N}(\frac{2}{\beta}-\frac{1}{2})2\e^{-u/2}\sum_i\frac{\fu_i-\bu_j}{(x_i-z_{t-u})^3} \d u\\
\label{eqn: character diff flow 3}&-\sqrt{\frac{2}{N\beta}}\e^{-u/2}\sum_i\frac{\fu_i - \bu_j}{(x_i-z_{t-u})^2}\d B_i.
\end{align}
We will integrate each of the terms over the time interval $[t/2, t]$ and bound them separately in the following series of lemmas.

The first term \eqref{eqn: character diff flow 1} is large, but in the following we will see that it is close to \beq
\pare{s_u(z_{t-u})-\msc(z_{t-u})}e^{-u/2}N\partial_t\bu_j. \eeq

\bel
We have, for the term \eqref{eqn: character diff flow 1}, 
\begin{align} \label{eqn: character diff flow 1 bound}
& \left| \int_{t/2}^t\pare{s_u(z_{t-u})-\msc(z_{t-u})}\e^{-u/2}\pare{\sum_i\frac{\fu_i-\bu_j}{(x_i-z_{t-u})^2} - N\partial_t\bu_j}\d u \right| \notag\\
\leq & \varphi^C \left( \frac{1}{ \eta Q_a (t, j)} + \frac{ N \rho_j}{Q_0 (t, j) } \right)
\end{align}
\eel
\proof 
 By \eqref{eqn: |s_u(z_t)-msc(z_t)|},
\begin{align*}
&\int_{t/2}^t\pare{s_u(z_{t-u})-\msc(z_{t-u})}\e^{-u/2}\pare{\sum_i\frac{\fu_i-\bu_j}{(x_i-z_{t-u})^2} - N\partial_t\bu_j}\d u\\
&\lesssim \int_0^{t/2}\frac{1}{N\eta_u}\pare{\sum_i\frac{|\fu_i(t-u)-\bu_i(t-u)|}{|x_i-z_u|^2}+\abs{\sum_i\frac{\bu_i(t-u)-\bu_j(t-u)}{(x_i-z_u)^2}-N\partial_t\bu_j(t-u)}}du\\
&=: I + II.
\end{align*}
Here $\eta_u$ is defined in \eqref{def: eta_u}. We now bound each of the terms $I$ and $II$. 

Using Property $(\pP_a)$ and \eqref{eqn:char-est}, we have with overwhelming probability,
\begin{align} \label{eqn:big-flow-a1}
I
 \leq \varphi^{C_1} \sum_i\left( \frac{1}{NQ_a(t,i)} + \frac{1}{N Q_0 (t, i)} \right) \int_0^{t/2}\frac{1}{\eta_u(|\gamma_i-\gamma_j|^2+\eta_u^2)} \d u.
\end{align}
Then, starting with the $Q_a$ term above,   and using \eqref{eqn:misc-time-1} to bound the time integral in the first inequality,
\begin{align}
& \sum_i\frac{1}{NQ_a(t,i)}\int_0^{t/2}\frac{1}{\eta_u(|\gamma_i-\gamma_j|^2+\eta_u^2)} \d u \notag\\
\lesssim & \sum_i\frac{1}{N\sqrt{\kappa_j+\eta}\pare{|\gamma_i-\gamma_j|^2+\eta^2}}\frac{1}{N(t+\rho_i)(Nt(t+\rho_i)\rho_i)^{1-a}} \notag \\
\lesssim & \frac{1}{\sqrt{\kappa_j+\eta}}\int \frac{\rhosc (y)}{\pare{|y-\gamma_j|^2+\eta^2}N(t+\rhosc (y))\pare{Nt(t+\rhosc (y))\rhosc (y)}^{1-a}} \d y \notag \\
\lesssim &\frac{1}{\eta}\frac{1}{Q_a(t,j)}.
\end{align}
We used \eqref{eqn:misc-11} in the last inequality. By the same argument the $Q_0$ term in \eqref{eqn:big-flow-a1} is bounded by $(\eta Q_0 (t, j) )^{-1} \leq ( \eta Q_a (t, j) )^{-1}$ as we are assuming that $B(t, j) \geq 1$. Therefore,
\beq \label{eqn:big-flow-a6}
I \leq \frac{\varphi^{C_1}}{\eta Q_a (t, j)} .
\eeq 

To deal with $II$, we first replace the factors $(x_i-z_u)$ in the sum by $(\gamma_i-\gamma_j)$:
\begin{align}
&\int_0^{\frac{t}{2}}\frac{1}{N\eta_u}\abs{\sum_i\frac{\bu_i-\bu_j}{(x_i-z_u)^2}-\frac{\bu_i-\bu_j}{(\gamma_i-\gamma_j)^2}} \d u \notag\\
\lesssim & \int_0^{\frac{t}{2}}\frac{1}{N\eta_u}\pare{\abs{\sum_{i:|\gamma_i-\gamma_j|>\eta_u}\frac{\bu_i-\bu_j}{(x_i-z_u)^2}-\frac{\bu_i-\bu_j}{(\gamma_i-\gamma_j)^2}} + \abs{\sum_{i:|\gamma_i-\gamma_j|\leq \eta_u}\frac{\bu_i-\bu_j}{(\gamma_i-\gamma_j)^2}}} \d u \notag \\
\lesssim &  \int_0^{\frac{t}{2}}  \frac{1}{N\eta_u}\pare{\sum_{i:|\gamma_i-\gamma_j|>\eta_u}\frac{|x_i-\gamma_i|+|z_u-\gamma_j|}{|\gamma_i-\gamma_j|^2Nt(t+\rho_i\wedge\rho_j)^2}+\sum_{i:|\gamma_i-\gamma_j|\leq\eta_u}\frac{1}{|\gamma_i-\gamma_j|Nt(t+\rho_i\wedge\rho_j)^2}} \d u \label{eqn:big-flow-a2}
\end{align}
In the last line above, we used regularity of $\bu$ \eqref{eqn:bua-dif}. For the sum inside, the first term is bounded as
\begin{align}
& \varphi^{-1} \sum_{i:|\gamma_i-\gamma_j|>\eta_u}\frac{|x_i-\gamma_i|+|z_u-\gamma_j|}{|\gamma_i-\gamma_j|^2Nt(t+\rho_i\wedge\rho_j)^2} \lesssim \sum_{i:|\gamma_i-\gamma_j|>\eta_u}\frac{\frac{1}{N\rho_i}+\eta_u}{|\gamma_i-\gamma_j|^2Nt(t+\rho_i\wedge\rho_j)^2} \notag \\
&\lesssim \int_{|y-\gamma_j|>\eta_u}\frac{1+N\rhosc (y)\eta_u}{(y-\gamma_j)^2Nt(t+\rhosc (y)\wedge \rho_j)^2} \d y \notag \\
&\lesssim \int_{y<\gamma_j-\eta_u}\frac{1+N\rhosc (y)\eta_u}{(y-\gamma_j)^2Nt(t+\rho_j)^2}dy + \mathbbm{1}_{\kappa_j>\eta_u}\int_{y>\gamma_j+\eta_u}\frac{1+N\rhosc(y)\eta_u}{(y-\gamma_j)^2Nt(t+\rhosc(y))^2}\d y \notag \\
&\lesssim \pare{\frac{1}{\eta_u Nt(t+\rho_j)^2}+\frac{N\sqrt{\kappa_j+\eta_u}}{Nt(t+\rho_j)^2}} + \mathbbm{1}_{\kappa_j>\eta_u}\pare{\frac{1}{\eta_u\wedge \kappa_j}\frac{1}{Nt(t+\rho_j)^2} + \frac{N\sqrt{\kappa_j+\eta_u}}{Nt(t+\rho_j)^2}}\notag \\
&\lesssim \frac{N\sqrt{\kappa_j+\eta_u}}{Nt(t+\rho_j)^2} \label{eqn:big-flow-a3}
\end{align}
In bounding the integrals in the second last inequality, the bounds for the $y < \gamma_j - \eta_u$ integral are easy. For the other integral first do $y < \gamma_{j/2}$ and pull out all the $\rhosc (y) \to \rho_j$ terms. Then for $y > \gamma_{j/2}$ use $(y- \gamma_j) \geq \kappa_j \geq \kappa_j + \eta_u$ and then use \eqref{eqn:misc-9}. 
In the last line, we used that $N\sqrt{\kappa_j+\eta_u}\eta_u \geq N\rho_j\eta \gg 1$ from the assumption $\eta \gg \frac{1}{N\rho_j}$.

Similarly, the second term in \eqref{eqn:big-flow-a2} is bounded as
\begin{align} \label{eqn:big-flow-a4}
\sum_{i:|\gamma_i-\gamma_j|\leq\eta_u}\frac{1}{|\gamma_i-\gamma_j|Nt(t+\rho_i\wedge\rho_j)^2} \leq \varphi \frac{N\sqrt{\kappa_j+\eta_u}}{Nt(t+\rho_j)^2} ,
\end{align}
by using \eqref{eqn:misc-12}. 
Then, using \eqref{eqn:big-flow-a3} and \eqref{eqn:big-flow-a4} we find that \eqref{eqn:big-flow-a2} is bounded by,
\beq \label{eqn:big-flow-a5}
\int_0^{\frac{t}{2}}\frac{1}{N\eta_u}\abs{\sum_i\frac{\bu_i-\bu_j}{(x_i-z_u)^2}-\frac{\bu_i-\bu_j}{(\gamma_i-\gamma_j)^2}} \d u \leq \varphi \frac{1}{Nt(t+\rho_j)^2}\int_0^{\frac{t}{2}}\frac{\sqrt{\kappa_j+\eta_u}}{\eta_u}du \leq \varphi^2 \frac{1}{N t (t + \rho_j)^2}
\eeq
where we used \eqref{eqn:misc-time-2}. By Proposition \ref{prop:deltbu}, the rest of the term $II$ can be bounded by,
\begin{align}
\int_0^{\frac{t}{2}}\frac{1}{N\eta_u}\abs{\sum_i\frac{\bu_i-\bu_j}{(\gamma_i-\gamma_j)^2}-N\partial_t\bu_j}\d u 
=\int_0^{\frac{t}{2}}\frac{1}{N\eta_u}\abs{\sum_i\frac{\bu_i-\bu_j}{(\gamma_i-\gamma_j)^2}-N\bar{\sB}\bu_j}\d u.
\end{align}
The sum inside the absolute value can be bounded as
\begin{align}
&\abs{\sum_i\frac{\bu_i-\bu_j}{(\gamma_i-\gamma_j)^2}-N\bar{\sB}\bu_j} \notag \\
&\leq \sum_{i:i\neq j}\abs{\frac{\bu_i-\bu_j}{(\gamma_i-\gamma_j)^2} -N\int_{\gamma_{i-\frac{1}{2}}}^{\gamma_{i+\frac{1}{2}}}\frac{\bu(y)-\bu_j}{(y-\gamma_j)^2}\rhosc(y)\d y} +N\int_{\gamma_{j+\frac{1}{2}}}^{\gamma_{j-\frac{1}{2}}}\frac{|\bu(y)-\bu_j|}{(y-\gamma_j)^2}\rhosc(y)\d y \notag \\
&\lesssim \sum_{i:i\neq j}\pare{\frac{1}{N\rho_j}+\frac{1}{N\rho_i}}\frac{|\bu_i-\bu_j|}{(\gamma_i-\gamma_j)^3}+ \sum_{i \neq j } \sup_{y \in [ \gamma_i,\gamma_{i+1} ] } \frac{ | \bu (y) - \bu_i|}{(\gamma_i - \gamma_j)^2} +N\int_{\gamma_{j+\frac{1}{2}}}^{\gamma_{j-\frac{1}{2}}}\frac{|\bu(y)-\bu_j|}{(y-\gamma_j)^2}\rhosc(y)\d y \notag \\
&\lesssim  \varphi \frac{N^2\kappa_j}{Q_0(t,j)}.
\end{align}
To get the last line, we applied \eqref{est:|Lu-d_t(u_bar)|,part1} and \eqref{est:|Lu-d_t(u_bar)|,part2} with parameter $\ell = \frac{1}{2}$.  Using \eqref{eqn:misc-time-3} to bound the time integral we find,
\begin{align}
\int_0^{\frac{t}{2}}\frac{1}{N\eta_u}\abs{\sum_i\frac{\bu_i-\bu_j}{(\gamma_i-\gamma_j)^2}-N\partial_t\bu_j}du \lesssim \varphi \frac{\rho_j}{\sqrt{\kappa_j+\eta}}\frac{1}{Nt(t+\rho_j)^2}.
\end{align}
Combining the above with \eqref{eqn:big-flow-a5} we find,
\begin{align}
II \leq \varphi^C \frac{1}{Nt(t+\rho_j)^2}.
\end{align}
We conclude the lemma from the above inequality and \eqref{eqn:big-flow-a6}. \qed

\bel
We have,
\begin{align} \label{eqn: character diff flow 2 bound}
& \abs{\int_{\frac{t}{2}}^t \frac{1}{N}(\frac{2}{\beta}-\frac{1}{2})2\e^{-u/2}\sum_i\frac{\fu_i-\bu_j}{(x_i-z_{t-u})^3}\d u}  \notag\\
\leq & \varphi^C \left(  \frac{1}{\eta Q_a(t,j)} + \frac{N\rho_j}{Q_0(t,j)} \right) .
\end{align}
\eel
\proof 
For the second term \eqref{eqn: character diff flow 2}, we have
\begin{align} \label{eqn:big-flow-a7}
\abs{\int_{\frac{t}{2}}^t \frac{1}{N}(\frac{2}{\beta}-\frac{1}{2})2\e^{-u/2}\sum_i\frac{\fu_i-\bu_j}{(x_i-z_{t-u})^3}\d u} \lesssim \frac{1}{N}\int_0^t\sum_i\frac{|\fu_i(t)-\bu_i(t)|+|\bu_i(t)-\bu_j(t)|}{|\gamma_i-\gamma_j|^3+\eta_u^3} \d u
\end{align}
Let us  deal with the first sum above. Using property $\pP_a$,  
\begin{align}
\frac{\varphi^{-C}}{N}\sum_i\frac{|u_i(t)-\bu_i(t)|}{|\gamma_i-\gamma_j|^3+\eta_u^3} &\lesssim \int \frac{\rhosc(y)}{|y-\gamma_j|^3+\eta_u^3}\frac{1}{N(t+\rhosc(y))\pare{Nt(t+\rhosc(y))\rhosc(y)}^{1-a}}\notag\\
+ &  \int \frac{\rhosc(y)}{|y-\gamma_j|^3+\eta_u^3}\frac{1}{N(t+\rhosc(y))\pare{Nt(t+\rhosc(y))\rhosc(y)}^{1}} \notag\\
&\lesssim \frac{\sqrt{\kappa_j+\eta_u}}{\eta_u^2}\frac{1}{Q_a(t,j)}.
\end{align}
Above we used \eqref{eqn:misc-13} and the fact that $B(t, j) \geq 1$. 

For the other sum in \eqref{eqn:big-flow-a7}, using \eqref{eqn:bua-dif}, 
\begin{align}
& \varphi^{-1} \frac{1}{N}\sum_i\frac{|\bu_i(t)-\bu_j(t)|}{|\gamma_i-\gamma_j|^3+\eta_u^3} \lesssim \int \frac{\rhosc(y)}{|y-\gamma_j|^3+\eta_u^3}\frac{|y-\gamma_j|}{Nt(t+\rho(y)\wedge \rho_j)^2} \notag\\
\lesssim & \int_{|y-\gamma_j|<\kappa_j/2}\frac{|y-\gamma_j|dy}{|y-\gamma_j|^3+\eta_u^3}\frac{\rho_j}{Nt(t+\rho_j)^2} + \int_{y>\gamma_j +\frac{\kappa_j}{2}}\frac{\rhosc (y)}{\kappa_j^3+\eta_u^3}\frac{\kappa_j \d y}{Nt(t+\rhosc (y))^2}\\
+& \int_{y<\gamma_j-\kappa_j/2}\frac{\kappa(y)}{\kappa(y)^3+\eta_u^3}\frac{\rhosc (y)}{Nt(t+\rho_j)^2} \d y \notag\\
\lesssim & \left( \frac{ N \kappa_j}{\eta_u} + \frac{ N \kappa_j^3}{\kappa_j^3+\eta_u^3} + \frac{ N \rho_j}{ \sqrt{ \kappa_j + \eta_u}} \right) \frac{1}{Q_0 (t, j)} \lesssim  \pare{\frac{N\kappa_j}{\eta_u}+\frac{N\rho_j}{\sqrt{\kappa_j+\eta_u}}}\frac{1}{Q_0(t,j)}
\end{align}
Above we used \eqref{eqn:misc-5a} to bound the second integral in the second line. Therefore,
\begin{align}
& \abs{\int_{\frac{t}{2}}^t \frac{1}{N}(\frac{2}{\beta}-\frac{1}{2})2\e^{-u/2}\sum_i\frac{\fu_i-\bu_j}{(x_i-z_{t-u})^3}\d u} \notag\\
\lesssim & \varphi^C \int_0^{t/2} \left( \frac{ \sqrt{ \kappa_j + \eta_u}}{\eta_u^2 Q_a (t, j) } + \left( \frac{N \kappa_j}{\eta_u} + \frac{ N \rho_j}{ \sqrt{ \kappa_j + \eta_u} } \right) \frac{1}{Q_0 (t, j) } \right)\d u\notag\\
\lesssim & \varphi^C \left( \frac{1}{\eta Q_a (t, j) } + \frac{ N \rho_j}{ Q_0 (t ,j) }  \right)
\end{align}
Above we used \eqref{eqn:misc-time-3}, \eqref{eqn:misc-time-4} and \eqref{eqn:misc-time-5}. This completes the proof. \qed

\bel
With overwhelming probability we have,
\begin{align} \label{eqn: character diff flow 3 bound}
& \left| \sqrt{\frac{2}{N\beta}}\e^{-u/2}\sum_i\frac{\fu_i - \bu_j}{(x_i-z_{t-u})^2}\d B_i \right| \notag\\
\leq & \varphi^{C_1} \left(  \frac{1}{\eta Q_a (t,j)} + \frac{ N \rho_j}{Q_0 (t, j)} \right).
\end{align}
\eel
\proof Let $\tau_2 >0$ be stopping time,
\beq
\tau_2 := \inf \{ u \in [0, 1] : \exists k : | \fu_k (u) - \bu_k (u) | > \varphi^{C_1} ( Q_a (u, k)^{-1} + Q_0 (u, k)^{-1} )\}.
\eeq
As we are assuming Property $\pP_a$ holds, we have $\tau_2 =t$ with overwhelming probability for some $C_1 >0$ sufficiently large. Let $\tau_1$ be the stopping time from \eqref{eqn:tau-1-def} and set
\beq
\tau = \tau_1 \wedge \tau_2.
\eeq
Then $\tau = 1$ with overwhelming probability.  Further, note that for $t < \tau$ we have that \eqref{eqn: |s_u(z_t)-msc(z_t)|} holds.

The quadratic variation of the martingale term on the LHS of \eqref{eqn: character diff flow 3 bound} is bounded by,
\beq \label{eqn:big-flow-a8}
\int_{\frac{t}{2}}^t\frac{2\e^{-u}}{N\beta}\sum_i\frac{|\fu_i-\bu_j|^2}{|x_i-z_{t-u}|^4}\d u \lesssim \frac{1}{N}\int_0^{\frac{t}{2}}\sum_i\frac{|\fu_i(u)-\bu_i(u)|^2+|\bu_i(u)-\bu_j(u)|^2}{|\gamma_i-\gamma_j|^4+\eta_u^4}\d u
\eeq
We begin by estimating the two sums inside the integral. We first have, for $u < \tau$, 
\begin{align}
& \varphi^{-C_1} \frac{1}{N}\sum_i \frac{|\fu_i(u)-\bu_i(u)|^2}{|\gamma_i-\gamma_j|^4+\eta_u^4}\notag \\
\lesssim & \int \frac{\rhosc(y)}{|y-\gamma_j|^4+\eta_u^4}\pare{\frac{1}{N(t+\rhosc(y))\pare{Nt(t+\rhosc(y))\rhosc(y)}^{1-a}}}^2 \notag \\
+ & \int \frac{\rhosc(y)}{|y-\gamma_j|^4+\eta_u^4}\pare{\frac{1}{N(t+\rhosc(y))\pare{Nt(t+\rhosc(y))\rhosc(y)}}}^2 
\end{align}
The first integral can be bounded by,

\begin{align}
 &\int \frac{\rhosc(y)}{|y-\gamma_j|^4+\eta_u^4}\pare{\frac{1}{N(t+\rhosc(y))\pare{Nt(t+\rhosc(y))\rhosc(y)}^{1-a}}}^2 \notag\\
\lesssim & \int_{|y-\gamma_j|<\kappa_j/2}\frac{\rho_j}{|y-\gamma_j|^4+\eta_u^4}\frac{1}{Q_a(t,j)^2}dy \notag\\
+& \int_{y>\gamma_j +\frac{\kappa_j}{2}}\frac{\rhosc(y)}{\kappa_j^4+\eta_u^4}\pare{\frac{1}{N(t+\rhosc(y))\pare{Nt(t+\rhosc(y))\rhosc(y)}^{1-a}}}^2 \d y \notag\\
+&\int_{y<\gamma_j-\kappa_j/2}\frac{\rhosc(y)}{\kappa(y)^4+\eta_u^4}\pare{\frac{1}{N(t+\rhosc(y))\pare{Nt(t+\rhosc(y))\rhosc(y)}^{1-a}}}^2 \d y \notag \\
&\lesssim \frac{\rho_j\kappa_j}{\eta_u^3(\eta_u+\kappa_j)}\frac{1}{Q_a(t,j)^2} + \mathbbm{1}_{t<\rho_j}\frac{t^3}{\kappa_j^4+\eta_u^4}\pare{\frac{1}{Nt(Nt^3)^{1-a}}}^2 +\frac{1}{(\kappa_j+\eta_u)^{\frac{5}{2}}}\frac{1}{Q_a(t,j)^2} \notag  \\
&\lesssim \frac{\sqrt{\kappa_j+\eta_u}}{\eta_u^3}\frac{1}{Q_a(t,j)^2} + \mathbbm{1}_{t<\rho_j}\frac{t^3}{\kappa_j^4+\eta_u^4}\pare{\frac{1}{Nt(Nt^3)^{1-a}}}^2. \label{eqn:big-flow-a9}
\end{align}
Above, the estimates for the integrals for $y < \gamma_j + \kappa_j/2$ are straightforward (the first and third integrals). For the integral with $y > \gamma_{j/2}$ we use \eqref{eqn:misc-14}. Combining this with the analogous argument when $a=0$ we therefore find,
\begin{align} \label{eqn:big-flow-a10}
\frac{1}{N}\sum_i \frac{|\fu_i(u)-\bu_i(u)|^2}{|\gamma_i-\gamma_j|^4+\eta_u^4} \leq &  \varphi^{C_1} \left( \frac{\sqrt{\kappa_j+\eta_u}}{\eta_u^3}\frac{1}{Q_a(t,j)^2} + \mathbbm{1}_{t<\rho_j}\frac{t^3}{\kappa_j^4+\eta_u^4}\pare{\frac{1}{Nt(Nt^3)^{1-a}}}^2 \right)\notag\\
+&  \varphi^{C_1} \1_{t<\rho_j} \frac{t^3}{ \kappa_j^4+\eta_u^4} \left( \frac{1}{Nt (N t^3) } \right)^2
\end{align}


For the other term in \eqref{eqn:big-flow-a8} we have, using \eqref{eqn:bua-dif}
\begin{align}
&\frac{1}{N}\sum_i\frac{|\bu_i(u)-\bu_j(u)|^2}{|\gamma_i-\gamma_j|^4+\eta_u^4} \notag \\
&\lesssim \int \frac{\rhosc(y)}{|y-\gamma_j|^4+\eta_u^4}\pare{\frac{|y-\gamma_j|}{Nt(t+\rhosc(y)\wedge \rho_j)^2}}^2 \d y \notag \\
&\lesssim \int_{|y-\gamma_j|<\kappa_j/2}\frac{|y-\gamma_j|^2}{|y-\gamma_j|^4+\eta_u^4}\frac{\rho_j}{(Nt(t+\rho_j)^2)^2} \d y + \int_{y>\gamma_j +\frac{\kappa_j}{2}}\frac{\kappa_j^2}{\kappa_j^4+\eta_u^4}\frac{\rhosc(y)}{(Nt(t+\rhosc (y))^2)^2}\d y \notag\\
&\quad +\int_{y<\gamma_j-\kappa_j/2}\frac{\kappa(y)^2\rho(y)}{\eta_u^4+\kappa(y)^4}\frac{1}{(Nt(t+\rho_j)^2)^2}\d y \notag\\
&\lesssim \frac{N^2\rho_j^9}{\eta_u(\eta_u^3+\kappa_j^3)}\frac{1}{Q_0(t,j)^2} + \pare{\frac{N^2\rho_j^9}{\kappa_j^4+\eta_u^4}\frac{1}{Q_0(t,j)^2} + \mathbbm{1}_{t<\rho_j}\frac{t^3\kappa_j^2}{\kappa_j^4+\eta_u^4}\frac{1}{(Nt^3)^2}} + \frac{N^2\rho_j^2}{\sqrt{\kappa_j+\eta_u}}\frac{1}{Q_0(t,j)^2} \notag \\
&\lesssim \frac{N^2\rho_j^9}{\eta_u(\eta_u^3+\kappa_j^3)}\frac{1}{Q_0(t,j)^2} + \mathbbm{1}_{t<\rho_j}\frac{t^3\kappa_j^2}{\kappa_j^4+\eta_u^4}\frac{1}{(Nt^3)^2} + \frac{N^2\rho_j^2}{\sqrt{\kappa_j+\eta_u}}\frac{1}{Q_0(t,j)^2}.
\end{align}
Above, the integrals for $y < \gamma_{j/2}$ are straightforward.  For the integral with $y > \gamma_{j/2}$ we use \eqref{eqn:misc-15}. 

For the quadratic variation of the stopped process we therefore have,
\begin{align} \label{eqn:big-flow-a11}
& \int_{\frac{t}{2} \wedge \tau}^{t \wedge \tau}\frac{2\e^{-u}}{N\beta}\sum_i\frac{|\fu_i-\bu_j|^2}{|x_i-z_{t-u}|^4}\d u \notag\\
\leq & \varphi^{C_1} \int_0^{t/2}   \left( \frac{\sqrt{\kappa_j+\eta_u}}{\eta_u^3}\frac{1}{Q_a(t,j)^2} + \mathbbm{1}_{t<\rho_j}\frac{t^3}{\kappa_j^4+\eta_u^4}\pare{\frac{1}{Nt(Nt^3)^{1-a}}}^2 \right) \d u \notag\\
+ & \varphi^{C_1}  \int_0^{t/2} \left(  \frac{N^2\rho_j^9}{\eta_u(\eta_u^3+\kappa_j^3)}\frac{1}{Q_0(t,j)^2} + \mathbbm{1}_{t<\rho_j}\frac{t^3\kappa_j^2}{\kappa_j^4+\eta_u^4}\frac{1}{(Nt^3)^2} + \frac{N^2\rho_j^2}{\sqrt{\kappa_j+\eta_u}}\frac{1}{Q_0(t,j)^2} \right) \d u \notag\\
\leq & \varphi^{C_1} \left( \frac{1}{\eta Q_a (t, j)} + \mathbbm{1}_{t<\rho_j}\frac{1}{(t\rho_j)(Nt\kappa_j)^a Q_a(t, j)} \right)^2 \notag\\
+ & \varphi^{C_1} \left( \frac{N^2 \rho_j^9}{(\kappa_j + \eta)^{7/2} Q_0 (t, j)^2} + \1_{t < \rho_j} \frac{t^4}{\kappa_j^4 (N t^3)^2} + \frac{N^2 \rho_j^2}{Q_0 (t, j)^2} \right) \notag\\
\leq & \varphi^{C_1} \left( \frac{1}{\eta Q_a (t, j)} + \frac{ N \rho_j}{Q_0 (t, j) } \right)^2.
\end{align}
In the first inequality we absorbed the term on the second line of \eqref{eqn:big-flow-a10}  into the second term on the third line of \eqref{eqn:big-flow-a11} (as if $t \leq \rho_j$ then the assumption $B(t, j) \geq 1$ implies that $(N t)^{-1} \leq \kappa_j$). In the second inequality we used \eqref{eqn:misc-time-6}, \eqref{eqn:misc-time-7}, \eqref{eqn:misc-time-8} and \eqref{eqn:misc-time-9}. We also used
\begin{align}
\1_{t < \rho_j}\int_0^t \frac{t^3}{ \kappa_j^4 + \eta_u^4} \left( \frac{1}{Nt (N t^3)^{1-a}} \right)^2 = \1_{t < \rho_j} \left(\frac{( N t^3)^a}{\rho_j^4 N^2 t^2} \right)^2
\end{align}
If $t < N^{-1/3}$ then the above term is bounded by
\beq
\1_{t < \rho_j} \left( \frac{1}{ t \rho_j (N t \kappa_j)^a Q_a (t, j)} \right)^2 = \1_{t < \rho_j} \left( \frac{1}{ t \rho_j (N t \rho_j^2 )^a N \rho_j (N t \rho_j^2)^{1-a} } \right)^2 = \1_{t < \rho_j} \left( \frac{1}{N^2 t^2 \rho_j^4} \right)^2
\eeq
If $N^{-1/3} \leq t \leq \rho_j$ then instead it is bounded by
\beq
\1_{t < \rho_j} \left( \frac{ N \rho_j}{Q_0 (t, j)} \right)^2 = \1_{t < \rho_j} \left(\frac{1}{N \rho_j^2 t} \right)^2 \geq \1_{t < \rho_j} \left( \frac{ N t^3}{N^2 t^2 \rho_j^4} \right)^2 \geq \1_{t < \rho_j} \left( \frac{ (N t^3)^a}{N^2 t^2 \rho_j^4} \right)^2.
\eeq
The first inequality uses $t \leq \rho_j$ and the second that $N t^3 \geq 1$. 

In the last inequality in \eqref{eqn:big-flow-a11} we absorbed the $\1_{t < \rho_j} Q_a^{-1}$ term into the $\eta^{-1} Q_a^{-1}$ term (as $t< \rho_j$ and $B(t, j) \geq$ implies $Nt \kappa_j \geq 1$ and the assumption $\eta \leq t^2 + t \rho_j$ gives $(t \rho_j)^{-1} \leq \eta^{-1}$). We also used that
\beq
\1_{t < \rho_j}\left( \frac{1}{Nt \rho_j^2} \right)^2 = \1_{t < \rho_j} \left( \frac{N \rho_j}{Q_0 (t, j) } \right)^2 
\eeq
With \eqref{eqn:big-flow-a11} in hand, the claim follows from Lemma \ref{lem:martingale}. \qed

\vspace{5 pt}

\noindent{\bf Proof of Theorem \ref{thm: character flow}}.  Integrating \eqref{eqn: character diff flow 0} over time $[t/2,t]$, and using \eqref{eqn: character diff flow 1 bound},\eqref{eqn: character diff flow 2 bound},\eqref{eqn: character diff flow 3 bound}, we get with overwhelming probability, 
\begin{align}
& \pare{f_t(z)-f_{t/2}(z_{t/2})} - \int_{\frac{t}{2}}^t\bu_j(u)\d\pare{\e^{-\frac{u}{2}}Ns_u(z_{t-u})} \notag\\
= & \int_{\frac{t}{2}}^t\pare{s_u(z_{t-u})-\msc{(z_{t-u})}}e^{-\frac{u}{2}}N\partial_u\bu_j(u) \d u\notag\\
+&\varphi^{C_1} \O \pare{ \frac{1}{\eta Q_a (t, j) } +\frac{N \rho_j}{Q_0(t,j)}}.
\end{align}

By definition,
\begin{align}
& \int_{\frac{t}{2}}^t\bu_j(u)\d\pare{\e^{-\frac{u}{2}}Ns_u(z_{t-u})} + \int_{\frac{t}{2}}^ts_u(z_{t-u})\e^{-\frac{u}{2}}N\partial_t\bu_j(u)du \notag\\
= & \e^{-t/2} s_t (z) N \bu_j (t) s_u(z_{t-u}) - \e^{-t/4} s_{t/2} (z_{t/2} ) N \bu_j (t/2)
\end{align}
Moreover, using that $\msc(z_{t-u}) = \e^{-\frac{t-u}{2}} \msc (z) $, we have
\begin{align}
 &\int_{\frac{t}{2}}^t\msc(z_{t-u})e^{-\frac{u}{2}}N\partial_u\bu_j(u) \d u =  \msc(z)\e^{-\frac{t}{2}}N\int_{\frac{t}{2}}^t\partial_u\bu_j(u) \d u \notag\\
 = & \e^{-t/2} \msc (z) N \bu_j (t)  - \e^{-t/4} \msc (z_{t/2} ) N \bu_j (t/2)
\end{align}

%

Applying \eqref{eqn: |s_u(z_t)-msc(z_t)|} and \eqref{eqn:buk-bound}, the difference of the $u=\frac{t}{2}$ terms in above two lines is bounded as 
\begin{align}
\abs{(s_u(z_{t-u})-\msc(z_{t-u})e^{-\frac{u}{2}}N\bu_j(u)\bigg\vert_{u=\frac{t}{2}}} \leq \frac{\varphi}{\eta_{\frac{t}{2}}N(\frac{t}{2}+\rho_j)} 
\lesssim \varphi \frac{N\rho_j}{Q_0(t,j)}.
\end{align}
Combining with Lemma \ref{lem: short-time increment}, we get the desired estimate. \qed

\subsection{Other lemma for averages}

\bel \label{lem: ave mod dyn, part 2}
Assume Property $(\pP_a)$. The following event holds with overwhelming probability. Suppose
\beq
\eta \geq \varphi \frac{1}{N\rho_j}, \quad \eta \leq \varphi^{-1}( t^2+t\rho_j), \quad \ell\geq \varphi N\eta\sqrt{\kappa_{k\vee R}+\eta}.
\eeq
Assume also 
\beq
B(t, k) \gg \varphi^{100}
\eeq
The term 
\begin{align} \label{eqn:big-flow-b}
\abs{\frac{1}{N}\im \sum_{i:|j-i|<\ell}\frac{\pare{\Av(s)\fu(s)}_i-\bu_k(s)}{x_i-z}}
\end{align}
is bounded as follows. If $k \in \Db$ we have that it is bounded above by (up to a $\varphi^C$ factor)
\begin{align}
 & \frac{1}{Q_a (t, k) } \left( \frac{1}{N \eta} + \frac{ N \eta \kappa_k }{\ell} \right)   + \frac{1}{Q_0 (t, k)} \left( R \rho_k + 
 \frac{\ell^2}{N \eta} \right)
\end{align}
and if $k \in \De$ we have that it is bounded above by (up to a $\varphi^C$ factor)
\beq 
\frac{1}{Q_a (t, j) } \left( \frac{1}{N \eta} + \frac{N \eta \kappa_{j \vee \ell}}{\ell} \right) + \frac{1}{Q_0 (t, j \wedge k) } \left( R\sqrt{ \kappa_{j \vee \ell} + \eta} + 
 + \frac{ \ell^2}{N \eta} \right)
\eeq
\eel
We begin with some preparation. By the definition \eqref{def: average operator}, the averaging operator $\Av$ can be written using a (Lipschitz) coefficient function $a$:
\begin{align}\label{def: a_j}
(\Av(t)w)_i = a_i w_i + (1-a_i)\bu_k(t) \mbox{ and } |a_i-a_j|\leq \frac{|i-j|}{R}.
\end{align}
Use this notation, we have the following decomposition
\begin{align}
&\nonumber \frac{1}{N}\im \sum_{i:|j-i|<\ell}\frac{\pare{\Av(s)\fu (s)}_i-\bu_k(s)}{x_i-z}\\
\label{eqn: ave mod dyn, part 2-1}&=\frac{a_j}{N}\im\sum_i\frac{\fu_i(s)-\bu_k(s)}{x_i-z}\\
\label{eqn: ave mod dyn, part 2-2}&- \frac{a_j}{N}\im\sum_{i:|i-j|>\ell}\frac{\fu_i(s)-\bu_i(s)}{x_i-z}\\
\label{eqn: ave mod dyn, part 2-3}&- \frac{a_j}{N}\im\sum_{i:|i-j|>\ell}\frac{\bu_i(s)-\bu_k(s)}{x_i-z}\\
\label{eqn: ave mod dyn, part 2-4}&+\frac{1}{N}\im\sum_{i:|j-i|\leq \ell}\frac{(a_i-a_j)(\fu_i(s)-\bu_k(s))}{x_i-z}
\end{align}

In the computations below, we will frequently use that $|x_i-z| \asymp |\gamma_i-\gamma_j|+\eta$ which is ensured by the assumption $\eta \geq \frac{\varphi}{N\rho_j}$. Note also that by Lemma \ref{lem:index-properties} we have $B(t, j) \geq \varphi^{C}$.  We will bound the terms above in the next series of lemmas. 
\bel
The term \eqref{eqn: ave mod dyn, part 2-1} 
\begin{align} \label{eqn: ave mod dyn, part 2-1 bound}
\frac{a_j}{N}\im\sum_i\frac{\fu_i(s)-\bu_k(s)}{x_i-z} &= a_j  \im\pare{\frac{e^{s/2}}{N}f_0(z_s)-\msc(z)\bu_j(s)} \notag\\
 +& \varphi^{C_1} \O \left( \frac{1}{N \eta Q_a (t, j) } +  \frac{ \sqrt{\kappa_j + \eta} R}{Q_0 (t, j \wedge k)} \right)
\end{align}
\eel
\proof We can write \eqref{eqn: ave mod dyn, part 2-1} as
\beq
\frac{a_j}{N}\im\sum_i\frac{\fu_i(s)-\bu_k(s)}{x_i-z} = a_j\im\pare{\frac{e^{s/2}}{N}f_s(z)-s_s(z)\bu_k(s)}.
\eeq
We have $\abs{\im s_s(z)}\lesssim  \sqrt{\kappa_j+\eta}$ with overwhelming probability. Using regularity of $\bu$ in \eqref{eqn:bua-dif}, we therefore have,
\beq
\abs{\im \pare{s_s(z)\bu_k(s) - s_s(z)\bu_j(s)}} \leq \varphi^C \frac{\sqrt{\kappa_j+\eta}|\gamma_j-\gamma_k|}{Nt(t+\rho_{k\wedge j})^2}\lesssim \varphi^C \sqrt{\kappa_j+\eta}\frac{R}{Q_0(t,j\wedge k)}.
\eeq
Combining this with the estimate about the characteristic flow (Theorem \ref{thm: character flow}), we have, 
\begin{align}
& \im\pare{\frac{e^{s/2}}{N}f_s(z)-s_s(z)\bu_k(s)} = \im\pare{\frac{e^{s/2}}{N}f_0(z_s)-\msc(z)\bu_j(s)+s_s(z)\pare{\bu_j(s)-\bu_k(s)}} \notag\\
+ & \varphi^{C_1} \O \left( \frac{1}{N \eta Q_a (t, j) } + \frac{ \rho_j}{Q_0 (t, j)} \right) \notag\\
= & \im\pare{\frac{e^{s/2}}{N}f_0(z_s)-\msc(z)\bu_j(s)} + \varphi^{C_1} \O \left( \frac{1}{N \eta Q_a (t, j) } + \frac{ \sqrt{\kappa_j + \eta} R}{Q_0 (t, j \wedge k)} \right)
\end{align}
This completes the proof. \qed

\bel
For the term \eqref{eqn: ave mod dyn, part 2-2} we have,
\begin{align}
\label{eqn: ave mod dyn, part 2-2 bound}
\left| \frac{a_j}{N}\im\sum_{i:|i-j|>\ell}\frac{\fu_i(s)-\bu_i(s)}{x_i-z} \right| \leq \varphi^{C_1} \left( \frac{\eta N \kappa_{j \vee \ell}}{\ell} \frac{1}{Q_a (t, j)} \right) .
\end{align}
\eel
\proof 
For the second term \eqref{eqn: ave mod dyn, part 2-2}.  We divide the sum into two part: for $i> j/2$, $|\gamma_i-\gamma_j|\gtrsim \frac{|i-j|}{N\rho_i}$. So 
\begin{align}
& \left| \im\sum_{i>j/2, |i-j|>l}\frac{\fu_i(s)-\bu_i(s)}{x_i(s)-z} \right|
\leq \varphi^{C}  \sum_{i>j/2, |i-j|>l}\frac{\eta}{|\gamma_i-\gamma_j|^2+\eta^2}\left( \frac{1}{Q_a(s,i)}+ \frac{1}{Q_0 (s, i)} \right) \notag \\
\leq &  \varphi^C \sum_{i>j/2, |i-j|>l}\frac{\eta}{\frac{|i-j|^2}{N^2\kappa_i}}\frac{1}{Q_a(s,j)} 
\leq  \varphi^C \frac{\eta N^2\kappa_{j\vee \ell}}{\ell} \frac{1}{Q_a(s,j)} 
\end{align}
In the first line we applied Property $(\pP_a)$. 
The sum over $i\leq j/2$ will only occur when $j>\ell$. In this case, $|\gamma_i-\gamma_j|\asymp \kappa_j$, so 
\begin{align}
& \left| \im\sum_{i<j/2}\frac{\fu_i(s)-\bu_i(s)}{x_i(s)-z} \right| \leq \sum_{i<j/2}\frac{\eta|\fu_i(s)-\bu_i(s)|}{\kappa_j^2}\notag\\
\leq & \varphi^C \frac{\eta}{\kappa_j^2}\sum_{i<j/2}\left( \frac{1}{N(s+\rho_i)(Ns(s+\rho_i)\rho_i)^{1-a}}+ \frac{1}{N(s+\rho_i)(Ns(s+\rho_i)\rho_i)} \right) \notag \\
\leq & \varphi^C \frac{\eta}{\kappa_j^2}\frac{j}{Q_a(s,j)},
\end{align}
which is smaller than the previous sum as $j \geq \ell$. Above, we used \eqref{eqn:misc-5a}. We conclude the proof. \qed  

\bel
For the term \eqref{eqn: ave mod dyn, part 2-3} we have
\begin{align} \label{eqn: ave mod dyn, part 2-3 bound}
- & \frac{a_j}{N}\im\sum_{i:|i-j|>\ell}\frac{\bu_i(s)-\bu_k(s)}{x_i-z} = -a_j \eta \del_s \bu_k (s) +  \varphi^{C_1} \O \left(  \frac{ N \eta \kappa_{j \vee \ell} + R \rho_{j \vee \ell}}{Q_0 (t, j \wedge k)} \right)
\end{align}
\eel
\proof 
For \eqref{eqn: ave mod dyn, part 2-3}, we compare it with a deterministic term $\eta\partial_s\bu_j(s)$ (which equals to $\eta\bar{\sB}\bu_j$ by Proposition \ref{prop:deltbu}).
\begin{align} \label{eqn:big-ave-a1}
&\abs{\frac{1}{N} \im\sum_{i:|i-j|>l}\frac{\bu_i-\bu_k}{x_i-z} - \eta\partial_s\bu_j(s)} \notag\\
\leq & \frac{1}{N} \abs{\im\sum_{i:|i-j|>l}\frac{\bu_i-\bu_k}{x_i-z}-\eta\sum_{i:|i-j|>l}\frac{\bu_i-\bu_k}{|\gamma_i-\gamma_j|^2}} + \eta\abs{\bar{\sB}\bu_j - \frac{1}{N}\sum_{i:|i-j|>l}\frac{\bu_i-\bu_k}{|\gamma_i-\gamma_j|^2}}
\end{align}
Using regularity of $\bu$ in \eqref{eqn:bua-dif}, and then \eqref{est:|Lu-d_t(u_bar)|,part1}, \eqref{est:|Lu-d_t(u_bar)|,part2} and \eqref{est:|B_bar u_bar_j - B_bar u_bar_j|} in the final inequality we obtain 
\begin{align}
&\abs{\bar{\sB}\bu_j -  \frac{1}{N}\sum_{i:|i-j|>l}\frac{\bu_i-\bu_k}{|\gamma_i-\gamma_j|^2}} \notag\\
\leq & \abs{\bar{\sB}\bu_j - \frac{1}{N}\sum_{i:|i-j|>\ell}\frac{\bu_i-\bu_j}{|\gamma_i-\gamma_j|^2}} + \frac{1}{N}\sum_{i:|i-j|>\ell}\frac{|\bu_j-\bu_k|}{|\gamma_i-\gamma_j|^2} \notag\\
\leq & \sum_{i:|i-j|> \ell}\abs{\int_{\gamma_i}^{\gamma_{i+1}}\frac{\bu(y)-\bu_j}{(y-\gamma_j)^2}\rhosc(y)\d y - \frac{\bu_i-\bu_j}{N|\gamma_i-\gamma_j|^2}}+ \int_{\gamma_{j+\ell}}^{\gamma_{j-\ell}}\frac{|\bu(y)-\bu_j|}{(y-\gamma_j)^2}\rhosc(y) \d y \notag\\
+ & \frac{1}{N}\sum_{i:|i-j|>\ell}\frac{|\bu_j-\bu_k|}{|\gamma_i-\gamma_j|^2} \notag\\
\lesssim & \frac{1}{N}\sum_{i:|i-j|> \ell}\pare{\pare{\frac{1}{N\rho_j}+\frac{1}{N\rho_i}}\frac{|\bu_i-\bu_j|}{(\gamma_i-\gamma_j)^3} + \sup_{y \in [\gamma_i, \gamma_{i+1}] } \frac{ |\bu (y) - \bu_i|}{(\gamma_i - \gamma_j)^2} +\frac{|\bu_j-\bu_k|}{|\gamma_i-\gamma_j|^2}} \notag\\
+&  \int_{\gamma_{j+\ell}}^{\gamma_{j-\ell}}\frac{|\bu(y)-\bu_j|}{(y-\gamma_j)^2}\rhosc(y)\d y \leq  \varphi^C  \frac{N\kappa_{j\vee \ell}}{Q_0(s,j\wedge k)}.
\end{align}
For the other term on the RHS of \eqref{eqn:big-ave-a1} we have, 
\beq
\abs{\im\sum_{i:|i-j|>\ell}\frac{\bu_i-\bu_k}{x_i-z}-\eta\sum_{i:|i-j|>\ell}\frac{\bu_i-\bu_k}{|\gamma_i-\gamma_j|^2}}\leq \varphi  \eta\sum_{i:|i-j|>\ell}\frac{\eta^2+\frac{|\gamma_i-\gamma_j|}{N\rho_i}}{|\gamma_i-\gamma_j|^4}\pare{|\bu_i-\bu_j|+|\bu_j-\bu_k|}
\eeq
For the summation over $i>j/2$, using spatial regularity of $\bu$ in \eqref{eqn:bua-dif},
\begin{align} \label{eqn:big-ave-a2}
&\eta\sum_{i:|i-j|>\ell, i > j/2}\frac{\eta^2+\frac{|\gamma_i-\gamma_j|}{N\rho_i}}{|\gamma_i-\gamma_j|^4}\pare{|\bu_i-\bu_j|+|\bu_j-\bu_k|} \notag \\
&\leq \varphi \eta\sum_{i>j/2, |i-j|>\ell}\frac{\eta^2+\frac{|\gamma_i-\gamma_j|}{N\rho_{i}}}{|\gamma_i-\gamma_j|^4}\pare{\frac{|\gamma_i-\gamma_j|}{Ns(s+\rho_j)^2}+\frac{|\gamma_j-\gamma_k|}{Ns(s+\rho_{j\wedge k})^2}} \notag \\
&\leq \eta\sum_{i>j/2, |i-j|>\ell} \pare{\frac{N\rho_j\eta^2}{|\gamma_i-\gamma_j|^3}+\frac{1}{|\gamma_i-\gamma_j|^2}}\frac{1}{Q_0(s,j)}+ \frac{\eta^2+\frac{|\gamma_i-\gamma_j|}{N\rho_{i}}}{|\gamma_i-\gamma_j|^4}\frac{R}{Q_0(s,j\wedge k)}\notag \\
&\lesssim \pare{\frac{(N\rho_j)(N\eta\rho_{j\vee \ell})^3}{\ell^2}+\frac{(N\rho_{j\vee l})^2\eta}{\ell}}\frac{1}{Q_0(s,j)}+\pare{\frac{(N\rho_{j\vee \ell})^4\eta^3}{\ell^3}+\frac{(N\rho_{j\vee \ell})^2\eta}{\ell^2}}\frac{R}{Q_0(s,j\wedge k)} \notag \\ 
& \leq \frac{ N^2 \eta \kappa_{j \vee \ell}}{Q_0 (t, j \wedge k)} + \frac{R N \rho_{j \vee \ell}}{Q_0 (t, j \wedge k)} 
\end{align}
In the fourth line we used \eqref{eqn:misc-8}. We simplified some errors using the assumption that $\ell \geq N \eta \sqrt{ \kappa_{k \vee R} + \eta}$.  For the summation over $i<j/2$ we have, since $| \gamma_i - \gamma_j | \asymp \kappa_j$, 
\begin{align}
&\eta\sum_{i:|i-j|>\ell, i < j/2}\frac{\eta^2+\frac{|\gamma_i-\gamma_j|}{N\rho_i}}{|\gamma_i-\gamma_j|^4}\pare{|\bu_i-\bu_j|+|\bu_j-\bu_k|}\notag \\
\leq & \frac{\eta}{ \kappa_j^4} \sum_{i:|i-j|>\ell, i < j/2} \left( \eta^2 + \frac{ \kappa_j}{N \rho_i } \right) \left( \frac{ \kappa_j}{N s (s + \rho_i)^2} + \frac{R}{Q_0 (s, j \wedge k ) } \right) \notag\\
\leq & \frac{1}{ Q_0 (s, j) } \left( \frac{ \eta^3 N^2}{ \kappa_j} + \frac{N \eta}{ \rho_j} \right)+ \frac{\eta}{ \kappa_j^4} \frac{ R }{ Q_0 (s, j \wedge k ) } \left( \eta^2 j + \kappa_j N^{-2/3} j^{2/3} \right) 
\end{align}
which is seen to be smaller than the error in the last line of \eqref{eqn:big-ave-a1}, using that $j \geq \ell$. We used \eqref{eqn:misc-9}. \qed


\bel
For the term \eqref{eqn: ave mod dyn, part 2-4} we have
\begin{align} \label{eqn: ave mod dyn, part 2-4 bound}
\left| \frac{1}{N}\im\sum_{i:|j-i|\leq \ell}\frac{(a_i-a_j)(\fu_i(s)-\bu_k(s))}{x_i-z} \right| \leq \varphi^{C_1} \left(  \frac{N \eta ( \kappa_{j \vee \ell} + \eta)}{R Q_a (s, j) } + \frac{ \ell^2}{N \eta Q_0 (s, j \wedge k) } \right)
\end{align}
\eel
\proof Using that $|a_j-a_i|\leq \frac{|i-j|}{R}$ (see \eqref{def: a_j}), Property $(\pP_a)$ and \eqref{eqn:bua-dif}, we have that the term \eqref{eqn: ave mod dyn, part 2-4}, is bounded by
\begin{align}
& \frac{1}{N}\sum_{i:|i-j|<\ell} \frac{\eta}{(\gamma_i-\gamma_j)^2+\eta^2}\frac{|i-j|}{R}\pare{|u_i(s)-\bar{u}_i(s)|+|\bar{u}_i(s)-\bar{u}_k(s)|} \notag \\
\lesssim & \frac{1}{N}\sum_{i:|i-j|<\ell}\frac{\eta}{(\gamma_i-\gamma_j)^2+\eta^2}\frac{|i-j|}{R}\pare{\frac{1}{Q_a(s,i)} + \frac{R}{N\rho_{i\wedge k}}\frac{1}{Ns(s+\rho_{i\wedge k})^2}} \notag \\
\lesssim & \frac{1}{N}\sum_{i:|i-j|<\ell}\left\{ \frac{\eta}{(\gamma_i-\gamma_j)^2+\eta^2}\frac{|i-j|}{R}\frac{1}{Q_a(s,i)} + \frac{|i-j|}{\eta Q_0(s,i\wedge k)} \right\}.
\end{align}
We divide the sum into three parts. First consider the sum over $j/2<i<2j$. In this case $i\asymp j$ and $|\gamma_i-\gamma_j|\asymp \frac{|i-j|}{N\rho_j}$. So, 
\begin{align}
&\frac{1}{N}\sum_{i:j/2\vee (j-\ell)<i<2j\wedge (j+\ell)}\left\{ \frac{\eta}{(\gamma_i-\gamma_j)^2+\eta^2}\frac{|i-j|}{R}\frac{1}{Q_a(s,i)} + \frac{|i-j|}{\eta Q_0(s,i\wedge k)} \right\} \notag \\
&\lesssim \frac{1}{N}\sum_{i:j/2\vee (j-\ell)<i<2j\wedge (j+\ell)}\left\{ \frac{\eta}{(\frac{|i-j|}{N\rho_j})^2+\eta^2}\frac{|i-j|}{R}\frac{1}{Q_a(s,j)} + \frac{|i-j|}{\eta Q_0(s,j\wedge k)} \right\} \notag \\
&\lesssim \pare{\frac{\rho_j}{R}\frac{1}{Q_a(s,j)}\sum_{i:j/2\vee (j-\ell)<i<2j\wedge (j+\ell)}\frac{1}{(\frac{|i-j|}{N\eta\rho_j})^2+1}\frac{|i-j|}{N\eta\rho_j}} + \frac{(j\wedge \ell)^2}{N\eta Q_0(s,j\wedge k)} \notag \\
&\lesssim \frac{N\eta\kappa_j}{R}\frac{1}{Q_a(s,j)}\int_0^{\frac{j\wedge \ell}{N\eta\rho_j}}\frac{x}{x^2+1}dx+\frac{(j\wedge \ell)^2}{N\eta Q_0(s,j\wedge k)} \notag \\
&\lesssim_\phi \pare{1 \wedge \frac{j\wedge \ell}{N\eta\rho_j}}^2\frac{N\eta\kappa_j}{R}\frac{1}{Q_a(s,j)} +\frac{(j\wedge \ell)^2}{N\eta Q_0(s,j\wedge k)}.
\end{align}

Second, consider the sum over $i<\frac{j}{2}$, which only exists if $j<\ell$. In this case, $|i-j|\asymp j$ and $|\gamma_i-\gamma_j|\asymp \kappa_j$. So,
\begin{align}
&\frac{1}{N}\sum_{i:(0\vee j-\ell)<i<j/2} \left\{ \frac{\eta}{(\gamma_i-\gamma_j)^2+\eta^2}\frac{|i-j|}{R}\frac{1}{Q_a(s,i)} + \frac{|i-j|}{\eta Q_0(s,i\wedge k)} \right\} \notag\\
&\lesssim \mathbbm{1}_{j<2\ell}\pare{\frac{1}{N}\frac{\eta}{\kappa_j^2+\eta^2}\frac{j}{R}\frac{j}{Q_a(s,j)} + \frac{j^2}{N\eta Q_0(s,j\wedge k)}}\notag \\
&= \mathbbm{1}_{j< \ell}\pare{\frac{N\eta\kappa_j^3}{\kappa_j^2+\eta^2}\frac{1}{R}\frac{1}{Q_a(s,j)} + \frac{j^2}{N\eta Q_0(s,j\wedge k)}}
\end{align}
In the last line, we use $j \asymp N\rho_j^3$.  We used \eqref{eqn:misc-9} and \eqref{eqn:misc-5a} above. 

Third, consider the sum over $i>2j$, which only exists if $j<\ell$. In this case, $|i-j|\asymp i$ and $|\gamma_i-\gamma_j|\asymp \kappa_i$. So,
\begin{align}
&\frac{1}{N}\sum_{i:2j<i<j+\ell} \left\{ \frac{\eta}{(\gamma_i-\gamma_j)^2+\eta^2}\frac{|i-j|}{R}\frac{1}{Q_a(s,i)} + \frac{|i-j|}{\eta Q_0(s,i\wedge k)} \right\} \notag \\
&\asymp \pare{\frac{1}{N}\sum_{i:2j<i<j+\ell}\frac{\eta}{(\frac{i}{N})^{\frac{4}{3}}+\eta^2}\frac{i}{R}\frac{1}{Q_a(s,i)}} + \mathbbm{1}_{j<l}\frac{l^2}{N\eta Q_0(s,j\wedge k)} \notag \\
&\lesssim \frac{N\eta}{R}\int_{\frac{2j}{N}}^{\frac{j+\ell}{N}}\frac{x}{x^{\frac{4}{3}}+\eta^2}dx\frac{1}{Q_a(s,j)}+ \mathbbm{1}_{j<\ell}\frac{\ell^2}{N\eta Q_0(s,j\wedge k)} \notag \\
&=\frac{N\eta^2}{R}\int_{\frac{2j}{N\eta^{3/2}}}^{\frac{j+\ell}{N\eta^{3/2}}}\frac{y}{y^{\frac{4}{3}}+1}dy\frac{1}{Q_a(s,j)}+ \mathbbm{1}_{j<l}\frac{\ell^2}{N\eta Q_0(s,j\wedge k)} \notag \\
&\lesssim \mathbbm{1}_{j<\ell}\pare{\pare{\mathbbm{1}_{\kappa_\ell>\eta}\frac{N\eta^2}{R}  \frac{ \kappa_\ell}{\eta} +\mathbbm{1}_{\kappa_l<\eta}\frac{l^2}{N\eta R}}\frac{1}{Q_a(s,j)}+\frac{\ell^2}{N\eta Q_0(s,j\wedge k)}}.
\end{align}
Using the assumption that $\ell \geq \varphi N \eta \sqrt{ \kappa_{k \vee R} + \eta}$ we find a total error of,
\beq
\frac{N \eta ( \kappa_{j \vee \ell} + \eta)}{R Q_a (s, j) } + \frac{ \ell^2}{N \eta Q_0 (s, j \wedge k) }
\eeq
which yields the claim. \qed


\subsubsection{Proof of Lemma \ref{lem: ave mod dyn, part 2}}

By the decomposition in \eqref{eqn: ave mod dyn, part 2-1} - \eqref{eqn: ave mod dyn, part 2-4}, and the estimates \eqref{eqn: ave mod dyn, part 2-1 bound}, \eqref{eqn: ave mod dyn, part 2-2 bound}, \eqref{eqn: ave mod dyn, part 2-3 bound} and \eqref{eqn: ave mod dyn, part 2-4 bound} we see that, up to error terms from these estimates, the quantity in question in \eqref{eqn:big-flow-b} is,
\begin{align}
a_j\pare{\im\pare{\frac{\e^{s/2}}{N}f_0(z_s)-\msc(z)\bu_j(s)} - \eta\partial_s\bu_j(s)} = \varphi^{C_1} \O\pare{\sqrt{\kappa_j+\eta}\frac{N\eta\rho_j}{Q_0(t,j)}}
\end{align}
where we used \eqref{eqn:f0-est}.  We conclude the proof by adding up the various errors in these five estimates. We simplified error terms using that $\eta \leq \kappa_\ell$ (and in the bulk that $k \geq R$). We also used $\ell^2 \geq (N \eta)^2 \kappa_{j \vee \ell}$ to simplify the $Q_0$ term.  \qed

\subsection{Proof of Lemma \ref{lem: ave interpolated dyn}} \label{sec:ave-proof} 
Applying Duhamel formula to \eqref{def: intermediate m} and using Proposition \ref{prop:deltbu}, we have
\beq
m_i(s) = \pare{\US(u,s)\Av(u) \fu (u)}_i + (\bu_k(s)-\bu_k(u)).
\eeq
Therefore, the sum in Lemma \ref{lem: ave interpolated dyn} can be decomposed into two parts:
\begin{align}
&\frac{1}{N\sqrt{\kappa_j+\eta}} \im \sum_{i: |j-i|<\ell}\frac{m_i(s)-\bu_k(s)}{x_i(s)-z}\nonumber\\
\label{eqn: ave mod dyn, part 1}=&\frac{1}{N\sqrt{\kappa_j+\eta}} \im \sum_{i:|j-i|<\ell}\frac{\pare{\US (u,s)\Av(u)\fu (u)}_i-\pare{\Av(s)\UB (u,s) \fu (u)}_i + (\bu_k(s)-\bu_k(u))}{x_i-z}\\
\label{eqn: ave mod dyn, part 2}+&\frac{1}{N\sqrt{\kappa_j+\eta}}  \im \sum_{i:|j-i|<\ell}\frac{\pare{\Av(s) \fu (s)}_i-\bu_k(s)}{x_i-z}.
\end{align}
Note that we used $u(s) = \UB(u,s)u(u)$ above. The terms \eqref{eqn: ave mod dyn, part 1} and \eqref{eqn: ave mod dyn, part 2} were bounded separately in Lemma \ref{lem: ave mod dyn, part 1} and Lemma \ref{lem: ave mod dyn, part 2} above. 

\section{Quantitative edge universality}

\subsection{Bounded density for Gaussian ensembles}

We recall here that the Gaussian $\beta$-ensemble is the measure on $\{ \lambda_1 > \lambda_2 > \dots > \lambda_N \} =: \Delta_N \subseteq \rr^N$ with density proportional to,
\beq
p_\beta ( \blam ) := \exp\left[ \beta N \left( \frac{1}{N} \sum_{i < j} \log ( \lambda_i - \lambda_j) - \sum_{i=1}^N \frac{ \lambda_i^2}{4} \right)\right] .
\eeq
We denote
\beq
\blam := ( \lambda_1, \dots, \lambda_N) , \qquad \hblam := ( \lambda_2, \dots, \lambda_N)
\eeq

We require the following.
\bep \label{prop:density} Let $\eps >0$. For all $N$ sufficiently large, the following holds. 
Let $I = [s, s+a]$ be an interval with $s \geq 0$, $a >0$. Then,
\beq
\pp\left[ N^{2/3} ( \lambda_1 - \lambda_2 ) \in I \right] \leq N^{\eps} a + \e^{ -  \varphi}.
\eeq
\eep

For any $\delta >0$, let  $\A_\delta \subseteq \Delta_N$ be set
\beq
\A_\delta := \left\{\blam \in \Delta_N : | \lambda_i - \gamma_i | \leq \frac{N^{\delta}}{N^{2/3} \hat{i} } , \forall i \in [\![1, N]\!] \right\}
\eeq
\bel \label{lem:density}
For any $\delta >0$ we have for all $N$ sufficiently large for all $\blam \in \A_\delta$ that,
\beq
p_\beta ( \blam) \leq 2 p_\beta ( \lambda_1 + b , \hblam ).
\eeq
for all $b \in [0, N^{-2/3-10\delta}]$
\eel
\proof Let us write
\beq \label{eqn:density-a1}
p_\beta = g ( \hblam ) \exp\left[ N \beta f ( \blam ) \right]
\eeq
with
\beq
f ( \blam ) := \frac{1}{N} \sum_{j>1} \log ( \lambda_1 - \lambda_j) - \frac{\lambda_1^2}{4}.
\eeq
We estimate, with $K = N^{5 \delta}$,
\begin{align}
& f ( \blam ) \leq \frac{1}{N} \sum_{j=1}^K \log ( \lambda_1 + b - \lambda_j ) + \frac{1}{N} \sum_{j>K} \log ( \lambda_1 - \lambda_j) - \frac{\lambda_1^2}{4} \notag\\
= & f ( \lambda_1 + b , \hblam ) - \int_{\lambda_1}^{\lambda_1+b} \frac{\d}{\d u } \left(\frac{1}{N} \sum_{j>K} \log ( u - \lambda_j) - \frac{u^2}{4}  \right) \d u \notag\\
= &f ( \lambda_1 + b , \hblam ) - \int_{\lambda_1}^{\lambda_1+b}  \left(\frac{1}{N} \sum_{j>K} \frac{1}{u - \lambda_j }  - \frac{u}{2}  \right) \d u
\end{align}
We estimate the derivative on the RHS. First we have for $u \in [\lambda_1, \lambda_1 + b]$,
\beq
\frac{u}{2} = 1 + \O ( N^{-1/2} )
\eeq
On the other hand, 
\begin{align}
& \left| \frac{1}{N} \sum_{j>K} \frac{1}{ \lambda_j - u} - N \int_{\gamma_j}^{\gamma_{j+1} } \frac{\rhosc (x) \d x}{ x - 2} \right| \notag\\
\lesssim & \frac{1}{N} \frac{N^{\delta}}{N^{2/3}} \sum_{j > K } \frac{1}{ | \gamma_j - 2|^2} \leq \frac{N^{\delta}}{N^{1/3}} \sum_{j > K} \frac{1}{j^{4/3}} \leq \frac{1}{N^{1/3}} ,
\end{align}
and 
\beq
\int_{\gamma_K}^2 \frac{\rhosc(x)}{|x-2|} \d x \lesssim \kappa_K^{1/2} \leq \frac{N^{3 \delta}}{N^{1/3}}.
\eeq
Therefore,
\beq
\sum_{j>K} \frac{1}{u - \lambda_j } = - \msc (2) + \O ( N^{3\delta-1/3} ) = 2 + \O ( N^{3\delta-1/3} ). 
\eeq
Since $b \leq N^{-10\delta-2/3}$ we see that,
\beq
f ( \blam ) \leq f ( \lambda_1 + b , \blam) + \O ( N^{-\delta-1} ).
\eeq
This, together with \eqref{eqn:density-a1} completes the proof. \qed

\vspace{5 pt} 

\noindent{\bf Proof of Proposition \ref{prop:density}.} By \cite[Theorem 2.4]{bourgade2014edge} we have for any $\delta >0$ that
\beq
\pp\left[ \blam \notin \A_\delta \right] \leq \e^{ - \varphi}.
\eeq
Let,
\beq
 \hA_\delta := \left\{ \hblam \in \Delta_{N-1} : | \lambda_i - \gamma_i| \leq \frac{N^{\delta}}{N^{2/3} \hat{i}^{1/3}} \forall i \in [\![ 2, N]\!] \right\}
 \eeq
Let $I = N^{2/3}[S, S+A]$ and write,
\begin{align}
&\pp\left[ N^{2/3} ( \lambda_2 - \lambda_1) \in I \right] \leq \pp\left[ N^{2/3} ( \lambda_2 - \lambda_1) \in I, \blam \in \A_\delta \right] + \e^{- \varphi} \notag\\
= & \frac{1}{Z_\beta} \int_{ \hblam \in \hA_\delta} \d \hblam \int_{\lambda_2+S}^{\lambda_2 +S +A } \d \lambda_1 p_\beta (\blam) \1_{\{  | \lambda_1 - \gamma_1 \leq N^{\delta-2/3} \} } + \e^{- \varphi}
\end{align}
Let $B = N^{-2/3-10\delta}$. 
By Lemma \ref{lem:density} for any $b \leq B$ we have, for any $\blam \in \A_\delta$,
\beq
p_\beta ( \blam) \leq 2 p_\beta ( \lambda_1 + b , \hblam).
\eeq
Let us assume that $A \leq N^{-2/3-11 \delta}$. Applying the above inequality we have for any $\hblam \in \hA_\delta$ and
\beq
\lambda_1 \in [\lambda_2 +S , \lambda_2 +S + A] \cap [ \gamma_1 - N^{-2/3+\delta}, \gamma_1 + N^{-2/3+\delta} ]
\eeq
that 
\beq
p_\beta ( \blam) \leq 2 p_\beta ( \lambda_2 + S + A + b , \hblam )
\eeq
for any $b \leq B/2$. Integrating over $b \in [0 ,B/2]$ we find,
\beq
p_\beta ( \blam) \leq \frac{4}{B} \int_{\lambda_2 + S + A}^{ \lambda_2 +S + A + B/2} p_\beta (u , \hblam ) \d u
\eeq
Therefore,
\begin{align}
&\frac{1}{Z_\beta} \int_{ \hblam \in \hA_\delta} \d \hblam \int_{\lambda_2+S}^{\lambda_2 +S +A } \d \lambda_1 p_\beta (\blam) \1_{\{  | \lambda_1 - \gamma_1 \leq N^{\delta-2/3} \} } \notag\\
\leq & \frac{1}{Z_\beta} \int_{ \hblam \in \hA_\delta} \d \hblam \int_{\lambda_2+S}^{\lambda_2+S+A} \d \lambda_1 \frac{4}{B} \int_{\lambda_2 + S + A}^{ \lambda_2 +S + A + B/2} p_\beta (u , \hblam ) \d u\notag\\
\leq & \frac{4 A}{B} \frac{1}{Z_\beta} \int_{ \hblam \in \hA_\delta} \d \hblam \int_{\lambda_2 + S + A}^{ \lambda_2 +S + A + B/2}  p_\beta (u , \hblam ) \d u \notag\\
\leq & \frac{4A}{B} \frac{1}{Z_\beta} \int_{\Delta_N} \d \blam p_\beta ( \blam) = \frac{4A}{B}.
\end{align}
This yields the claim when $a \leq N^{-11 \delta}$. The general case follows from adjusting the $N^{\eps}$ prefactor in the error. \qed

\subsection{Proof of Theorem \ref{thm:gap-univ}}

Let $H_t$ be a  Wigner matrix of the form,
\beq \label{eqn:gde}
H_t = \e^{-t/2} W + \sqrt{1-\e^{-t}} G
\eeq
where $t = N^{-\omega}$, $G$ is a GOE matrix and $W$ is another generalized Wigner matrix. If $\{ \lambda_i \}_{i=1}^N$ are the eigenvalues $H$ then Theorem \ref{thm:main-homog} together with \eqref{eqn:bua-dif} guarantee the existence of a coupling to the eigenvalues $\{ \mu_i \}_{i=1}^N$ of a GOE matrix so that
\beq
\left| N^{2/3} ( \lambda_1 - \lambda_2) - N^{2/3} ( \mu_1 - \mu_2) \right| \leq \frac{N^{5 \omega}}{N}
\eeq
with overwhelming probability. Combined with Proposition \ref{prop:density} this yields,
\beq
\sup_{r \in \rr} \left| \pp\left[ N^{2/3} ( \lambda_1 - \lambda_2) \geq r \right] - \pp\left[ N^{2/3} ( \mu_1 - \mu_2) \geq r \right] \right| \leq \frac{N^{6 \omega}}{N}
\eeq
for all $N$ sufficiently large, for matrices of the form \eqref{eqn:gde}. 

Now, given a smooth Wigner matrix $X$ as in Definition \ref{def:smooth}, [??] of [??] guarantees us a matrix $H_t$ of the form in \eqref{eqn:gde} (with $t = N^{-\omega}$) s.t.
\beq
\mathrm{TV} (H_t, X) \leq C N^{-10}
\eeq
where $\mathrm{TV}$ denotes the total variation distance between, in this case, two probability measures on the space of symmetric $N \times N$ matrices. Therefore,
\beq
\sup_{ r \in \rr} \left| \pp\left[ N^{2/3} ( \lambda_2 (X) - \lambda_1 (X) ) \geq r \right] - \pp\left[ N^{2/3} ( \lambda_2 (H_t) - \lambda_1 (H_t) ) \geq r \right] \right| \leq C N^{-2}.
\eeq
This completes the proof. 
 \qed
 
\subsection{Lower bound}

We will prove the lower bound, Theorem \ref{thm:main-lb} only in the case $\beta=1$, real symmetric Wigner matrices. The case $\beta=2$ is similar. 

Let the coupled DBM in \eqref{eqn:def-dbm} be the eigenvalues of a Wigner matrix $X$ and a GOE matrix $Y$. Define, for $W = X, Y$,
\begin{align}
\xi^W (z, t) := \frac{1}{\Im[ \msc (z) ] }\Im \left[ \sum_{i=1}^N \log ( \lambda_i (W) - z_t ) - N \int \log (x -z_t) \rhosc (x) \d x \right]
\end{align}
\bep
Let $z = \gamma_1 + \i N^{-2/3+\eps_1}$ and $c>0$. For all $t\in (c, 1)$ we have
\beq \label{eqn:LSS-prop-1}
\bu_1 (t) = \frac{1}{ \Im[ \msc(z) ]} \frac{\e^{t/2}}{N} \Im\left[ f_0 (z_t) \right] + \O ( N^{2\eps_1-4/3} )
\eeq
and so
\beq \label{eqn:LSS-prop-2}
\lambda_1 (t) - \mu_1(t) = \frac{ \xi^X (z, t)- \xi^Y(z, t)}{N} + \O ( N^{-4/3+3 \eps_1} ) ,
\eeq

\eep
\proof The first estimate \eqref{eqn:LSS-prop-1} follows from \eqref{eqn:f0-est} and \eqref{eqn:buk-time-bd}. The second follows from and Theorem \ref{thm:main-homog}, \eqref{eqn:LSS-prop-1} and the fact that  
\beq
\del_\nu \sum_{i=1}^N \log ( \xnu_i(0) - z_t)=  f_0 (z_t).
\eeq
 \qed
 
\bel
Let $f : \rr \to \rr$ be the function,
\beq
f(x) := \frac{1}{ \Im[\msc(z)]} \Im\left[ \log (x - z_t ) \right].
\eeq
Then $f$ is smooth for all $|x| < 2 + \delta$ for some $\delta >0$.
\eel
\proof By \eqref{eqn:char-bourgade} we have that $\Re[z_t] > 2 + 2 \delta$ for some $\delta >0$. The claim then follows from the fact that,
\beq
f(x) = \frac{1}{\Im[ \msc (z) ] } \arctan\left( \frac{\Im[z_t]}{\Re[z_t]-x} \right)
\eeq
 and that $\Im[z_t] \asymp \Im[ \msc (z)]$. \qed

 \bel \label{lem:dif-LSS}
We have
 \beq
 \ee[\xi^X]- \ee[\xi^Y] = \e^{-2t} s_4 + \O ( N^{-1/3+5 \eps_1} ).
 \eeq
 \eel 
 \proof From \cite[(1.17)]{landon2024some} we have that for any real symmetric Wigner matrix $W$,
 \begin{align}
 \e^{t/2} \ee[ \xi^W] &= - \frac{1}{2 \pi} \int_{-2}^2 \frac{f(x)}{\sqrt{4-x^2}} \d x + \frac{f(2)+f(-2)}{4} \notag\\
 + & \frac{s_4}{2\pi} \int_{-2}^2 f(x) \frac{x^4-4x^2+2}{\sqrt{4-x^2}} \d x + \O ( N^{-1/4} ).
 \end{align}
 By \cite[(A.16)]{baik2016fluctuations} we have
 \begin{align}
 & \frac{1}{2 \pi} \int_{-2}^2 \log (z-x) \frac{x^4-4x^2+2}{\sqrt{4-x^2}} \d x \notag\\
 = &\frac{1}{8}(z^3-2z)\sqrt{z^2-4}  - \frac{z^4}{8} + \frac{z^2}{2} - \frac{1}{4} \notag\\
 = & \frac{1}{8} (z^3-2z)(\sqrt{z^2-4}-z) +\frac{z^2}{4} - \frac{1}{4} \notag\\
 =& \frac{1}{4}(z^3-2z)\msc(z) + \frac{z^2}{4} - \frac{1}{4} 
 \end{align}
 Then,
 \begin{align}
 &\frac{1}{ \Im[ \msc (z)]} \Im\left[  \frac{1}{4}(z_t^3-2z_t)\msc(z_t) + \frac{z_t^2}{4} - \frac{1}{4} \right] \notag\\
 = & \frac{1}{4} \Re[ z_t^3-2z_t] \frac{\Im[ \msc(z_t)]}{\Im[ \msc(z)]} +\frac{1}{4}\Re[ \msc(z_t)] \frac{\Im[ z_t^3 -2 z_t]}{\Im [\msc(z)]}+ \frac{1}{2} \Re[z_t] \frac{\Im[z_t]}{\Im[ \msc(z)]} \notag\\
 = & \frac{\e^{-t/2}}{4} \left(8 \cosh^3(t/2)-4 \cosh(t/2) -24 \sinh(t/2) \cosh^2(t/2) + 4 \sinh(t/2) \right) \notag\\
 + & 2 \sinh(t/2) \cosh(t/2) + \O ( N^{-1/3+5\eps_1} ) \notag\\
 = & \e^{-2t} + \O ( N^{-1/3+5\eps_1} ).
 \end{align}
 In the third line we used $\Re[z_t] = \cosh(t/2) + \O (N^{-1/3+ \eps_1 } )$, $\Im[z_t] = 2 \Im[ \msc (z) ] \sinh(t/2) + \O ( N^{-2/3+\eps_1} )$ and $\msc(z_t) = \e^{-t/2} \msc (z)$ and the final equality is just simplifying the expression on the previous lines.  The claim now follows. \qed
 
 \bel \label{lem:lb-1}
 Let $c_1>0$ and let $W$ be a Wigner matrix with eigenvalues $\{ \lambda_i \}_{i=1}^N$ so that
 \beq
 \sup_{ r \in \rr} \left| \pp\left[ N^{2/3} ( \lambda_1 - 2 ) \geq r \right] -\pp\left[ N^{2/3} ( \mu_1 - 2 ) + b N^{-1/3} \geq r \right] \right| \leq \frac{1}{N^{1/3+c_1}}
 \eeq
 for some $b \in \rr$. 
 Then,
 \beq
\left|  \ee[N^{2/3} ( \lambda_1 - 2 ) ] - \ee[ N^{2/3} ( \mu_1 - 2 ) + b N^{-1/3} ] \right| \leq N^{-1/3+\eps - c_1}
 \eeq
 for any $\eps >0$. 
 \eel
 \proof Let $X = N^{2/3} ( \lambda_1 -2) $ and $Y = N^{2/3} ( \mu_1 - 2)+b N^{-1/3}$. Using the rigidity estimates and an easy inequality, $\ee[ X^2] + \ee[ Y^2] \leq N^3$ we have,
 \beq
 \ee[X] = \int_0^\infty \left( \pp[ X > r ] - \pp[ X < -r ] \right) \d r = \int_0^{N^{\eps}} \left( \pp[ X > r ] - \pp[ X < -r ] \right) \d r + \O (N^{-2} ). 
 \eeq
 The claim then follows from using the same expression for $\ee[Y]$ and taking the difference. \qed
 
 \subsubsection{Proof of Theorem \ref{thm:main-lb}}
 
 From \eqref{eqn:LSS-prop-2}, and Lemma \ref{lem:dif-LSS} we have
 \beq
 \ee[ N^{2/3} ( \lambda_1 -2 ) ] - \ee[ N^{2/3} ( \mu_1 -2 ) ] = N^{-1/3} \e^{-2t} s_4 + \O ( N^{-2/3+\eps} )
 \eeq
 for any $\eps >0$. The claims then follow from applying this to Lemma \ref{lem:lb-1}. \qed

\appendix

\section{Proof of Lemma \ref{lem:im-m-ell}} \label{sec:im-m-ell-proof} 

For $z$ as in the statement of the lemma we have $\Im[ \msc (z) ] \asymp \sqrt{ \kappa_j + \eta}$ and that, with overwhelming probability,
\beq
|s(z) - \msc (z) | \leq \frac{\varphi^{1/100}}{N \eta} \leq \varphi^{-1/4} \sqrt{\kappa_j+\eta}
\eeq
and so it suffices to prove that
\beq
\frac{1}{N} \sum_{ i : |i-j| \geq \ell} \frac{\eta}{(x_i-\gamma_j)^2 + \eta^2} \leq \varphi^{-1/4} \sqrt{\kappa_j + \eta}
\eeq
with overwhelming probability. Let us first consider the terms where $i < j/2$ (which arise only if $j \geq \ell$). Then 
\beq
\frac{1}{N} \sum_{i<j/2} \frac{\eta}{(x_i - \gamma_j)^2 + \eta^2} \lesssim \frac{ \kappa_j^{3/2} \eta}{ \kappa_j^2 + \eta^2} \lesssim \sqrt{\kappa_j+\eta} \frac{ \eta \kappa_j^{1/2}}{ \kappa_\ell^{3/2}} \leq \varphi^{-1} \sqrt{\kappa_j+\eta}. 
\eeq
using $j \geq \ell$ in the last inequality. For the terms with $i > 2j$ we have,
\begin{align}
\frac{1}{N} \sum_{i > 2j, i > j+\ell} \frac{\eta}{ (x_i - \gamma_j)^2 + \eta^2} \lesssim N^{1/3} \sum_{i > j + \ell} \frac{\eta}{ i^{4/3} + (N^{2/3} \eta)^2}
\end{align}
We can bound this above by
\beq
 N^{1/3} \sum_{i > j + \ell} \frac{\eta}{ i^{4/3} + (N^{2/3} \eta)^2} \leq  \frac{ N^{1/3} \eta}{ (j + \ell)^{1/3} } \leq  \sqrt{\kappa_j + \eta} \frac{ \eta^{1/2}}{\kappa_\ell^{1/2}} \lesssim \varphi^{-1/4} \sqrt{\kappa_j+\eta}
\eeq
using that $\kappa_\ell \gtrsim \varphi^{2/3} \eta$ by assumption. We are left with the sum where $j/2 < i < 2j$ this is bounded by
\begin{align}
N \kappa_j \sum_{ j/2 < i < j, |i-j| > \ell} \frac{\eta}{(i-j)^2 + (N \rho_j \eta)^2} \lesssim \frac{\kappa_j \eta}{ \kappa_\ell^{3/2}} \leq \sqrt{\kappa_j+\eta} \frac{ \eta \kappa_j^{1/2}}{\kappa_\ell^{3/2}} \leq \varphi^{-1} \sqrt{\kappa_j + \eta}
\end{align}
by assumption. This completes the proof. \qed

\section{Martingale result}

We use the following large deviation result for continuous martingales from \cite[Appendix B.6, equation (18)]{shorack2009empirical}:

\bel\label{lem:martingale}
Let $M$ be a continuous martingale such that $\langle M \rangle_t \leq \mu$ almost surely. Then,
\beq
\pp \left[ \sup_{0\leq u \leq t}|M_u|\geq \lambda \right] \leq 2e^{-\frac{\lambda^2}{2\mu}}.
\eeq
\eel

\section{Various estimates of sums and integrals}

\bel
For $ (N \rho_j)^{-1} \leq a \leq 1$ and $b >1 $ we have,
\beq \label{eqn:misc-1}
\frac{1}{N} \sum_{i=1}^N \frac{1}{N \rho_i } \frac{1}{ ( | \gamma_i - \gamma_j| + a )^b} \leq \frac{C_b}{Na^{b-1}} .
\eeq
If $b=1$ then the same estimate holds but with $C_b$ replaced by $C \log N$. 
\eel
\proof We just do the case $b>1$, the case $b=1$ is similar. Let us first consider the case that $\kappa_j > a$. We first consider the part of the sum where $i < j/2$. We have,
\begin{align}
\frac{1}{N} \sum_{i < j/2} \frac{1}{N \rho_i} \frac{1}{ (| \gamma_i - \gamma_j|  + a)^b} \lesssim \frac{1}{N^{5/3}  ( \kappa_j +a)^{b}} \sum_{i < j/2} \frac{1}{i^{1/3}} \lesssim \frac{ \kappa_j}{ N ( \kappa_j + a)^b} \leq \frac{ C}{ N a^{b-1}}.
\end{align}
Now for $i > 2j$ we have,
\begin{align}
\frac{1}{N} \sum_{ i > 2 j } \frac{1}{N \rho_i} \frac{1}{ ( | \gamma_i - \gamma_j | + a )^b} \lesssim \frac{1}{N^{5/3}} \sum_{ i > 2 j } \frac{1}{ i^{1/3}} \frac{1}{ ( (i/N)^{2/3} + a )^b}
\end{align}
If $\kappa_j > a $ then we estimate this by
\beq
\frac{1}{N^{5/3}} \sum_{ i > 2 j } \frac{1}{ i^{1/3}} \frac{1}{ ( (i/N)^{2/3} + a )^b} \lesssim \frac{1}{N^{5/3}} \sum_{i > 2 j } \frac{1}{ i^{1/3} (i/N)^{2 b/3} } \lesssim \frac{1}{N \kappa_j^{b-1}} \lesssim \frac{1}{N a^{b-1}}.
\eeq
Finally, for the terms $j/2 \leq i \leq 2j$ we have,
\begin{align}
& \frac{1}{N} \sum_{j/2 < i < 2j} \frac{1}{N \rho_i } \frac{1}{ ( | \gamma_i - \gamma_j| + a )^b}  \lesssim \frac{1}{N^{5/3} j^{1/3}} \sum_{j/2 < i < 2j} \frac{1}{ (a + \frac{ |i-j|}{j^{1/3} N^{2/3}} )^b} \notag\\
\lesssim & \frac{1}{N^{5/3} j^{1/3}} \frac{j^{1/3} N^{2/3}}{a^{b-1}} \lesssim \frac{1}{N a^{b-1}} .
\end{align}
This completes the case that $\kappa_j > a$. If $a \geq \kappa_j$, we first consider the sum where $i \leq C_1 N a^{3/2}$ for some $C_1 >0$. This contributes,
\beq
\frac{1}{N} \sum_{i < C_1 N a^{3/2}} \frac{1}{N \rho_i } \frac{1}{ ( | \gamma_i - \gamma_j| + a )^b} \lesssim \frac{1}{N^{5/3} a^b} \sum_{i < C_1 N a^{3/2}} \frac{1}{i^{1/3}} \lesssim \frac{1}{ N a^{b-1}}. 
\eeq 
By taking $C_1$ sufficiently large we guarantee that $| \gamma_i - \gamma_j| \geq c i^{2/3} / N^{2/3}$ for $i \geq C_1 N a^{3/2}$. Then,
\beq
\frac{1}{N} \sum_{i \geq C_1 N a^{3/2}} \frac{1}{N \rho_i } \frac{1}{ ( | \gamma_i - \gamma_j| + a )^b}  \lesssim \frac{1}{N^{5/3}} \sum_{i \geq C_1 N a^{3/2}} \frac{1}{i^{1/3} (i/N)^{2b/3}} \lesssim \frac{1}{N a^{b-1}}
\eeq
This completes the case $a > \kappa_j$ and the proof. \qed

\bel
For $(N \rho_j )^{-1} \leq a \leq 1$ and $b >1$ we have,
\beq \label{eqn:misc-2}
\frac{1}{N} \sum_{i=1}^N \frac{1}{N \rho_i^2} \frac{1}{ (| \gamma_i - \gamma_j| + a )^b} \leq \frac{ C_b}{N a^{b-1} \rho_j}
\eeq
\eel
\proof The part of the sum with $i \geq j/2$ can easily be bounded using \eqref{eqn:misc-1} using that $\rho_i \geq c \rho_j$ for such $i$. For the part of the sum with $i < j/2$ we bound,
\begin{align}
 & \frac{1}{N} \sum_{i<j/2} \frac{1}{N \rho_i^2} \frac{1}{ (| \gamma_i - \gamma_j-2| + a )^b}  \lesssim \frac{1}{ ( | \gamma_j| + a)^b} \frac{1}{N^{4/3}} \sum_{i < j/2} \frac{1}{i^{2/3}}  \notag\\
\lesssim & \frac{1}{N} \frac{ \rho_j}{ a^{b-1} \kappa_j} \lesssim \frac{1}{N a^{b-1} \rho_j } ,
\end{align}
which yields the claim. \qed

\bel 
Let $A \geq \varphi^{1/2}$. On the event that \eqref{eqn:rig} holds we have,
\beq \label{eqn:misc-3}
\frac{1}{N}\left| \sum_{j : |i-j| > A } \frac{1}{ x_i -x _j} - \frac{1}{ \gamma_i - \gamma_j} \right| \leq  \frac{\varphi}{N^{1/3}} \left( \frac{ i^{1/3}}{A} + \frac{1}{ i^{1/3} A^{1/3}} \right)
\eeq
\eel
\proof We have,
\begin{align}
\left| \sum_{j : |i-j| > A } \frac{1}{ x_i -x _j} - \frac{1}{ \gamma_i - \gamma_j} \right| \leq \frac{\varphi^{1/2}}{N^{1/3}} \sum_{ j : |i-j| > A } \left( i^{-1/3} + j^{-1/3} \right) \frac{j^{2/3} + i^{2/3}}{(i-j)^2} 
\end{align}
For $j > i$ we can bound the sum by,
\begin{align}
 & \sum_{ j : j > i+ A } \left( i^{-1/3} + j^{-1/3} \right) \frac{j^{2/3} + i^{2/3}}{(i-j)^2}   \lesssim \sum_{j > i + A } \frac{ j^{2/3}}{i^{1/3} (i-j)^2} \notag\\
\lesssim &  \sum_{j > i + A } \frac{ (j-i)^{2/3} + i^{2/3}}{i^{1/3}(i-j)^2} \lesssim  \left( \frac{ i^{1/3}}{A} + \frac{1}{ i^{1/3} A^{1/3}} \right).
\end{align}
For $j < i$ we split the sum into the cases that $j < i/2$ or $ j \geq i/2$. For the case that $j \geq i/2$ we have,
\begin{align}
 \sum_{ i/2 \leq j \leq i-A}   \left( i^{-1/3} + j^{-1/3} \right) \frac{j^{2/3} + i^{2/3}}{(i-j)^2} 
\lesssim \sum_{j \leq i - A} \frac{ i^{1/3}}{(i-j)^2} \lesssim \frac{ i^{1/3}}{A}.
\end{align}
For the case that $j \leq i/2$ we have, since $(i-j) \geq i/2$ in this case,
\begin{align}
\sum_{j \leq (i-A) \wedge i/2}  \left( i^{-1/3} + j^{-1/3} \right) \frac{j^{2/3} + i^{2/3}}{(i-j)^2} \lesssim \sum_{j \leq i/2} \frac{1}{j^{1/3} A^{1/3} i} \lesssim  \frac{1}{ i^{1/3} A^{1/3}}.
\end{align}
This completes the proof. \qed

\bel
We have
\beq \label{eqn:gamma-rho}
| \gamma_a - \gamma_b| \asymp \frac{ |a-b|}{N} \frac{1}{ \rho_a + \rho_b}
\eeq
\eel
\proof This follows from a straightforward estimate using
\beq
\int_{\gamma_a}^{\gamma_b} \rhosc (x) \d x = \frac{a-b}{N}. 
\eeq

\bel
We have
\beq \label{eqn:misc-4}
\int_{-2}^{\gamma_{j+\ell}} \frac{ \rhosc (y)}{ (y- \gamma_j)^2} \d y \lesssim  \frac{ \rho_{j+\ell}^2}{ \rho_\ell^3}
\eeq
\eel
\proof When $y > \gamma_{2j + \ell}$ we have that $(y - \gamma_j) \asymp y$ and so this part of the integral contributes $\O ( \rho_{j+\ell}^{-1})$. For the remaining part of the integral,
\begin{align}
 & \int_{\gamma_{2j+\ell}}^{\gamma_{j+\ell}} \frac{ \rhosc (y)}{ (y- \gamma_j)^2} \d y \lesssim \rho_{j+\ell} \int_{\gamma_{2j+\ell}}^{\gamma_{j+\ell}} \frac{1}{ (y- \gamma_j)^2} \d y \notag\\
\lesssim & \frac{ \rho_{j+\ell}}{|\gamma_j - \gamma_{j+\ell} | } \lesssim \frac{ \rho_{j+\ell}^2}{\rho_l^3}
\end{align}
using \eqref{eqn:gamma-rho}. \qed

\bel
We have  for $j \geq \ell$,
\beq \label{eqn:misc-5}
\int_{\gamma_{j-\ell}}^2 \frac{ ( \rhosc(y))^a}{ (y-\gamma_j)^2 N^{2-a} (s  + \rhosc (y))^{2-a} s^{1-a}} \d y  \lesssim \frac{ N \kappa_j}{ \ell N(t + \rho_j )[Nt (t+\rho_j ) \rho_j]^{1-a} } 
\eeq
\eel
\proof The part of the integral where $y > \gamma_{j/2}$ contributes,
\begin{align}
\int_{ \gamma_{(j-\ell)} \vee \gamma_{j/2}}^2 \frac{ ( \rhosc(y))^a}{ (y-\gamma_j)^2 N^{2-a} (s  + \rhosc (y))^{2-a} s^{1-a}} \d y \lesssim \frac{1}{ \rho_j^4 N^{2-a} s^{1-a}} \int_{\gamma_{j/2}}^2 \frac{ \rho(y)^a}{(s+ \rhosc (y))^{2-a}} \d y
\end{align}
If $ s> \rho_j$, then we can bound
\beq
\frac{1}{ \rho_j^4 N^{2-a} s^{1-a}} \int_{\gamma_{j/2}}^2 \frac{ \rho(y)^a}{(s+ \rhosc (y))^{2-a}} \d y \lesssim \frac{1}{\rho_j^4 N^{2-a} s^{1-a} } \frac{ \rho_j^{2+a}}{(s+ \rho_j)^{2-a}}
\eeq
which is sufficient. If $s < \rho_j$, then,
\begin{align}
& \int_{\gamma_{j/2}}^2 \frac{ \rho(y)^a}{(s+ \rhosc (y))^{2-a}} \d y = \int_{s^2}^2 \frac{ \rho(y)^a}{(s+ \rhosc (y))^{2-a}} \d y+\int_{\gamma_{j/2}}^{s^2} \frac{ \rho(y)^a}{(s+ \rhosc (y))^{2-a}} \d y \notag\\
\leq & s^{2a} + \int_0^{\kappa_j} y^{a-1} \d y = s^{2a} + \rho_j^{2a} \leq \rho_j^{2a}
\end{align}
which is sufficient. We record here the resulting estimate we have proven,
\beq \label{eqn:misc-5a}
\int_{\gamma_{j/2}}^2 \frac{ \rho(y)^a}{(s+ \rhosc (y))^{2-a}} \d y \lesssim \frac{ \rho_j^{2+a}}{(s+\rho_j)^{2-a}} ,
\eeq
for future use. 

The remaining part of the integral (non-zero only if $j \geq 2 \ell$ which we now assume) gives,
\begin{align}
& \int_{ \gamma_{j-\ell}}^{\gamma_{j/2}} \frac{ ( \rhosc(y))^a}{ (y-\gamma_j)^2 N^{2-a} (s  + \rhosc (y))^{2-a} s^{1-a}} \d y \asymp \frac{ \rho_j^a}{N^{2-a} (s+\rho_j)^{2-a} s^{1-a}} \int_{\gamma_{j-\ell}}^{\gamma_{j/2}} \frac{1}{ (y- \gamma_j)^2} \notag\\
\lesssim & \frac{ \rho_j^a}{N^{2-a} (s+\rho_j)^{2-a} s^{1-a}} \frac{\rho_j}{\rho_\ell^3}
\end{align}
which completes the proof. 

\bel
\beq \label{eqn:misc-6}
\int_{\gamma_{j+l}}^{\gamma_j-N^{-D}}\frac{\rho(y)dy}{Ns(s+\rho_j)^2|y-\gamma_j|} + \int_{\gamma_j+N^{-D}}^{\gamma_{j-l}}\frac{\rho(y)dy}{Ns(s+\rho(y))^2|y-\gamma_j|}\lesssim \frac{ \rho_{j+\ell}}{N s (s+ \rho_j)^2}
\eeq
\eel
\proof For the first integral, the part of the integral where $y > \gamma_{2j}$ is bounded using
\beq
\int_{\gamma_{j+\ell}}^{\gamma_{2j}} \frac{ \rhosc(y)}{ |y-\gamma_j|} \leq \int_{\kappa_j}^{\kappa_{j+\ell}} y^{-1/2} \leq \rho_{j+\ell}
\eeq
The part of the integral where $y < \gamma_j$, just bound $\rhosc (y) \leq \rho_j$. 

For the second integral, the part with $\gamma_j > y > \gamma_{j/2}$ gives $\frac{ \rho_j}{N s (s+\rho_j)^2}$. For the rest, if $ s> \rho_j$, then the bound is easy. Otherwise, using $|y- \gamma_j| \geq \rho_j^2$ for $y< \gamma_{j/2}$ gives
\beq
\int_{\gamma_{j/2}}^{\gamma_{j-l}}\frac{\rho(y) \d y}{Ns(s+\rho(y))^2|y-\gamma_j|} \lesssim \frac{1}{N s \rho_j^2} \int_{\kappa_{j-\ell}}^{ \kappa_{j/2}} \frac{1}{ \sqrt{y}} \d y \lesssim \frac{\rho_j}{Ns (\rho_j + s)^2}
\eeq
which finishes the proof. \qed

\bel Assume $|j- k| \leq \ell$ so that $|y- \gamma_j| \asymp | y - \gamma_k|$
\beq \label{eqn:misc-7}
\int_{[-2,\gamma_{j+\varphi^{10} \ell}]\cup [\gamma_{j-\varphi^{10} \ell},2]}\frac{|\gamma_k-\gamma_j|}{|y-\gamma_k|^2Ns(s+\rho_k\wedge \rho(y))^2}\rho(y) \d y \leq \varphi  \frac{ \kappa_j + \kappa_\ell}{N s (s + \rho_{j \wedge k})^2 \rho_{j \wedge k}} 
\eeq
\eel
\proof We start with the part of the integral with $y \in [-2, \gamma_{\varphi^{10} ( j + k )} ]$. Here $|y - \gamma_k | \asymp |\kappa(y)|$ and so we get,
\begin{align}
& \int_{-2}^{ \gamma_{\varphi^{10} ( j +k ) }} \frac{|\gamma_k-\gamma_j|}{|y-\gamma_k|^2Ns(s+\rho_k\wedge \rho(y))^2}\rho(y) \d y \lesssim \frac{|\gamma_k - \gamma_j| }{Ns (s + \rho_k)^2} \int_{-2}^{ \gamma_{j+k}} \frac{\rhosc(y)}{\kappa(y)^2} \d y \notag\\
\lesssim &  \frac{ | \gamma_j - \gamma_k|}{N s ( s + \rho_k )^2 \rho_{j+k} }
\end{align}
and then we have $| \gamma_j - \gamma_k | \lesssim \kappa_j + \kappa_\ell$ since $|j-k| \leq \ell$. For the part of the integral with $y \in [ \gamma_{ \varphi^{10} ( j +k ) }, \gamma_{j  + \ell} ]$ we use $|y - \gamma_j | \gtrsim | \gamma_k - \gamma_j|$ to get
\beq
\int_{ \gamma_{ \varphi^{10} (j+k)}}^{ \gamma_{j + \varphi \ell}}  \frac{|\gamma_k-\gamma_j|}{|y-\gamma_k|^2Ns(s+\rho_k\wedge \rho(y))^2}\rho(y) \d y \leq \frac{ \rho_j + \rho_k}{N s (s+ \rho_k)^2} \log (N)
\eeq
which is sufficient. For the part of the integral $y \in [ \gamma_{j - \varphi \ell} , 2]$ (which is non-zero only if  $j \geq \varphi^{1/2} \ell$ in which case $j \asymp k$) we can use $|y - \gamma_k | \asymp \kappa_j$ to get,
\beq
\int_{ \gamma_{j- \varphi \ell}}^2 \frac{|\gamma_k-\gamma_j|}{|y-\gamma_k|^2Ns(s+\rho_k\wedge \rho(y))^2}\rho(y) \d y \lesssim \frac{ | \gamma_k - \gamma_j |}{ \kappa_j^2 Ns} \int_{ \gamma_{j- \varphi \ell}}^2 \frac{ \rhosc (y) }{ (s + \rhosc (y))^2} \d y
\eeq
If $s > \rho_j$, then we bound the integral on RHS by $ \frac{ \rho_j^3}{s^2}$ and we're done. Otherwise,
\begin{align}
 & \int_{ \gamma_{j- \varphi \ell}}^2 \frac{ \rhosc (y) }{ (s + \rhosc (y))^2} \d y \leq \int_0^{ \kappa_j} \frac{1}{ y^{1/2}} \leq \rho_j \lesssim \frac{ \rho_j^3}{ ( \rho_j + s)^2}
\end{align}
which is sufficient. \qed

\bel
We have,
\beq \label{eqn:misc-8} 
\sum_{ i : i > j/2, |i-j| \geq \ell} \frac{1}{ \rho_i^b | \gamma_i - \gamma_j |^a} \lesssim N^{b/3+2a/3} \frac{ j^{(a-b)/3} + \ell^{(a-b)/3}}{ \ell^{a-1}}
\eeq
for $b \geq 0 $ and $a > 1.5$
\eel
\proof Using $| \gamma_i - \gamma_j| \gtrsim \frac{ |i-j|}{ N^{2/3} i^{1/3}}$ we have,
\begin{align}
 & \sum_{ i : i > j/2, |i-j| \geq \ell} \frac{1}{ \rho_i^b | \gamma_i - \gamma_j |^a} \lesssim N^{b/3 + 2a/3} \sum_{ i : i > j/2 , |i-j| \geq \ell} \frac{ i^{(a-b)/3}}{|i-j|^a} \notag \\
 \lesssim & N^{b/3 + 2a/3} \sum_{i : i >  j/2, |i-j| \geq \ell} \frac{ j^{(a-b)/3} + |i-j|^{(a-b)/3}}{|i-j|^a} \lesssim N^{b/3+2a/3} \frac{ j^{(a-b)/3} + \ell^{(a-b)/3}}{ \ell^{a-1}}
\end{align}
\qed

\bel
We have, for $1 \geq b \geq 0$,
\beq \label{eqn:misc-9}
\sum_{ i <j/2} \frac{1}{\rho_i^b} \frac{1}{ (s + \rho_i)^2} \lesssim N \frac{ \rho_j^{3-b}}{(s + \rho_j)^2}
\eeq
\eel
\proof If $s > \rho_j$ then
\beq
\sum_{ i <j/2} \frac{1}{\rho_i^b} \frac{1}{ (s + \rho_i)^2} \lesssim \frac{1}{s^2} N^{b/3} \sum_{i < j/2} i^{-b/3} \lesssim \frac{ N^{b/3} j^{1-b/3}}{s^2} \lesssim N \frac{ \rho_j^{3-b}}{(s + \rho_j)^2}
\eeq
Otherwise, we may assume $N s^3 \leq j$ and so 
\begin{align}
 \sum_{ i <j/2} \frac{1}{\rho_i^b} \frac{1}{ (s + \rho_i)^2}  \lesssim N^{b/3+2/3} \sum_{ i < j/2} \frac{1}{i^{b/3+2/3}} \leq N^{b/3+2/3} j^{1/3-b/3} = N \rho_j^{1-b}
\end{align} 
which yields the claim. \qed

\bel
We have for $\eta \geq  (N \rho_j)^{-1}$,
\beq \label{eqn:misc-10}
\frac{1}{N} \sum_{i} \frac{1}{ (\gamma_i - \gamma_j)^2 + \eta^2} \frac{1}{ \rho_i^2} \lesssim \frac{1}{ \eta \sqrt{ \eta+ \kappa_j}}
\eeq
\eel
\proof First consider the sum $i <j/2$. Then,
\begin{align} 
& \frac{1}{N} \sum_{i<j/2} \frac{1}{ (\gamma_i - \gamma_j)^2 + \eta^2} \frac{1}{ \rho_i^2}  \lesssim  \frac{1}{N ( \eta^2 + \kappa_j^2)} \sum_{i < j/2} N^{2/3} i^{-2/3} 
\lesssim  \frac{ \rho_j}{ \eta^2 + \kappa_j^2} 
\end{align}
which is enough.  Consider now when $i / 2j$. If $\kappa_j > \eta$ we bound,
\beq
\frac{1}{N} \sum_{i > 2j} \frac{1}{ ( \gamma_i - \gamma_j)^2 + \eta^2} \frac{1}{ \rho_i^2} \lesssim \frac{1}{N \rho_j^2} \sum_{i > 2j} \frac{N^{4/3}}{i^{4/3}} \lesssim \frac{1}{ \rho_j^3}
\eeq
which is sufficient if $\rho_j > \sqrt{\eta}$. If $\eta > \kappa_j$ then,
\begin{align}
& \frac{1}{N} \sum_{i > 2j} \frac{1}{ ( \gamma_i - \gamma_j)^2 + \eta^2} \frac{1}{ \rho_i^2}  = \frac{1}{N} \sum_{N \eta^{3/2} > i > 2j} \frac{1}{ ( \gamma_i - \gamma_j)^2 + \eta^2} \frac{1}{ \rho_i^2} + \frac{1}{N} \sum_{i > N \eta^{3/2}} \frac{1}{ ( \gamma_i - \gamma_j)^2 + \eta^2} \frac{1}{ \rho_i^2}  \notag\\
\lesssim & \frac{1}{N \eta^2} \sum_{ i < N \eta^{3/2}} \frac{N^{2/3}}{i^{2/3}} + \frac{1}{N} \sum_{i > N \eta}^{3/2} \frac{N^2}{i^2} \lesssim \frac{1}{ \eta^{3/2}} 
\end{align}
which is sufficient. Finally we have the part of the sum where $j/2 < i < 2j$ which is bounded by if $\eta > \kappa_j$,
\begin{align}
\frac{1}{N} \sum_{j/2 < i < 2j} \frac{1}{ ( \gamma_i - \gamma_j)^2 + \eta^2} \frac{1}{ \rho_i^2} \lesssim \frac{1}{N \rho_j^2} \frac{j}{\eta^2} \lesssim \frac{ \rho_j}{\eta^2}
\end{align}
which is sufficient. If $\eta < \kappa_j$ then we can bound it by,
\begin{align}
& \frac{1}{N} \sum_{j/2 < i < 2j} \frac{1}{ ( \gamma_i - \gamma_j)^2 + \eta^2} \frac{1}{ \rho_i^2} \lesssim \frac{1}{N \rho_j^2} \sum_{|i-j| < j } \frac{(N \rho_j)^2}{ |i-j|^2 + ( \eta N \rho_j)^2} \notag\\
\lesssim & N \int_\rr \frac{1}{ x^2 + (N \eta \rho_j)^2} \lesssim \frac{1}{ \eta \rho_j}
\end{align}
which is sufficient.  \qed

\bel
If $\eta \geq (N \rho_j)^{-1}$ we have
\beq \label{eqn:misc-11}
\int_{-2}^2 \frac{\rho(y)}{\pare{|y-\gamma_j|^2+\eta^2}N(t+\rho(y))\pare{Nt(t+\rho(y))\rho(y)}^{1-a}} \d y \lesssim \frac{ \sqrt{ \kappa_j+\eta}}{\eta} \frac{1}{N(t+\rho_j) [ Nt \rho_j (t + \rho_j)]^{1-a}} 
\eeq
\eel
\proof The part of the integral with $y < \gamma_{j/2}$ is easy as we bound it by, (as $0 \leq a \leq 1$)
\begin{align}
& \int_{-2}^{\gamma_{j/2}}  \frac{\rho(y)}{\pare{|y-\gamma_j|^2+\eta^2}N(t+\rho(y))\pare{Nt(t+\rho(y))\rho(y)}^{1-a}} \d y \notag\\
\lesssim & \frac{1}{N(t+\rho_j)( Nt \rho_j (t + \rho_j) )^{1-a} } \int_{-2}^2 \frac{ \rhosc (y)}{ (y-\gamma_j)^2 + \eta^2} 
\end{align}
and the last integral is $\frac{ \sqrt{ \kappa_j + \eta}}{\eta}$ as its $ \Im[ \msc]$. For the part of the integral with $y >  \gamma_{j/2}$ we have,
\begin{align}
 & \int_{\gamma_{j/2}}^{2}  \frac{\rho(y)}{\pare{|y-\gamma_j|^2+\eta^2}N(t+\rho(y))\pare{Nt(t+\rho(y))\rho(y)}^{1-a}} \d y  \notag\\
\lesssim & \frac{1}{ \kappa_j^2 + \eta^2} \int_{\gamma_{j/2}}^2 \frac{ \rhosc(y)}{ N (t + \rhosc(y) ) ( N t \rhosc (y) (t+ \rhosc(y))^{1-a}}  \d y \lesssim \frac{1}{ \kappa_j^2 + \eta^2} \frac{1}{(Nt)^{1-a} N} \frac{ \rho_j^{2+a}}{(s+\rho_j)^{2-a}}
\end{align}
where we used \eqref{eqn:misc-5a} in the second inequality. This is sufficient. \qed

\bel Let $\eta > (N \rho_j)^{-1}$. 
We have,
\beq \label{eqn:misc-12} 
\sum_{i \neq j :|\gamma_i-\gamma_j|\leq\eta}\frac{1}{|\gamma_i-\gamma_j|Nt(t+\rho_i\wedge\rho_j)^2} \lesssim \frac{ N \sqrt{ \kappa_j + \eta}}{N t (t + \rho_j)^2}
\eeq
\eel
\proof When $i >j/2$ we simply bound,
\begin{align}
 & \sum_{i \neq j :|\gamma_i-\gamma_j|\leq\eta, i\geq j/2 }\frac{1}{|\gamma_i-\gamma_j|Nt(t+\rho_i\wedge\rho_j)^2} \leq \frac{1}{N t (t + \rho_j)^2} \sum_{j/2 \leq i \leq j + N \eta^{3/2}} \frac{ N^{2/3} i^{1/3}}{ |i-j|} \notag\\
\leq & \log N \frac{ N \rho_j + N \sqrt{\eta}}{N t (t + \rho_j)^2}
\end{align}
where we just bounded $i \leq j + N \eta^{3/2}$ for the numerator. THen for $i < j/2$ just use $| \gamma_i - \gamma_j| \geq c \kappa_j$ and \eqref{eqn:misc-9}. \qed

\bel
Assume $\eta > (N \rho_j)^{-1}$. We have,
\begin{align} \label{eqn:misc-13}
\int \frac{\rho(y)}{|y-\gamma_j|^3+\eta^3}\frac{1}{N(t+\rho(y))\pare{Nt(t+\rho(y))\rho(y)}^{1-a}} \lesssim \frac{ \sqrt{ \kappa_j + \eta}}{\eta^2} \frac{1}{Nt (t+\rho_j) (N t \rho_j (t+\rho_j) )^{1-a}} 
\end{align}
\eel
\proof For $y > \gamma_{j/2}$ bound $|y-\gamma_j| + \eta \geq \kappa_j + \eta$ and use \eqref{eqn:misc-5a}. For $y < \gamma_{2j}$ bound it as,
\begin{align}
& \int_{-2}^{\gamma_{2j}} \frac{\rho(y)}{|y-\gamma_j|^3+\eta^3}\frac{1}{N(t+\rho(y))\pare{Nt(t+\rho(y))\rho(y)}^{1-a}} \notag\\
\lesssim & \frac{1}{N t (t + \rho_j)( N t \rho_j (t + \rho_j))^{1-a} } \frac{1}{\eta} \int \frac{ \rhosc (y)}{ (y-\gamma_j)^2 + \eta^2}
\end{align}
and use asymptotics for $\Im[ \msc]$. For the rest,
\begin{align}
 & \int_{\gamma_{2j}}^{\gamma_{j/2}} \frac{\rho(y)}{|y-\gamma_j|^3+\eta^3}\frac{1}{N(t+\rho(y))\pare{Nt(t+\rho(y))\rho(y)}^{1-a}} \notag\\
\leq & \frac{\rho_j}{N(t + \rho_j) ( Nt \rho_j (t + \rho_j ) )^{1-a}} \eta^{-1} \int_{\rr} \frac{1}{ (y- \gamma_j)^2 + \eta^2} 
\end{align}
which is enough. \qed 

\bel
We have,
\begin{align} \label{eqn:misc-14}
 & \int_{\gamma_j}^2 \rhosc(y) \left( \frac{1}{ N (t + \rhosc (y)) ( (N t (t + \rhosc(y) ) \rhosc(y) )^{1-a}} \right)^2 \notag\\
 \lesssim & \1_{ t < \rho_j }   t^3 \left( \frac{1}{ N t (N t^3)^{1-a} } \right)^2  +   \left( \frac{1}{N(t + \rho_j) ( Nt \rho_j (t+\rho_j ) )^{1-a}} \right)^2 \rho_j^3
\end{align}
\eel
\proof If $t > \rho_j$ we can bound it by
\beq
\left( \frac{1}{ Nt (N t^2)^{1-a} } \right)^2 \int_{0}^{\kappa_j} y^{a-1/2} \lesssim \left( \frac{\rho_j^{a-1}}{ Nt (N t^2)^{1-a} } \right)^2  \rho_j^3 \lesssim\left( \frac{1}{N(t+\rho_j)( (Nt\rho_j (t + \rho_j) )^{1-a}} \right)^2 \rho_j^3
\eeq
If $ t < \rho_j$ we first bound the integral $y \in [2, 2-t^2]$ by,
\beq
\left( \frac{1}{ Nt (N t^2)^{1-a} } \right)^2 \int_{0}^{t^2} y^{a-1/2} \lesssim \left( \frac{1}{ Nt (N t^2)^{1-a} } \right)^2 t^{2a+1} \lesssim t^3 \left( \frac{1}{ N t (N t^3)^{1-a} } \right)^2 
\eeq
and the integral $y \in [2-t^2, \gamma_j]$ by,
\begin{align}
\left( \frac{1}{N (N t)^{1-a} } \right)^2 \int_{t^2}^{ \kappa_j} y^{2a-2.5} 
\end{align}
If $a \geq 0.75$ this gives
\beq
\left( \frac{1}{N (N t)^{1-a} } \right)^2  \rho_j^{4a-3} \lesssim \left( \frac{1}{N(t + \rho_j) ( Nt \rho_j (t+\rho_j ) )^{1-a}} \right)^2 \rho_j^3
\eeq
and if $a \leq 0.75$ this gives
\beq
\left( \frac{1}{N (N t)^{1-a} } \right)^2   t^{4a-3} \lesssim t^3 \left( \frac{1}{ N t (N t^3)^{1-a} } \right)^2 
\eeq
\qed

\bel
We have,
\begin{align} \label{eqn:misc-15}
 & \int_{\gamma_j}^2 \frac{\rho(y)}{(Nt(t+\rho(y))^2)^2}dy 
\lesssim \1_{t < \rho_j} \frac{t^3}{ (N t^3)^2} + \frac{ \rho_j^3}{ ( N t (t + \rho_j)^2 )^2}
\end{align}
\eel
\proof If $t > \rho_j$ we easily bound it by
\beq
\frac{ \rho_j^3}{ ( N t (t + \rho_j)^2 )^2}
\eeq
If $ t < \rho_j$ we bound the integral $y \in [2, 2-t^2]$ by 
\beq
\frac{ t^3}{(N t^3)^2}
\eeq
and the integral for $y \in [2-t^2, \gamma_j]$ by
\beq
\frac{1}{ (N t)^2 } \int_{t^2}^{ \kappa_j} y^{-1.5} = \frac{t^3}{ (N t^3)^2}
\eeq
\qed

\section{Time integrals for flows}

For the proofs in this section, recall that 
\beq
\eta_u = u^2 + u \sqrt{\kappa+\eta} + \eta. 
\eeq
In particular, if $u < \frac{\eta}{\sqrt{\kappa_j+\eta}}$ then $\eta_u \asymp \eta$. If $\frac{\eta}{\sqrt{\eta+\kappa_j}} \leq u \leq \sqrt{\kappa_j+\eta}$ then $\eta_u \asymp u \sqrt{\kappa_j+\eta}$. If $u \geq \sqrt{\kappa_j+\eta}$ then $\eta_u \asymp u^2$. All of the estimates in this section can be deduced by breaking up the time integral into these three regimes and estimating the result. We give the details for the first few estimates, the others are the similar.

\bel
We have,
\beq \label{eqn:misc-time-1}
\int_0^{t/2}\frac{du}{\eta_u\pare{|\gamma_i-\gamma_j|^2+\eta_u^2}} \leq \varphi  \frac{1}{\sqrt{\kappa_j+\eta}(|\gamma_i-\gamma_j|^2+\eta^2)}
\eeq
\eel
\proof 
\begin{align*}
&\int_0^{t/2}\frac{du}{\eta_u\pare{|\gamma_i-\gamma_j|^2+\eta_u^2}}\\
&\lesssim \int_{\sqrt{\kappa_j+\eta}}^t\frac{du}{u^2(|\gamma_i-\gamma_j|^2+u^4)} + \int_{\frac{\eta}{\sqrt{\kappa_j+\eta}}}^{\sqrt{\kappa_j+\eta}}\frac{du}{u\sqrt{\kappa_j+\eta}(|\gamma_i-\gamma_j|^2+(u\sqrt{\kappa_j+\eta})^2))}\\
&\quad+\int_0^{\frac{\eta}{\sqrt{\kappa_j+\eta}}}\frac{du}{\eta(|\gamma_i-\gamma_j|^2+\eta^2)}\\
&\lesssim \frac{1}{\sqrt{\kappa_j+\eta}(|\gamma_i-\gamma_j|^2+(\kappa_j+\eta)^2)} + \frac{\log\frac{\kappa_j+\eta}{\eta}}{\sqrt{\kappa_j+\eta}(|\gamma_i-\gamma_j|^2+\eta^2)} +\frac{1}{\sqrt{\kappa_j+\eta}(|\gamma_i-\gamma_j|^2+\eta^2)}\\
&\leq \varphi \frac{1}{\sqrt{\kappa_j+\eta}(|\gamma_i-\gamma_j|^2+\eta^2)}.
\end{align*}
\qed

\bel
We have,
\beq \label{eqn:misc-time-2}
\int_0^{\frac{t}{2}}\frac{\sqrt{\kappa_j+\eta_u}}{\eta_u} \d u \leq \varphi
\eeq
\eel
\proof \begin{align*}
&\int_0^{\frac{t}{2}}\frac{\sqrt{\kappa_j+\eta_u}}{\eta_u}du\\
&\lesssim \int_{\sqrt{\kappa_j+\eta}}^t\frac{1}{u}+\int_{\frac{\eta}{\sqrt{\kappa_j+\eta}}}^{\sqrt{\kappa_j+\eta}} \frac{\sqrt{\kappa_j+u\sqrt{\kappa_j+\eta}}}{u\sqrt{\kappa_j+\eta}} + \int_0^{\frac{\eta}{\sqrt{\kappa_j+\eta}}}\frac{\sqrt{\kappa_j+\eta}}{\eta}\\
&\leq \varphi
\end{align*}

To bound the second term in the second line, we note that $\sqrt{\kappa_j+u\sqrt{\kappa_j+\eta}}\asymp \sqrt{\kappa_j+\eta}$. 
\qed

\bel
We have,
\beq
\label{eqn:misc-time-3}
\int_0^{t/2}\frac{\d u}{\eta_u} \leq \varphi \frac{1}{\sqrt{\kappa_j+\eta}},
\eeq
\eel
\proof We have, 
\begin{align} \label{eqn:time-int-a1}
\int_0^{t/2}\frac{\d u}{\eta_u} \lesssim \int_{\sqrt{\kappa_j+\eta}}^t\frac{\d u}{u^2} + \int_{\frac{\eta}{\sqrt{\kappa_j+\eta}}}^{\sqrt{\kappa_j+\eta}}\frac{\d u}{u\sqrt{\kappa_j+\eta}}+\int_0^{\frac{\eta}{\sqrt{\kappa_j+\eta}}}\frac{\d u}{\eta} \lesssim \varphi \frac{1}{\sqrt{\kappa_j+\eta}},
\end{align}
\qed

The proofs of the estimates below are very similar to the above proofs and so we omit them. 

\bel
We have,
\beq \label{eqn:misc-time-4}
\int_0^t \frac{\d u}{\sqrt{\kappa_j+\eta_u}}\leq \varphi \left(  1\wedge \frac{t}{\rho_j} \right) 
\eeq
and
\beq \label{eqn:misc-time-5}
\int_0^t \frac{\sqrt{\kappa_j + \eta_u} \d u}{\eta_u^2} \lesssim \frac{1}{ \eta} 
\eeq
\eel

\bel
We have,
\beq \label{eqn:misc-time-6}
\int_0^t \frac{\sqrt{ \kappa_j + \eta_u}}{\eta_u^3} \d u \lesssim \frac{1}{\eta^2}
\eeq
and
\beq \label{eqn:misc-time-7}
\int_0^t \frac{ \d u}{ \kappa_j^4 + \eta_u^4} \lesssim \frac{t}{\kappa_j^4}
\eeq
and
\beq \label{eqn:misc-time-8}
\int_0^t \frac{ \d u}{ \sqrt{ \kappa_j+\eta_u}} \lesssim  1
\eeq
and
\beq \label{eqn:misc-time-9}
\int_0^t \frac{1}{ \eta_u ( \kappa_j^3 + \eta_u^3)} \d u \lesssim  \frac{1}{ ( \kappa_j + \eta)^{7/2}}
\eeq
\eel

\vspace{5 pt}

\noindent{\bf Acknowledgements.} The work of B.L. and T.X. is supported by an NSERC Discovery Grant and a Connaught New Researcher award. B.L. thanks Paul Bourgade for helpful discussions about this project.


\bibliography{mybib}{}
\bibliographystyle{abbrv}

\end{document}